\documentclass[a4paper,10pt,oneside]{article} 
\usepackage{geometry}
\usepackage[T1] {fontenc}
\usepackage [utf8] {inputenc}
\usepackage{longtable}
\usepackage[english] {babel}
\usepackage{graphicx}
\usepackage{authblk}
\usepackage{bbm}
\usepackage{amsmath,amssymb,amsthm, xcolor, listings}
\usepackage{mathrsfs}
\usepackage{tikz}
\usetikzlibrary{patterns}
\usetikzlibrary{decorations.pathreplacing}
\usepackage{paralist}
\usepackage{fancyhdr}
\usepackage{mathtools}
\usepackage{titlesec}
\usepackage{enumitem}
\usepackage{float,placeins}
\usepackage{url}
\usepackage{hyperref}
\usepackage{bm}
\usepackage{mathdots}
\usepackage{eurosym}

\newtheorem{theorem}{Theorem}[section]
\newtheorem{proposition}{Proposition}[section]
\newtheorem{lemma}{Lemma}[section]
\newtheorem{cor}{Corollary}[section]

\theoremstyle{definition}
\newtheorem{definition}{Definition}[section]
\newtheorem{remark}{Remark}[section]

\numberwithin{equation}{section}
\begin{document}
\title{Critical configurations and tube of typical trajectories \\for the Potts and Ising models with zero external field}

\author[
         {}\hspace{0.5pt}\protect\hyperlink{hyp:email1}{1},\protect\hyperlink{hyp:affil1}{a}
        ]
        {\protect\hypertarget{hyp:author1}{Gianmarco Bet}}

\author[
         {}\hspace{0.5pt}\protect\hyperlink{hyp:email2}{2},\protect\hyperlink{hyp:affil1}{a},\protect\hyperlink{hyp:corresponding}{$\dagger$}
        ]
        {\protect\hypertarget{hyp:author2}{Anna Gallo}}

\author[
         {}\hspace{0.5pt}\protect\hyperlink{hyp:email3}{3},\protect\hyperlink{hyp:affil1}{a},\protect\hyperlink{hyp:affil2}{b}
        ]
        {\protect\hypertarget{hyp:author3}{Francesca R.~Nardi}}

\affil[ ]{
          \small\parbox{365pt}{
             \parbox{5pt}{\textsuperscript{\protect\hypertarget{hyp:affil2}{a}}}Università degli Studi di Firenze,
            \enspace
             \parbox{5pt}{\textsuperscript{\protect\hypertarget{hyp:affil1}{b}}}Eindhoven University of Technology
            }
          }

\affil[ ]{
          \small\parbox{365pt}{
             \parbox{5pt}{\textsuperscript{\protect\hypertarget{hyp:email1}{1}}}\texttt{\footnotesize\href{mailto:gianmarco.bet@unifi.it}{gianmarco.bet@unifi.it}},
             \parbox{5pt}{\textsuperscript{\protect\hypertarget{hyp:email2}{2}}}\texttt{\footnotesize\href{mailto:anna.gallo1@stud.unifi.it}{anna.gallo1@stud.unifi.it}},
             \parbox{5pt}{\textsuperscript{\protect\hypertarget{hyp:email3}{3}}}\texttt{\footnotesize\href{mailto:francescaromana.nardi@unifi.it}{francescaromana.nardi@unifi.it}}
            }
          }

\affil[ ]{
          \small\parbox{365pt}{
             \parbox{5pt}{\textsuperscript{\protect\hypertarget{hyp:corresponding}{$\dagger$}}}Corresponding author
            }
          }

\date{\today}

\maketitle

\begin{abstract}
We consider the ferromagnetic $q$-state Potts model with zero external field in a finite volume evolving according to Glauber-type dynamics described by the Metropolis algorithm in the low temperature asymptotic limit. Our analysis concerns
the multi-spin system that has $q$ stable equilibria. Focusing on grid graphs with periodic boundary conditions, we study the tunneling 
between two stable states and from one stable state to the set of all other stable states.
In both cases we identify the set of gates for the transition and prove that this set has to be crossed with high probability during the transition. Moreover, we identify the tube of typical paths and prove that the probability to deviate from it during the transition is exponentially small.

\\
 \emph{Acknowledgment:} The research of Francesca R.~Nardi was partially supported by the NWO Gravitation Grant 024.002.003--NETWORKS and by the PRIN Grant 20155PAWZB ``Large Scale Random Structures''.

\end{abstract}

\maketitle

\tableofcontents
\section{Introduction}\label{intro}
Metastability is a phenomenon that occurs when a physical system is close to a first order phase transition. More precisely, the phenomenon of metastability occurs when a system is trapped for a long time in a state different from the stable state, the so-called \textit{metastable state}.
After a long (random) time or due to random fluctuations the system makes a sudden transition from the metastable state to the stable state. When this happens, the system is said to display \textit{metastable behavior}. 
On the other hand, when the system lies on the phase coexistence line, it is of interest to investigate its \textit{tunneling behavior}, i.e., how the system transitions between the two (or more) stable states.
Since metastability occurs in several physical situations, such as supercooled liquids and supersaturated gases, many models for metastable behavior have been formulated throughout the years. Tipically, the evolution of the physical system is approximated by a stochastic process, and broadly speaking the following three main issues are investigated. The first is the study of the \textit{first hitting time} at which the process starting from a metastable state visits a stable state. The second issue is the study of the so-called set of \textit{critical configurations}, i.e., the set of those configurations that are crossed by
the process during the transition from the metastable state to the stable state. The final issue is the study of the typical trajectories that the system follows during the transition from the metastable state to the stable state. This is the so-called \textit{tube of typical paths}. The same three issues are investigated when a system displays tunneling behavior, except that in this
case one is interested in the transition from one stable state to another stable state.\\

In this paper study the tunneling behavior of the $q$-state Potts model on a two-dimensional discrete torus.
At each site $i$ of the lattice lies a spin with value $\sigma(i)\in\{1,\dots,q\}$. To each spin configuration we associate an energy such that configura-tions where neighboring spins have the same value are energetically favored. A model that
satisfies this condition is said to be ferromagnetic. 
The $q$-state Potts model is an extension of the classical Ising model from $q=2$ to an arbitrary number of spins states.
We study the $q$-state ferromagnetic Potts model with zero external magnetic field ($h=0$) in the limit of large inverse temperature $\beta\to\infty$. When the external magnetic field is zero, the system lies on a coexistence line. 
The stochastic evolution is described by a \textit{Glauber-type dynamics}, that consists of a single-spin flip Markov chain on a finite state space $\mathcal X$ with transition probabilities given by the Metropolis algorithm and with stationary distribution given by the \textit{Gibbs measure} $\mu_\beta$, see \eqref{gibbs}. 
We consider the setting where there is no external magnetic field, and so to each configuration $\sigma\in\mathcal X$ we associate an energy $H(\sigma)$ that only depends on the local interactions between nearest-neighbor spins. In the low-temperature regime $\beta\gg1$ there are $q$ stable states, corresponding to the configurations where all spins are equal. In this setting, the metastable states are not interesting since they do not have a clear physical interpretation, hence we focus our attention on the tunneling behavior between stable configurations. \\

The goal of this paper is to investigate the second and third issues introduced above for the tunneling behavior of the system. We describe the set of \textit{minimal gates}, which have the physical meaning of  ``critical configurations'', and the \textit{tube of typical paths} for three different types of transitions. More precisely, we study the transition from any stable configuration $\mathbf r$
(a) to some other stable configuration $\mathbf s\neq\mathbf r$ under the constraint that the path followed does not intersect the cycle of other stable configurations, (b) to \textit{any} other stable configuration, and (c) to some other stable configuration $\mathbf s\neq\mathbf r$. In Section \ref{modinddefnotgates} (resp.~Section \ref{deftubemodind}) we introduce the notion of \textit{minimal restricted-gates} (resp.~\textit{restricted-tube of typical paths}) to denote the minimal gates (resp.~tube of typical paths) for the transition (a). 
\begin{figure}
\centering
\begin{tikzpicture} [scale=0.7, transform shape]
\draw [fill=gray,lightgray] (0,0) rectangle (1.5,2.7);
\draw [fill=gray,lightgray] (1.2,0.3) rectangle (1.8,1.2);
\draw[step=0.3cm,color=black] (0,0) grid (3.6,2.7);
\end{tikzpicture}\ \ \ \ \ \ \ \
\begin{tikzpicture}[scale=0.7, transform shape]
\draw [fill=gray,lightgray] (0,0) rectangle (2.4,2.7);
\draw [fill=gray,lightgray] (2.4,1.8) rectangle (2.7,1.2);
\draw[step=0.3cm,color=black] (0,0) grid (3.6,2.7);
\end{tikzpicture}
\caption{\label{figintrogates} Example of configurations belonging to the set of the minimal-restricted gates between the stable configurations $\mathbf r$ and $\mathbf s$. We color white the vertices with spin $r$, gray those vertices with spin $s$.}
\end{figure}
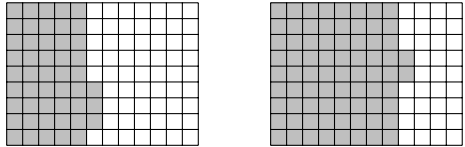
Let us now briefly describe our approach. First we  study the energy landscape between two stable configurations. Roughly speaking, we prove that the set of minimal-restricted gates for any transition (a) contains those configurations in which all the spins are $r$ (respectively $s$) except those, which have spins $s$ (respectively $r$), in a strip that wraps around the shortest side of the torus and that also has a bar attached to one of the two vertical sides, see Figure \ref{figintrogates}. We build on this result to describe the set of minimal gates for the transitions (b) and (c). Next we describe the tube of typical trajectories for the transitions (a), (b) and (c). Once again we first describe the restricted-tube of typical paths between two stable configurations, and then we lean on this result to describe the tube of typical paths for the transitions (b) and (c). We show that the restricted-tube of typical paths between the stable states $\mathbf r$ and $\mathbf s$ includes the minimal-restricted gates, as well as configurations with one or more clusters of spin $s$ (respectively $r$) --of at most a certain size, which we identify-- surrounded by spins $r$ (respectively $s$). See Subsection \ref{moddeptubo} and Figures \ref{examplesprincbound}--\ref{thirdexamplesprincbound} below. As a special case of our general results we retrieve the minimal gates and the tube of typical paths for the Ising model with zero external magnetic field. We give these respectively at the end of Section \ref{secgatesmoddef} and of Section \ref{tube}.\\
\textit{Related work.} Our work concludes the study of the metastability of the Potts model in the low-temperature regime first initiated in \cite{nardi2019tunneling}, where the authors derive the asymptotic behavior of the \textit{first hitting time} associated with the transitions (b) and (c) above. They obtain convergence results in probability, in expectation and in distribution. They also investigate the \textit{mixing time}, which describes the rate of convergence of the process to its stationary distribution $\mu_\beta$. They further show that, as $\beta\to\infty$, the mixing time grows as $\text{exp}{(c\ell\beta)}$, where $c>0$ is some constant constant factor and $\ell$ is the smallest side length of $\Lambda$.

\noindent In \cite{kim2021metastability} the authors study the $q$-Potts model with zero external magnetic field in two and three dimensions. Their manuscript appeared on ArXiv roughly at the same time as ours. While there are some overlaps between the two papers, the works of the two groups were carried out independently of each other. They find sharp estimates for the tunneling time (the so-called prefactor) using the potential-theoretic approach and extending the results in \cite{nardi2019tunneling}. In the second part of their manuscript, the authors focus on the two-dimensional setting. However, the results and the strategy adopted to prove them are quite different from ours. For this type of transitions, they define the set of the \textit{gateway configurations} in Definition 8.1 that is a model-dependent set visited with probability tending to one as $\beta\to\infty$, see \cite[Corollary 8.9]{kim2021metastability}. 
This set is different from the classical model-independent definition of gate and of union of minimal gates, see \cite{manzo2004essential} and following papers. The gates and minimal gates have the physical meaning of critical configurations. In Theorems \ref{mingatescond}, \ref{mingatesNOcond} and \ref{mingatessingh0}, we give an explicit and thorough geometric description of the minimal gates and union of all minimal gates for the three transitions (a), (b), (c). Moreover, in Corollaries \ref{corgatecond}, \ref{corgateNOcond}, \ref{corgateconv2} we prove that these sets are crossed with probability tending to one. On the other hand, \cite{kim2021metastability} does not give a complete geometrical characterization of the set of gateway configurations and this set is not a minimal gate. The description of the gateway configurations is suitable to allow them to compute the prefactor. Finally, in our paper we analyse the third issue of metastability by indentifying precisely the tube of typical trajectories. This analysis is absent in \cite{kim2021metastability}.\\
In \cite{bet2021metastability}, the authors consider the q-Potts model with non-zero external field and analyze separately the case of positive and negative external magnetic field. In the first scenario there are $q-1$ stable configurations and a unique metastable state. In the second scenario there are $q-1$ degenerate-metastable configurations and only one global minimum. In both cases the authors describe the asymptotic behavior of the first hitting time from the metastable to the stable state as $\beta\to\infty$, the mixing time, the spectral gap and they identify geometrically the set of gates for these transitions. 

We adopt the statistical mechanics framework known as  \textit{pathwise approach}. This is a set of techniques that rely on a detailed knowledge of the energy landscape and on ad hoc large deviations estimates to give a quantitative answer to the three issues of metastability which we described above. The pathwise approach was first introduced in 1984 \cite{cassandro1984metastable} and then developed in \cite{olivieri1995markov},\cite{olivieri1996markov},\cite{olivieri2005large}, \cite{catoni1997exit}. We adopt the convention of listing citations in order of publication date. This approach was further expanded in \cite{manzo2004essential},\cite{cirillo2013relaxation},\cite{cirillo2015metastability},\cite{nardi2016hitting},\cite{fernandez2015asymptotically},\cite{fernandez2016conditioned} to distinguish the study of the transition time and of the gates from the study of typical paths. In \cite{olivieri2005large},\cite{manzo2004essential} the pathwise approach was expanded and refined with the aim of finding answers valid with maximal generality and of reducing as much as possible the number of model dependent inputs necessary to study the metastable behaviour of a system.

The pathwise approach was applied in \cite{arous1996metastability},\cite{cirillo1998metastability},\cite{cirillo1996metastability},\cite{kotecky1994shapes},\cite{nardi1996low},\cite{neves1991critical},\cite{neves1992behavior},\cite{olivieri2005large} to study the metastable behaviour of Ising-like models with Glauber dynamics. Moreover, the approach was used in \cite{hollander2000metastability},\cite{den2003droplet},\cite{gaudilliere2005nucleation},\cite{apollonio2021metastability},\cite{nardi2016hitting},\cite{zocca2019tunneling} to find the transition time and the gates for Ising-like and hard-core models with Kawasaki and Glauber dynamics. The pathwise approach was also applied to study probabilistic cellular automata (parallel dynamics) in \cite{cirillo2003metastability},\cite{cirillo2008competitive},\cite{cirillo2008metastability},\cite{procacci2016probabilistic},\cite{dai2015fast}. 

On the other hand, the so-called \emph{potential-theoretical approach} exploits a suitable Dirichlet form and spectral properties of the transition matrix to study the hitting time of metastable dynamics. One of the advantages of this approach is that it makes possible the estimation of the expected value of the transition time up to the (lower-order)
coefficient that multiplies the (leading-order) exponential term. The coefficient is known in the literature as the \textit{pre-factor}. These results are grounded in a detailed knowledge of the critical configurations and on the configurations connected to them in one step of the the dynamics, see \cite{bovier2002metastability},\cite{bovier2004metastability},\cite{bovier2016metastability},\cite{cirillo2017sum}. This method was applied to find the pre-factor for Ising-like models and the hard-core model in \cite{bashiri2019on},\cite{boviermanzo2002metastability},\cite{cirillo2017sum},\cite{bovier2006sharp},\cite{den2012metastability},\cite{jovanovski2017metastability},\cite{den2018metastability} for Glauber and Kawasaki dynamics and in \cite{nardi2012sharp},\cite{bet2020effect} for parallel dynamics.
Recently, other approaches have been developed in \cite{beltran2010tunneling},\cite{beltran2012tunneling},\cite{gaudillierelandim2014} and in \cite{bianchi2016metastable}. These approaches are particularly well-suited to find the pre-factor when dealing with the tunnelling between two or more stable states.

The outline of the paper is as follows. At the beginning of Section \ref{modeldescription}, we define the model. In Section \ref{secgatesmoddef} we give a list of definitions that are necessary to state our main results on the set of minimal restricted-gates and on the set of minimal gates. In Section \ref{subgateres} we give the main results for the minimal restricted-gates for the transition (a). In Section \ref{subgateset} and Section \ref{subgatetarget} we state our main results for the minimal gates for the transitions (b) and (c), respectively. Next, at the beginning of Section \ref{tube} we expand the list of definitions in order to state the main results on the restricted-tube and on the tube of typical paths. More precisely, in Section \ref{mainrestuberes} we state the main results on the restricted-tube of typical paths. See Sections \ref{mainrestubeset} and \ref{mainrestubetarget} for the main results on the tube of typical paths for the transitions (b) and (c), respectively. In Section \ref{minimalresgatessection} we prove some useful lemmas that allow us to complete the proof of the main results stated in Section \ref{subgateres}. In Section \ref{proofmaingates} we carry out the proof of the main results introduced in Section \ref{subgateset} and Section \ref{subgatetarget}. Finally, in Section \ref{proofmaintube} we describe the typical paths between two Potts stable states and we prove the main results given in Section \ref{mainrestube}. In the Appendix \ref{sectionappendix} we give some more detailed proofs.

\section{Model description}\label{modeldescription}
In the $q$-state Potts model each spin lies on the vertices of a finite two-dimensional rectangular lattice $\Lambda=(V,E)$, where $V=\{0,\dots,K-1\}\times\{0,\dots,L-1\}$ is the vertex set and $E$ is the edge set, namely the set of the pairs of vertices whose spins interact with each other. We consider periodic boundary conditions, that is, we identify
each pair of vertices lying on opposite sides of the rectangle, so that we end up with a two-dimensional torus. Using this representation, two vertices $v,w\in V$ are said to be
nearest-neighbors when they share an edge of $\Lambda$. Without loss of generality, we assume
$K<L$ and $L \ge 3$.
Let $S= \{1,\dots,q\}$ be the set of spin values. To each vertex $v\in V$ is associated a spin value $\sigma(v)\in S$, and $\mathcal X := S^V$ denotes the set of spin
configurations.  

To each configuration $\sigma \in \mathcal{X}$ we associate the energy $H(\sigma)$ given by  
\begin{align}\label{hamiltoniangeneral}
H(\sigma):=-J_c \sum_{(v,w)\in E} \mathbbm{1}_{\{\sigma(v)=\sigma(w)\}}, \ \ \sigma\in\mathcal{X},
\end{align}
where $J_c$ is the \textit{coupling} or \textit{interation constant}. The function $H:\mathcal X\to \mathbb R$ is called \textit{Hamiltonian} or \textit{energy function} $H: \mathcal{X} \to \mathbb{R}$. In particular, there is no external magnetic field and $H$ is just a sum of the local interactions between nearest-neighbor spins. The
Potts model is said to be ferromagnetic when $J_c>0$, and antiferromagnetic otherwise. In this paper we focus on the ferromagnetic Potts model, and we set $J_c=1$ without loss of generality. We denote by $\mathcal{X}^s$ the set of the global minima of the Hamiltonian \eqref{hamiltoniangeneral} in $\mathcal X$.  \\

The \textit{Gibbs measure} for the $q$-state Potts model on $\Lambda$ is a probability distribution on $\cal X$ given by
\begin{align}\label{gibbs}
\mu_\beta(\sigma):=\frac{e^{-\beta H(\sigma)}}{\sum_{\sigma'\in\mathcal{X}}e^{-\beta H(\sigma')}},
\end{align}
where $\beta>0$ is the inverse temperature. In particular, when $J_c>0$ in the low-temperature regime $\beta\gg1$, $\mu_\beta$ is concentrated on the the global minima of $H$. 
In our setting,
by simple algebraic calculations we prove that these are the configurations with constant spin values, and we denote them by $\mathbf{1},\dots,\mathbf{q}\in\mathcal X$. For example, $\mathbf1(v) = 1$ for any $v\in V$. 
\begin{proposition}[Identification of the stable configurations]\label{propsetstabzero}
The set of the global minima of the Hamiltonian \eqref{hamiltoniangeneral} is
\begin{align}
\mathcal X^s=\{\mathbf 1,\dots,\mathbf q\}.
\end{align}
\end{proposition}

The spin system evolves according to a Glauber-type dynamics. This is described by a single-spin update Markov chain $\{X_t^\beta\}_{t\in\mathbb{N}}$ on the state space $\mathcal{X}$ with the following transition probabilities: for $\sigma, \sigma' \in \mathcal{X}$,
\begin{align}\label{metropolisTP}
P_\beta(\sigma,\sigma'):=
\begin{cases}
Q(\sigma,\sigma')e^{-\beta [H(\sigma')-H(\sigma)]^+}, &\text{if}\ \sigma \neq \sigma',\\
1-\sum_{\eta \neq \sigma} P_\beta (\sigma, \eta), &\text{if}\ \sigma=\sigma',
\end{cases}
\end{align}
where $[n]^+:=\max\{0,n\}$ is the positive part of $n$ and $Q$ is the \textit{connectivity matrix} defined by
\begin{align}\label{Qmatrix}
Q(\sigma,\sigma'):=
\begin{cases}
\frac{1}{q|V|}, &\text{if}\ |\{v\in V: \sigma(v) \neq \sigma'(v)\}|=1,\\
0, &\text{if}\ |\{v\in V: \sigma(v) \neq \sigma'(v)\}|>1,
\end{cases}
\end{align}
for any $\sigma, \sigma' \in \mathcal{X}$.
The matrix $Q$ is symmetric and irreducible, i.e., for all $\sigma, \sigma' \in \mathcal{X}$, there exists a finite sequence of configurations $\omega_1,\dots,\omega_n \in \mathcal{X}$ such that $\omega_1=\sigma$, $\omega_n=\sigma'$ and $Q(\omega_i,\omega_{i+1})>0$ for $i=1,\dots,n-1$.  Hence, the resulting stochastic dynamics defined by \eqref{metropolisTP} is reversible with respect to the Gibbs measure \eqref{gibbs}. We shall refer to the triplet $(\mathcal{X},H,Q)$ as the \textit{energy landscape}. \\

The dynamics defined above belongs to the class of Metropolis dynamics. More precisely, given a configuration $\sigma \in \mathcal{X}$, at each step
\begin{itemize}
\item[1.] a vertex $v \in V$ and a spin value $s \in S$ are selected independently and uniformly at random;
\item[2.] the spin at $v$ is updated to spin $s$ with probability
\begin{align}\begin{cases}
1, &\text{if}\ H(\sigma^{v,s})-H(\sigma) \le 0,\\
e^{-\beta[H(\sigma^{v,s})-H(\sigma)]},&\text{if}\ H(\sigma^{v,s})-H(\sigma) > 0,
\end{cases}\end{align}
\end{itemize}
where $\sigma^{v,s}$ is the configuration obtained from $\sigma$ by updating the spin in the vertex $v$ to $s$, i.e., 
\begin{align}\label{confspinflip}
\sigma^{v,s}(w):=
\begin{cases}
\sigma(w)\ &\text{if}\ w\neq v,\\
s\ &\text{if}\ w=v.
\end{cases}
\end{align}
We say that $\sigma\in\mathcal X$ \textit{communicates} with another configuration $\bar\sigma\in\mathcal X$ if there exist a vertex $v\in V$ and a spin value $s\in S$, such that $\sigma(v)\neq s$ and $\bar\sigma=\sigma^{v,s}$. 
Hence, at each step the update of vertex $v$ depends on the neighboring spins of $v$ and on the energy difference
\begin{align}\label{energydiff}
H(\sigma^{v,s})-H(\sigma)= \sum_{w \sim v} (\mathbbm{1}_{\{\sigma(v)=\sigma(w)\}}-\mathbbm{1}_{\{\sigma(w)=s\}}).
\end{align}
\section{Minimal restricted-gates and minimal gates}\label{secgatesmoddef}
In this section we introduce our main results on the set of minimal restricted-gates and the one of minimal gates for the transition either from a Potts stable configuration to the other Potts stable states or from a Potts stable state to another Potts stable configuration. In order to state these main results, we need to give some notations and definitions which are used throughout the next sections.
\subsection{Definitions and notations}\label{defnotgates}
We will denote the edge that links the vertices $v$ and $w$ as $(v,w)\in\ E$. Each $v\in V$ is naturally identified by its coordinates $(i,j)$, where $i$ and $j$ denote respectively the row and the column where $v$ lies. Moreover, the collection of vertices with first coordinate equal to $i=0,\dots,K-1$ is denoted as $r_i$, which is the $i$-th row of $\Lambda$. The collection of those vertices with second coordinate equal to $j=0,\dots,L-1$ is denoted as $c_j$, which is the $j$-th column of $\Lambda$.

\subsubsection{Model-independent definitions and notations}\label{modinddefnotgates}
We now give a list of model-independent definitions and notations that will be useful in formulating our main results.
\begin{itemize}
\item[-] A \textit{path} $\omega$ is a finite sequence of configurations $\omega_0,\dots,\omega_n \in \mathcal{X}$ such that $Q(\omega_i,\omega_{i+1})>0$ for $i=0,\dots,n-1$. A path from $\omega_0=\sigma$ to $\omega_n=\sigma'$ is denoted as $\omega: \sigma \to \sigma'$. Finally, $\Omega_{\sigma,\sigma'}$ denotes the set of all paths between $\sigma$ and $\sigma'$.
\item[-] The \textit{height} of a path $\omega$ is
\begin{align}\label{height}
\Phi_\omega:=\max_{i=0,\dots,n} H(\omega_i).
\end{align}
\item[-] For any pair $\sigma, \sigma' \in \mathcal{X}$, the \textit{communication height} $\Phi(\sigma,\sigma')$ between $\sigma$ and $\sigma'$ is the minimal energy across all paths $\omega:\sigma \to \sigma'$. Formally, 
\begin{align}\label{comheight}
\Phi(\sigma,\sigma'):=\min_{\omega:\sigma \to \sigma'} \Phi_\omega = \min_{\omega:\sigma \to \sigma'} \max_{\eta \in \omega} H(\eta). 
\end{align}
More generally, the communication energy between any pair of non-empty disjoint subsets $\mathcal{A},\mathcal{B} \subset \mathcal{X}$ is
\begin{align}
\Phi(\mathcal{A},\mathcal{B}):=\min_{\sigma \in \mathcal{A},\ \sigma' \in \mathcal{B}} \Phi(\sigma,\sigma').
\end{align}
\item[-] The bottom $\mathscr{F}(\mathcal{A})$ of a non-empty set $\mathcal{A}\subset\mathcal{X}$ is the set of \textit{global minima} of $H$ in $\mathcal{A}$, i.e.,
\begin{align}\label{bottom}
\mathscr{F}(\mathcal{A}):=\{\eta \in \mathcal{A}:H(\eta)=\min_{\sigma \in \mathcal{A}}H(\sigma)\}. 
\end{align}
In particular, 
\begin{align}\label{stableset}
\mathcal{X}^s:=\mathscr{F}(\mathcal{X})
\end{align}
is the set of the \textit{stable states}.
\item[-] The set of \textit{optimal paths} between $\sigma, \sigma' \in\mathcal{X}$ is defined as
\begin{align}\label{optpaths}
\Omega_{\sigma,\sigma'}^{opt}:=\{\omega\in\Omega_{\sigma,\sigma'}:\ \max_{\eta\in\omega} H(\eta)=\Phi(\sigma,\sigma')\}.
\end{align}
In other words, the optimal paths are those that realize the min-max in \eqref{comheight} between $\sigma$ and $\sigma'$.
\item[-] The set of \textit{minimal saddles} between $\sigma, \sigma' \in \mathcal{X}$ is defined as
\begin{align}\label{saddles}
\mathcal S(\sigma,\sigma'):=\{\xi\in\mathcal{X}:\exists\omega\in\Omega_{\sigma,\sigma'}^{opt},\ \xi\in\omega:\ \max_{\eta\in\omega} H(\eta)=H(\xi)\}.
\end{align}
\item[-] We say that $\eta\in\mathcal S(\sigma,\sigma')$ is an \textit{essential saddle} if either
\begin{itemize}
\item there exists $\omega\in\Omega_{\sigma,\sigma'}^{opt}$ such that $\{\text{argmax}_\omega H\}=\{\eta\}$ or
\item there exists $\omega\in\Omega_{\sigma,\sigma'}^{opt}$ such that $\{\text{argmax}_\omega H\}\supset\{\eta\}$ and $\{\text{argmax}_{\omega'} H\}\not\subseteq\{\text{argmax}_\omega H\}\backslash \{\eta\}$ for all $\omega'\in\Omega_{\sigma,\sigma'}^{opt}$.
\end{itemize}
A saddle $\eta\in\mathcal S(\sigma,\sigma')$ that is not essential is said to be \textit{unessential}. 
\item[-] Given $\sigma, \sigma' \in \mathcal{X}$, we say that $\mathcal{W}(\sigma,\sigma')$ is a \textit{gate} for the transition from $\sigma$ to $\sigma'$ if  $\mathcal{W}(\sigma,\sigma')\subseteq\mathcal S(\sigma,\sigma')$ and $\omega\cap\mathcal{W}(\sigma,\sigma')\neq\varnothing$ for all $\omega\in\Omega_{\sigma,\sigma'}^{opt}$.
\item[-] We say that  $\mathcal{W}(\sigma,\sigma')$ is a \textit{minimal gate} for the transition from $\sigma$ to $\sigma'$ if it is a minimal (by inclusion) subset of $\mathcal S(\sigma,\sigma')$ that is visited by all optimal paths, namely, it is a gate and for any $\mathcal{W}'\subset\mathcal{W}(\sigma,\sigma')$ there exists $\omega'\in\Omega_{\sigma,\sigma'}^{opt}$ such that $\omega'\cap\mathcal{W}'=\varnothing$. We denote by $\mathcal{G}=\mathcal{G}(\sigma,\sigma')$ the union of all minimal gates for the transition from $\sigma$ to $\sigma'$.
\item[-] Let $\sigma,\sigma'\in\mathcal X^s, \sigma\neq\sigma'$, we define \textit{restricted-gate} for the transition from $\sigma$ to $\sigma'$ a subset $\mathcal W_{\text{\tiny{RES}}}(\sigma,\sigma')\subset\mathcal S(\sigma,\sigma')$ which is intersected by all $\omega\in\Omega_{\sigma,\sigma'}^{opt}$ such that $\omega\cap\mathcal C^\eta=\varnothing$ for any $\eta\in\mathcal X^s\backslash\{\sigma,\sigma'\}$, where $\mathcal C^\eta$ denotes the valley of the energy landscape whose bottom is the stable state $\eta$ and whose depth is $\Phi(\sigma,\sigma')$. Such a path is said to be a \textit{restricted-path} between $\sigma$ and $\sigma'$.

We say that a restricted-gate $\mathcal W_{\text{\tiny{RES}}}(\sigma,\sigma')$ for the transition from $\sigma$ to $\sigma'$ is a \textit{minimal restricted-gate} if for any $\mathcal{W}'\subset\mathcal{W}_{\text{\tiny{RES}}}(\sigma,\sigma')$ there exists $\omega'\in\Omega_{\sigma,\sigma'}^{opt}$ such that $\omega'\cap\mathcal{W}'=\varnothing$. We denote by $\mathcal{F}(\sigma,\sigma')$ the union of all minimal restricted-gates for the transition from $\sigma$ to $\sigma'$. Note that all gates are restricted gates in the case of the Ising model with zero external magnetic field, for which $q=2$, see Corollary \ref{corgate}.
\end{itemize}
\subsubsection{Model-dependent definitions and notations}\label{moddepdefnot}
\begin{figure}[h!]
\centering
\begin{tikzpicture}[scale=0.6, transform shape]
\draw [fill=gray,lightgray] (0.6,0) rectangle (1.5,2.4);
\draw[step=0.3cm,color=black] (0,0) grid (3.6,2.7);
\draw (1.8,-0.5) node {\Large(a)};
\end{tikzpicture}\ \ \ \ \ \ \ \ \
\begin{tikzpicture}[scale=0.6, transform shape]
\draw [fill=gray,lightgray] (0.3,0) rectangle (1.5,2.1);
\draw [fill=gray,lightgray] (1.5,0.6) rectangle (1.8,1.8);
\draw[step=0.3cm,color=black] (0,0) grid (3.6,2.7);
\draw (1.8,-0.5) node {\Large(b)};
\end{tikzpicture}\ \ \ \ \ \ \ \ \
\begin{tikzpicture}[scale=0.6, transform shape]
\draw [fill=gray,lightgray] (0,0) rectangle (1.5,2.7);
\draw [fill=gray,lightgray] (1.2,0.3) rectangle (1.8,1.2);
\draw[step=0.3cm,color=black] (0,0) grid (3.6,2.7);
\draw (1.8,-0.5) node {\Large(c)};
\end{tikzpicture}
\caption{\label{esempioBRbar} Examples of configurations which belong to $\bar R_{3,8}(r, s)$ (a), $\bar B_{4,7}^4(r, s)$ (b) and $\tilde{\mathcal B}_{6,9}^{6}(r, s)$ (c). For semplicity we color white the vertices whose spin is $r$ and we color gray the vertices whose spin is $s$.}
\end{figure}
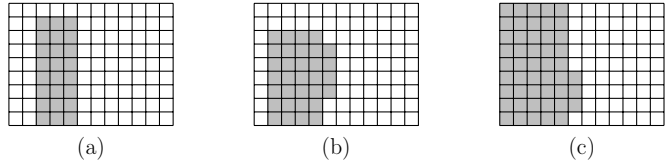

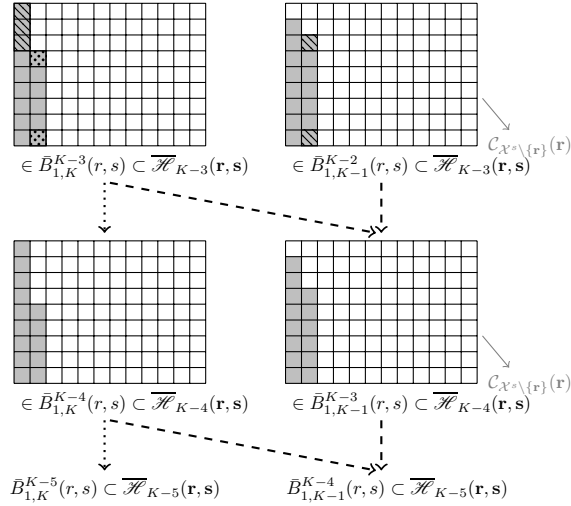
\begin{figure}[h!]
\centering
\begin{tikzpicture}[scale=0.7, transform shape]
\draw [fill=gray,lightgray] (0,0) rectangle (0.3,2.7);
\draw [fill=gray,lightgray] (0.3,0) rectangle (0.6,1.5);
\draw [fill=gray,lightgray] (5.1,0) rectangle (5.4,2.4);
\draw [fill=gray,lightgray] (5.4,0) rectangle (5.7,1.8);
\draw [fill=gray,lightgray] (0,4.5) rectangle (0.3,7.2);
\draw [fill=gray,lightgray] (0.3,4.5) rectangle (0.6,6.3);
\draw [fill=gray,lightgray] (5.1,4.5) rectangle (5.4,6.9);
\draw [fill=gray,lightgray] (5.4,4.5) rectangle (5.7,6.6);
\draw[pattern=north west lines, pattern color=black] (0,6.3) rectangle (0.3,7.2);
\draw[pattern=crosshatch dots, pattern color=black] (0.3,6) rectangle (0.6,6.3);
\draw[pattern=crosshatch dots, pattern color=black] (0.3,4.5) rectangle (0.6,4.8);
\draw[pattern=north west lines, pattern color=black] (5.4,4.5) rectangle (5.7,4.8);
\draw[pattern=north west lines, pattern color=black] (5.4,6.3) rectangle (5.7,6.6);
\draw[step=0.3cm,color=black] (0,4.5) grid (3.6,7.2) node[right] at (0,4.1){$\in\bar B_{1,K}^{K-3}(r,s)\subset\overline{\mathscr H}_{K-3}(\mathbf r,\mathbf s)$};
\draw[step=0.3cm,color=black] (5.1,4.5) grid (8.7,7.2) node[right] at (5.1,4.1){$\in\bar B_{1,K-1}^{K-2}(r,s)\subset\overline{\mathscr H}_{K-3}(\mathbf r,\mathbf s)$};
\draw[step=0.3cm,color=black] (0,0) grid (3.6,2.7) node[right] at (0,-0.4){$\in\bar B_{1,K}^{K-4}(r,s)\subset\overline{\mathscr H}_{K-4}(\mathbf r,\mathbf s)$};
\draw[step=0.3cm,color=black] (5.1,0) grid (8.7,2.7) node[right] at (5.1,-0.4){$\in\bar B_{1,K-1}^{K-3}(r,s)\subset\overline{\mathscr H}_{K-4}(\mathbf r,\mathbf s)$};
\draw [->,dotted,thick] (1.7, 3.8) -- (1.7,2.85);
\draw [->,dashed,thick] (1.8, 3.8) -- (6.8,2.85);
\draw [->,dashed,thick] (6.9, 3.8) -- (6.9, 2.85);
\draw [->,gray] (8.8,0.9) -- (9.3,0.3) node[below] at (9.7,0.3){$\mathcal C_{\mathcal X^s\backslash\{\mathbf r\}}(\mathbf r)$}; 
\draw [->,gray] (8.8,5.4) -- (9.3,4.8) node[below] at (9.7,4.8){$\mathcal C_{\mathcal X^s\backslash\{\mathbf r\}}(\mathbf r)$};
\draw [->,dotted,thick] (1.7, -0.7) -- (1.7,-1.65) node[below] at (1.9,-1.65){$\bar B_{1,K}^{K-5}(r,s)\subset\overline{\mathscr H}_{K-5}(\mathbf r,\mathbf s)$};
\draw [->,dashed,thick] (1.8, -0.7) -- (6.8,-1.65) node[below] at (7.2,-1.65){$\bar B_{1,K-1}^{K-4}(r,s)\subset\overline{\mathscr H}_{K-5}(\mathbf r,\mathbf s)$};
\draw [->,dashed,thick] (6.9, -0.7) -- (6.9, -1.65);
\draw [->,white] (-0.1,0.9) -- (-0.8,0.3) node[below] at (-1,0.4){\large$\mathcal C_{\mathcal X^s\backslash\{\mathbf r\}}(\mathbf r)$}; 
\draw [->,white] (-0.1,5.4) -- (-0.8,4.8) node[below] at (-1,4.9){\large$\mathcal C_{\mathcal X^s\backslash\{\mathbf r\}}(\mathbf r)$};
\end{tikzpicture}
\caption{\label{step2} Example of configuations belonging to $\overline{\mathscr H}(\mathbf r,\mathbf s)$ on a $9\times 12$ grid $\Lambda$. White vertices have spin $r$, gray vertices have spin $s$. By flipping a spin $s$ to $r$ both among those with dots and among those with lines, the path can pass to another configuration belonging to $\overline{\mathscr H}(\mathbf r,\mathbf s)$.}
\end{figure}

Next we give some further model-dependent notations in order to be able to state our main results.
The definitions hold for any $q$-Potts configuration $\sigma\in\mathcal X$ and any two different spin values $r,s\in\{1,\dots,q\}$.
\begin{itemize}
\item[-] $\bar R_{a,b}(r, s)$ denotes the set of those configurations in which all the vertices have spins equal to $r$, except those, which have spins $s$, in a rectangle $a\times b$, where $a$ is the horizontal length and $b$ is the vertical length, see Figure \ref{esempioBRbar}(a);
\item[-] $\bar B_{a,b}^h(r, s)$ denotes the set of those configurations in which all the vertices have spins $r$, except those, which have spins $s$, in a rectangle $a\times b$, where $a$ is the horizontal length and $b$ is the vertical length, with a bar $1\times h$ adjacent to one of the sides of length $b$, with $1\le h\le b-1$, see Figure \ref{esempioBRbar}(b).
\item[-] Analogously, we set $\tilde R_{a,b}(r, s)$ and $\tilde B_{a,b}^h(r, s)$ interchanging the role of spins $r$ and $s$, see Figure \ref{esempioBRbar}(c).
\end{itemize}
Note that 
\begin{align}\label{equivalenzarbk}
\bar R_{a,K}(r, s)\equiv\tilde R_{L-a,K}(r, s)\ \ \ \text{and}\ \ \ \bar B_{a,K}^h(r, s)\equiv\tilde B_{L-a-1,K}^{K-h}(r, s).
\end{align}
Next, we define sets of configurations that are crucial to describe the gate. We show
their location on the energy landscape in Figure \ref{ciclieplateau}.
\begin{itemize}
\item[-] We set
\begin{align}\label{definsP}
\overline{\mathscr P}(\mathbf r,\mathbf s):=\bar B_{1,K}^{K-1}(r, s),\ \ \ \ \ \widetilde{\mathscr P}(\mathbf r,\mathbf s):=\tilde B_{1,K}^{K-1}(r, s).
\end{align}
We refer the reader to Figure \ref{figdaPaQ} for an example of a configuration belonging to $\overline{\mathscr P}(\mathbf r,\mathbf s)$.
\item[-] We define
\begin{align}\label{definsQ}
\overline{\mathcal Q}(\mathbf r,\mathbf s):=\bar R_{2,K-1}(r, s)\cup\bar B_{1,K}^{K-2}(r, s), \\ \widetilde{\mathcal Q}(\mathbf r,\mathbf s):=\tilde R_{2,K-1}(r, s)\cup\tilde B_{1,K}^{K-2}(r, s).
\end{align}
\item[-] For any $i=1,\dots,K-3$,
 we define
\begin{align}\label{definsHbar}
\overline{\mathscr H}_i(\mathbf r,\mathbf s):=\bar B_{1,K}^i(r, s)\cup\bigcup_{h=i+1}^{K-2} \bar B_{1,K-1}^h(r, s),\ \ \overline{\mathscr H}(\mathbf r,\mathbf s) := \bigcup_{i=1}^{K-3}\overline{\mathscr H}_i(\mathbf r,\mathbf s)
\\ \label{definsHtilde}
\widetilde{\mathscr H}_i(\mathbf r,\mathbf s):=\tilde B_{1,K}^i(r, s)\cup\bigcup_{h=i+1}^{K-2} \tilde B_{1,K-1}^h(r, s),\ \ \widetilde{\mathscr H}(\mathbf r,\mathbf s) := \bigcup_{i=1}^{K-3}\widetilde{\mathscr H}_i(\mathbf r,\mathbf s).
\end{align}
We refer the reader to Figure \ref{figdaPaQ} and Figure \ref{step2} for an example of a configurations belonging to $\overline{\mathcal Q}(\mathbf r,\mathbf s)$ and to $\overline{\mathscr H}(\mathbf r,\mathbf s)$. See Figure \ref{passdaQaQH} and Figure \ref{stepfinale} for other examples of this type of configurations.
\end{itemize}
\begin{figure}[h!]
\centering
\begin{tikzpicture}[scale=0.7, transform shape]
\draw [fill=gray,lightgray] (0,0) rectangle (0.3,2.7);
\draw [fill=gray,lightgray] (0.3,0) rectangle (0.6,2.1);
\draw [fill=gray,lightgray] (4.8,0) rectangle (5.4,2.4);
\draw [fill=gray,lightgray] (2.4,3.6) rectangle (2.7,6.3);
\draw [fill=gray,lightgray] (2.7,3.6) rectangle (3,6);
\draw[pattern=north west lines, pattern color=black] (2.4,6) rectangle (2.7,6.3);
\draw[pattern=crosshatch dots, pattern color=black] (2.7,3.6) rectangle (3,3.9);
\draw[pattern=crosshatch dots, pattern color=black] (2.7,5.7) rectangle (3,6);
\draw[step=0.3cm,color=black] (2.4,3.6) grid (6,6.3) node[below] at (4.2,3.6){$\in\overline{\mathscr P}(\mathbf r,\mathbf s)$};
\draw[step=0.3cm,color=black] (0,0) grid (3.6,2.7) node[right] at (0,-0.4){$\in\bar B_{1,K}^{K-2}(r,s)\subset\overline{\mathcal Q}(\mathbf r,\mathbf s)$};

\draw[step=0.3cm,color=black] (4.8,0) grid (8.4,2.7) node[right] at (4.6,-0.4){$\in\bar R_{2,K-1}(r,s)\subset\overline{\mathcal Q}(\mathbf r,\mathbf s)$};
\draw [->,dotted,thick] (2.2, 3.6) -- (1.5,2.9);

\draw [->,dashed,thick] (6.2,3.6) -- (6.9,2.9);

\end{tikzpicture}
\caption{\label{figdaPaQ} Example of configuations belonging to $\overline{\mathscr P}(\mathbf r,\mathbf s)$ and $\overline{\mathcal Q}(\mathbf r,\mathbf s)$ on  a $9\times 12$ grid $\Lambda$. Gray vertices have spin value $s$, white vertices have spin value $r$. By flipping to $r$ a spin $s$ among those with the lines, the path enters into $\bar B_{1,K}^{K-2}(r,s)\subset\overline{\mathcal Q}(\mathbf r,\mathbf s)$; instead, by flipping to $r$ a spin $s$ among those with dots, the path goes to $\bar R_{2,K-1}(r,s)\subset\overline{\mathcal Q}(\mathbf r,\mathbf s)$.}
\end{figure}
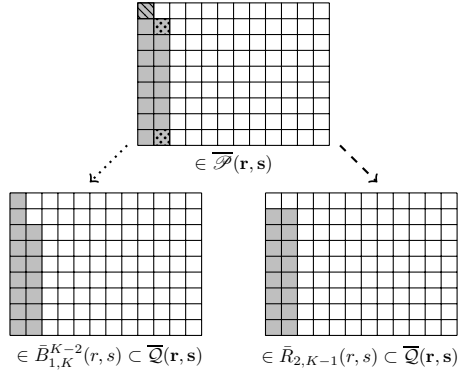\FloatBarrier

\FloatBarrier
\begin{itemize}
\item[-] Finally, we set 
\begin{align}\label{definsWjh}
&\mathcal{W}_j^{(h)}(\mathbf r,\mathbf s):=\bar B_{j,K}^h(r, s)=\tilde B_{L-j-1,K}^{K-h}(r, s)\ \ \ \text{for}\ j=2,\dots,L-3,\\
&\mathcal{W}_j(\mathbf r,\mathbf s):=\bigcup_{h=1}^{K-1} \mathcal{W}_j^{(h)}(\mathbf r,\mathbf s).
\end{align}
\end{itemize}
We refer to Figure \ref{esempioBRbar}(c) for an example of configuration belonging to $\mathcal W_5^{(3)}(\mathbf r,\mathbf s)$.

\subsection{Main results}
We are now ready to state the main results on minimal restricted-gates between two stable configurations and on minimal gates for the transition between a stable configuration and the other stable states and between a stable configuration and another one.

\subsubsection{Minimal restricted-gates between two Potts stable configurations}\label{subgateres}
In Section \ref{minimalresgatessection} we study the energy landscape between two given stable configurations $\mathbf r,\mathbf s\in\mathcal{X}^s$, $\mathbf s\neq\mathbf r$ and describe the set of all minimal restricted-gates for the transition between them. We recall that these gates are said to be ``restricted'' because they are gates for the transition from $\mathbf r$ to $\mathbf s$ following an optimal path $\omega\in\Omega_{\mathbf r,\mathbf s}^{opt}$ such that $\omega\cap(\mathcal X^s\backslash\{\mathbf r,\mathbf s\})=\varnothing$. 
We give two representations of the restricted-gates within the energy landscape. In Figure \ref{ciclieplateau} we give a side view of the energy landscape between two stable configurations \textbf{r} and \textbf{s}, and we draw the restricted-gates corresponding to the transition between these two configurations. In Figure \ref{figq5alto} we give a top-down view of the energy landscape between several stable configurations. Figure \ref{ciclieplateau} then is a side view of any one of the four arms in Figure \ref{figq5alto}. Accordingly, studying the restricted-gates between, say \textbf{1} and \textbf{2}, corresponds to focusing on only those paths that cross the right part of Figure \ref{figq5alto}. The following results identify the minimal restricted-gates between two stable configurations.

\begin{theorem}[Minimal-restricted gates]\label{mingatescond}
Consider the $q$-state Potts model on a $K\times L$ grid $\Lambda$ with $\max\{K,L\}\ge3$ and with periodic boundary conditions. For every $\mathbf r,\mathbf s\in\mathcal{X}^s$, $\mathbf s\neq\mathbf r$, the following sets are minimal restricted-gates for the transition $\mathbf r\to\mathbf s$:
\begin{itemize}
\item[\emph{(a)}]  $\overline{\mathscr P}(\mathbf r,\mathbf s)$ and $\widetilde{\mathscr P}(\mathbf r,\mathbf s)$;
\item[\emph{(b)}] $\overline{\mathcal Q}(\mathbf r,\mathbf s)$ and $\widetilde{\mathcal Q}(\mathbf r,\mathbf s)$;
\item[\emph{(c)}] $\overline{\mathscr H}_i(\mathbf r,\mathbf s)$ and $\widetilde{\mathscr H}_i(\mathbf r,\mathbf s)$ for any $i=1,\dots,K-3$;
\item[\emph{(d)}] $\mathcal W_j^{(h)}(\mathbf r,\mathbf s)$ for any $j=2,\dots,L-3$ and any $h=1,\dots,K-1$.
\end{itemize}
\end{theorem}
\begin{theorem}[Union of all minimal-restricted gates]\label{mingatescondset} 
For every $\mathbf r,\mathbf s\in\mathcal{X}^s$, $\mathbf s\neq\mathbf r$, the union of all minimal restricted-gates for the transition $\mathbf r\to\mathbf s$ is given by
\begin{align}\label{setmingatestwostable}
\mathcal F(\mathbf r,\mathbf s)=\bigcup_{j=2}^{L-3} \mathcal W_j(\mathbf r,\mathbf s)\cup\overline{\mathscr H}(\mathbf r,\mathbf s)\cup\widetilde{\mathscr H}(\mathbf r,\mathbf s)\cup\overline{\mathcal Q}(\mathbf r,\mathbf s)\cup\widetilde{\mathcal Q}(\mathbf r,\mathbf s)\cup\overline{\mathscr P}(\mathbf r,\mathbf s)\cup \widetilde{\mathscr P}(\mathbf r,\mathbf s).
\end{align}
\end{theorem}

 Given a non-empty subset $\mathcal{A} \subset \mathcal{X}$ and a configuration $\sigma \in \mathcal{X}$, we define $\tau_\mathcal{A}^\sigma := \text{inf}\{t>0: \ X_t^\beta \in \mathcal{A}\}$ as the \textit{first hitting time} of the subset $\mathcal{A}$ for the Markov chain $\{X_t^\beta\}_{t \in \mathbb{N}}$ starting from $\sigma$ at time $t=0$.
\begin{cor}[Crossing the gates]\label{corgatecond}
Consider any $\mathbf r, \mathbf s\in\mathcal{X}^s$ and the transition from $\mathbf r$ to $\mathbf s$. Then, the following properties hold:
\begin{itemize}
\item[\emph{(a)}] $\lim_{\beta\to\infty} \mathbb P_\beta(\tau^{\mathbf r}_{\overline{\mathscr P}(\mathbf r,\mathbf s)}<\tau^{\mathbf r}_{\mathcal X^s\backslash\{\mathbf r\}}|\tau^{\mathbf r}_{\mathbf s}<\tau^{\mathbf r}_{\mathcal X^s\backslash\{\mathbf r,\mathbf s\}})$

 \hspace{30pt} $=\lim_{\beta\to\infty} \mathbb P_\beta(\tau^{\mathbf r}_{\widetilde{\mathscr P}(\mathbf r,\mathbf s)}<\tau^{\mathbf r}_{\mathcal X^s\backslash\{\mathbf r\}}|\tau^{\mathbf r}_{\mathbf s}<\tau^{\mathbf r}_{\mathcal X^s\backslash\{\mathbf r,\mathbf s\}})=1$;
\item[\emph{(b)}] $\lim_{\beta\to\infty} \mathbb P_\beta(\tau^{\mathbf r}_{\overline{\mathcal Q}(\mathbf r,\mathbf s)}<\tau^{\mathbf r}_{\mathcal X^s\backslash\{\mathbf r\}}|\tau^{\mathbf r}_{\mathbf s}<\tau^{\mathbf r}_{\mathcal X^s\backslash\{\mathbf r,\mathbf s\}})$

\hspace{30pt}$=\lim_{\beta\to\infty} \mathbb P_\beta(\tau^{\mathbf r}_{\widetilde{\mathcal Q}(\mathbf r,\mathbf s)}<\tau^{\mathbf r}_{\mathcal X^s\backslash\{\mathbf r\}}|\tau^{\mathbf r}_{\mathbf s}<\tau^{\mathbf r}_{\mathcal X^s\backslash\{\mathbf r,\mathbf s\}})=1$;
\item[\emph{(c)}] $\lim_{\beta\to\infty} \mathbb P_\beta(\tau^{\mathbf r}_{\overline{\mathscr H}_i(\mathbf r,\mathbf s)}<\tau^{\mathbf r}_{\mathcal X^s\backslash\{\mathbf r\}}|\tau^{\mathbf r}_{\mathbf s}<\tau^{\mathbf r}_{\mathcal X^s\backslash\{\mathbf r,\mathbf s\}})$

\hspace{30pt}$=\lim_{\beta\to\infty} \mathbb P_\beta(\tau^{\mathbf r}_{\widetilde{\mathscr H}_i(\mathbf r,\mathbf s)}<\tau^{\mathbf r}_{\mathcal X^s\backslash\{\mathbf r\}}|\tau^{\mathbf r}_{\mathbf s}<\tau^{\mathbf r}_{\mathcal X^s\backslash\{\mathbf r,\mathbf s\}})=1$, for any $i=1,\dots,K-3$;
 \item[\emph{(d)}] $\lim_{\beta\to\infty} \mathbb P_\beta(\tau^{\mathbf r}_{\mathcal W_j^{(h)}(\mathbf r,\mathbf s)}<\tau^{\mathbf r}_{\mathcal X^s\backslash\{\mathbf r\}}|\tau^{\mathbf r}_{\mathbf s}<\tau^{\mathbf r}_{\mathcal X^s\backslash\{\mathbf r,\mathbf s\}})=1$ for any $j=2,\dots,L-3$, $h=1,\dots,K-1$.
\end{itemize}
\end{cor}
 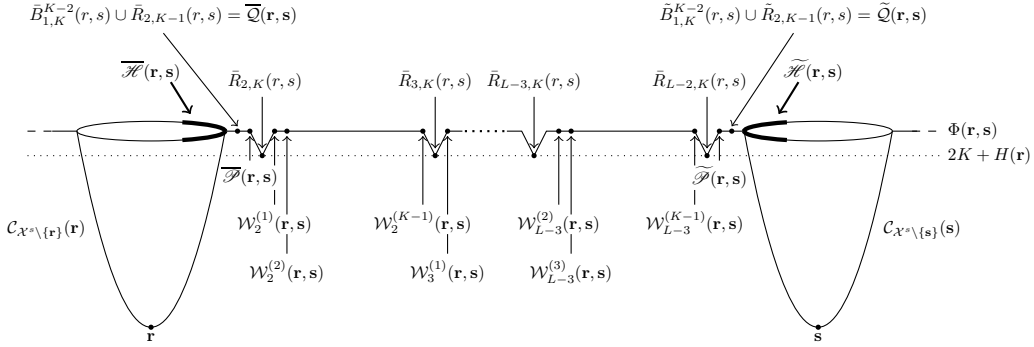
\begin{figure}[h!]
\begin{minipage}[c]{1.\textwidth}
    \centering
    \makebox[0pt]{%
\begin{tikzpicture}[scale=0.65, transform shape]
\draw[dashed] (-2.5,4) -- (-2,4);
\draw (-2,4) -- (-1.5,4);
\draw (0,0) parabola (1.5,4);
\draw (0,0) parabola (-1.5,4) ;
\draw (1.5,4) -- (2,4) -- (2.25,3.5) -- (2.5,4) -- (5.5,4) -- (5.75,3.5) -- (6,4) -- (6.3,4);
\draw[dotted,thick] (6.35,4) -- (7.15,4);
\draw (7.2,4) -- (7.5,4) -- (7.75,3.5) -- (8,4) -- (11,4) -- (11.25,3.5) -- (11.5,4) -- (12,4);
\draw (13.5,0) parabola (12,4);
\draw (13.5,0) parabola (15,4);
\draw (15,4) -- (15.5,4);
\draw[dashed] (15.4,4) -- (16,4);
\draw (0,4) ellipse (1.5cm and 0.2cm);
\draw[ultra thick] (1.5,4) arc(0:65:1.5cm and 0.2cm);
\draw[ultra thick] (1.5,4) arc(0:-65:1.5cm and 0.2cm);
\draw (13.5,4) ellipse (1.5cm and 0.2cm);
\draw[ultra thick] (12,4) arc(180:115:1.5cm and 0.2cm);
\draw[ultra thick] (12,4) arc(180:245:1.5cm and 0.2cm);
\draw[<-,thick] (0.8,4.35) -- (0.4,5);\draw(0,4.9) node[above] {$\overline{\mathscr H}(\mathbf r,\mathbf s)$};
\draw[<-,thick] (12.7,4.35) -- (13.1,5);\draw (13.4,4.9) node[above] {$\widetilde{\mathscr H}(\mathbf r,\mathbf s)$};
\fill[color=black] (0,0) circle (1.5pt) node[below] {$\mathbf r$};
\fill[color=black] (13.5,0) circle (1.5pt) node[below] {$\mathbf s$};
\draw[<-] (2,3.9) -- (2,3.4) node[below] {$\overline{\mathscr P}(\mathbf r,\mathbf s)$};
\fill[color=black] (2,4) circle (1.5pt); 
\fill[color=black] (1.75,4) circle (1.5pt);
\draw[<-] (11.5,3.9) -- (11.5,3.4) node[below] {$\widetilde{\mathscr P}(\mathbf r,\mathbf s)$};
\fill[color=black] (11.5,4) circle (1.5pt);
\fill[color=black] (11.75,4) circle (1.5pt) node[above] {};
\fill[color=black] (2.25,3.5) circle (1.5pt); \draw [<-] (2.25,3.65) -- (2.25,4.7) node[above] {$\ \bar R_{2,K}(r,s)$};
\fill[color=black] (5.75,3.5) circle (1.5pt); \draw[<-] (5.75,3.65) -- (5.75,4.7) node[above] {$\bar R_{3,K}(r,s)$};
\fill[color=black] (7.75,3.5) circle (1.5pt); \draw[<-] (7.75,3.65) -- (7.75,4.7) node[above] {$\bar R_{L-3,K}(r,s)$};
\fill[color=black] (11.25,3.5) circle (1.5pt) (11.1,3.5); \draw[<-] (11.25,3.65) -- (11.25,4.7);\draw(11.1,4.7) node[above] {$\bar R_{L-2,K}(r,s)$};
\draw[<-] (1.75,4.1) -- (0.25,6) node[above] {$\bar B_{1,K}^{K-2}(r,s)\cup\bar R_{2,K-1}(r,s)=\overline{\mathcal Q}(\mathbf r,\mathbf s)$};
\draw[<-] (11.75,4.1) -- (13,6) node[above] {$\tilde B_{1,K}^{K-2}(r,s)\cup\tilde R_{2,K-1}(r,s)=\widetilde{\mathcal Q}(\mathbf r,\mathbf s)$};
\fill[color=black] (2.5,4) circle (1.5pt) (2.75,4) circle (1.5pt) (5.5,4) circle (1.5pt) (11,4) circle (1.5pt) (6,4) circle (1.5pt) (8.25,4) circle(1.5pt) (8.5,4) circle (1.5pt);
\draw(2.5,2.8) -- (2.5,2.5) [<-] (2.5,3.9) -- (2.5,3.3);
\draw(2.75,1.9) -- (2.75,1.5) [<-] (2.75,3.9) -- (2.75,2.3);
\draw [<-] (5.5,3.9) -- (5.5,2.5);
\draw (6,1.9) -- (6,1.5) [<-] (6,3.9) -- (6,2.3);
\draw [<-] (8.25,3.9) -- (8.25,2.5);
\draw (8.5,1.9) -- (8.5,1.5) [<-] (8.5,3.9) -- (8.5,2.3);
\draw(11,2.8) -- (11,2.5) [<-] (11,3.9) -- (11,3.3);
\draw (2.5,2.5) node[below] {$\mathcal W_2^{(1)}(\mathbf r,\mathbf s)$};
\draw (2.75,1.5) node[below] {$\mathcal W_2^{(2)}(\mathbf r,\mathbf s)$};
\draw (5.5,2.5) node[below] {$\mathcal W_2^{(K-1)}(\mathbf r,\mathbf s)$};
\draw (6,1.5) node[below] {$\mathcal W_3^{(1)}(\mathbf r,\mathbf s)$};
\draw (8.25,2.5) node[below] {$\mathcal W_{L-3}^{(2)}(\mathbf r,\mathbf s)$};
\draw (8.5,1.5) node[below] {$\mathcal W_{L-3}^{(3)}(\mathbf r,\mathbf s)$};
\draw (11,2.5) node[below] {$\mathcal W_{L-3}^{(K-1)}(\mathbf r,\mathbf s)$};
\draw (-1.2,2) node[left] {$\mathcal C_{\mathcal X^s\backslash\{\mathbf r\}}(\mathbf r)$};
\draw (14.7,2) node[right] {$\mathcal C_{\mathcal X^s\backslash\{\mathbf s\}}(\mathbf s)$};
\draw (16,4) node[right] {$\Phi(\mathbf r,\mathbf s)$};
\draw[dotted] (-2.5,3.5) -- (16,3.5);
\draw (16,3.5) node[right] {$2K+H(\mathbf r)$};
\end{tikzpicture}%
    }\par
     \end{minipage}
     \caption{\label{ciclieplateau} Focus on the energy landscape between $\mathbf r$ and $\mathbf s$ and example of some essential saddles for the transition $\mathbf r\to\mathbf s$ following an optimal path which does not pass through other stable states.}
     \end{figure}\FloatBarrier
\subsubsection{Minimal gates for the transition from a stable state to the other stable states}\label{subgateset}
Using the results about the minimal restricted-gates, in Theorem \ref{mingatesNOcond} we identify geometrically all the sets of minimal gates for the transition from a stable configuration to the other stable states. While in Theorem \ref{setmingatesNOcond} we identify the union of all minimal gates for the same transition.
We assume $q>2$, otherwise when $q=2$, $|\mathcal X^s|=2$ and Theorems \ref{mingatesNOcond}--\ref{setmingatesNOcond} coincide with Theorems \ref{mingatescond}--\ref{mingatescondset}.

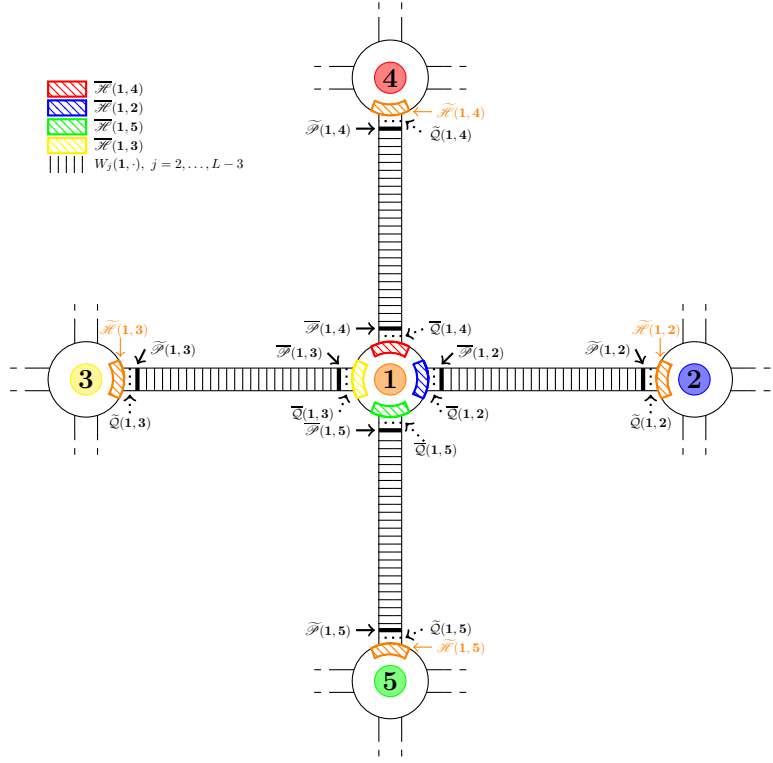
\begin{figure}[h!]
\begin{minipage}[c]{1\textwidth}
    \centering
    \makebox[0pt]{%
\begin{tikzpicture}[scale=0.5, transform shape]
\draw (0,0) circle (1cm) (0,8) circle (1cm) (8,0) circle (1cm) (-8,0) circle (1cm) (0,-8) circle (1cm);
\draw [thick,blue] (8,0) circle (0.4cm); \fill [blue!50!white] (8,0) circle (0.4cm);
\draw [thick,yellow] (-8,0) circle (0.4cm);  \fill [yellow!50!white]  (-8,0) circle (0.4cm);
\draw [thick, red] (0,8) circle (0.4cm);  \fill [red!50!white] (0,8) circle (0.4cm); 
\draw [thick,green] (0,-8) circle (0.4cm);  \fill [green!50!white]  (0,-8) circle (0.4cm);
\draw [thick,orange] (0,0) circle (0.4cm); \fill [orange!50!white] (0,0) circle (0.4cm);
\draw (0,0) node{\huge\textbf{$\mathbf 1$}} (8,0) node{\huge\textbf{$\mathbf 2$}} (-8,0) node{\huge\textbf{$\mathbf 3$}} (0,8) node{\huge\textbf{$\mathbf 4$}} (0,-8) node{\huge\textbf{$\mathbf 5$}};
\draw (-7.05,-0.3) -- (-0.95,-0.3) (-7.05,0.3) -- (-0.95,0.3) (0.95,-0.3) -- (7.05,-0.3) (0.95,0.3) -- (7.05,0.3);
\draw (0.3,7.05) -- (0.3,0.95) (-0.3,7.05) -- (-0.3,0.95) (0.3,-7.05) -- (0.3,-0.95) (-0.3,-7.05) -- (-0.3,-0.95);
\draw (-8.95,0.3) -- (-9.4,0.3) (-8.95,-0.3) -- (-9.4,-0.3) (8.95,0.3) -- (9.4,0.3) (8.95,-0.3) -- (9.4,-0.3);
\draw (0.3,8.95) -- (0.3,9.4) (-0.3,8.95) -- (-0.3,9.4) (0.3,-8.95) -- (0.3,-9.4) (-0.3,-8.95) -- (-0.3,-9.4);
\draw[dashed] (9.4,0.3) -- (10,0.3) (9.4,-0.3) -- (10,-0.3) (-9.4,0.3) -- (-10,0.3) (-9.4,-0.3) -- (-10,-0.3);
\draw[dashed] (0.3,9.4) -- (0.3,10) (-0.3,9.4) -- (-0.3,10) (0.3,-9.4) -- (0.3,-10) (-0.3,-9.4) -- (-0.3,-10);
\draw (7.7,0.95) -- (7.7,1.4) (8.3,0.95) -- (8.3,1.4) (-7.7,0.95) -- (-7.7,1.4) (-8.3,0.95) -- (-8.3,1.4);
\draw (-7.7,-0.95) -- (-7.7,-1.4) (-8.3,-0.95) -- (-8.3,-1.4) (7.7,-0.95) -- (7.7,-1.4) (8.3,-0.95) -- (8.3,-1.4);
\draw[dashed] (-7.7,-1.4) -- (-7.7,-2) (-8.3,-1.4) -- (-8.3,-2)  (7.7,-1.4) -- (7.7,-2) (8.3,-1.4) -- (8.3,-2);
\draw[dashed] (-7.7,1.4) -- (-7.7,2) (-8.3,1.4) -- (-8.3,2)  (7.7,1.4) -- (7.7,2) (8.3,1.4) -- (8.3,2);
\draw (-0.95,7.7) -- (-1.4,7.7) (-0.95,8.3) -- (-1.4,8.3) (0.95,7.7) -- (1.4,7.7) (0.95,8.3) -- (1.4,8.3);
\draw (-0.95,-7.7) -- (-1.4,-7.7) (-0.95,-8.3) -- (-1.4,-8.3) (0.95,-7.7) -- (1.4,-7.7) (0.95,-8.3) -- (1.4,-8.3);
\draw[dashed] (-1.4,7.7) -- (-2,7.7) (-1.4,8.3) -- (-2,8.3) (1.4,7.7) -- (2,7.7) (1.4,8.3) -- (2,8.3);
\draw[dashed] (-1.4,-7.7) -- (-2,-7.7) (-1.4,-8.3) -- (-2,-8.3) (1.4,-7.7) -- (2,-7.7) (1.4,-8.3) -- (2,-8.3);

\draw[ultra thick] (-6.65,0.3) -- (-6.65,-0.3)  (-1.35,0.3) -- (-1.35,-0.3) (1.35,0.3) -- (1.35,-0.3) (6.65,0.3) -- (6.65,-0.3);
\draw[ultra thick] (0.3,6.65) -- (-0.3,6.65) (0.3,-6.65) -- (-0.3,-6.65) (0.3,1.35) -- (-0.3,1.35) (0.3,-1.35) -- (-0.3,-1.35);

\draw[dotted,thick] (-6.85,0.3) -- (-6.85,-0.3);
\draw[dotted,thick]  (-1.15,0.3) -- (-1.15,-0.3) (1.15,0.3) -- (1.15,-0.3) (6.85,0.3) -- (6.85,-0.3);
\draw[dotted, thick] (0.3,6.85) -- (-0.3,6.85) (0.3,-6.85) -- (-0.3,-6.85) (0.3,1.15) -- (-0.3,1.15) (0.3,-1.15) -- (-0.3,-1.15);

\draw[<-,thick] (-6.65,0.4) -- (-6.45,0.8); \draw (-6.45,0.8) node[right] {$\widetilde{\mathscr P}(\mathbf 1,\mathbf 3)$};
\draw[<-,orange] (-7.1,0.6) -- (-7.1,1.1); \draw (-7,1) node[above,orange] {$\widetilde{\mathscr H}(\mathbf 1,\mathbf 3)$};
\draw[<-,thick] (6.65,0.4) -- (6.45,0.8); \draw (6.45,0.8) node[left] {$\widetilde{\mathscr P}(\mathbf 1,\mathbf 2)$};
\draw[<-,orange] (7.1,0.6) -- (7.1,1.1); \draw (7,1) node[above,orange] {$\widetilde{\mathscr H}(\mathbf 1,\mathbf 2)$};
\draw[<-,thick] (-0.4,6.65) -- (-0.9,6.65) node[left] {$\widetilde{\mathscr P}(\mathbf 1,\mathbf 4)$};
\draw[<-,thick] (-0.4,-6.65) -- (-0.9,-6.65) node[left] {$\widetilde{\mathscr P}(\mathbf 1,\mathbf 5)$};
\draw[<-,dotted,thick] (0.4,-6.85) -- (0.9,-6.6) node[right] {$\widetilde{\mathcal Q}(\mathbf 1,\mathbf 5)$};
\draw[<-,thick] (-0.4,-1.35) -- (-0.9,-1.35) node[left] {$\overline{\mathscr P}(\mathbf 1,\mathbf 5)$};
\draw[<-,dotted,thick] (6.85,-0.4) -- (6.85,-0.9); \draw (6.85,-0.8) node[below] {$\widetilde{\mathcal Q}(\mathbf 1,\mathbf 2)$};
\draw[<-,dotted,thick] (0.4,1.15) -- (0.9,1.35) node[right] {$\overline{\mathcal Q}(\mathbf 1,\mathbf 4)$};
\draw[<-,dotted,thick] (0.4,-1.15) -- (0.9,-1.65); \draw (1.2,-1.55) node[right,below] {$\overline{\mathcal Q}(\mathbf 1,\mathbf 5)$};
\draw[<-,dotted,thick] (1.15,-0.4) -- (1.35,-0.9) node[right] {$\overline{\mathcal Q}(\mathbf 1,\mathbf 2)$};
\draw[<-,dotted,thick] (-1.15,-0.4) -- (-1.35,-0.9) node[left] {$\overline{\mathcal Q}(\mathbf 1,\mathbf 3)$};
\draw[<-,thick] (-1.35,0.4) -- (-1.65,0.7) node[left] {$\overline{\mathscr P}(\mathbf 1,\mathbf 3)$};
\draw[<-,thick] (1.35,0.4) -- (1.65,0.7) node[right] {$\overline{\mathscr P}(\mathbf 1,\mathbf 2)$};
\draw[<-,dotted,thick] (-6.85,-0.4) -- (-6.85,-0.9); \draw (-6.85,-0.8) node[below] {$\widetilde{\mathcal Q}(\mathbf 1,\mathbf 3)$};
\draw[<-,dotted,thick] (0.4,6.85) -- (0.9,6.5) node[right] {$\widetilde{\mathcal Q}(\mathbf 1,\mathbf 4)$};
\draw[<-,thick] (-0.4,1.35)--(-0.9,1.35)node[left] {$\overline{\mathscr P}(\mathbf 1,\mathbf 4)$};
\draw[<-,orange] (0.6,7.1) -- (1.1,7.1) node[right,orange]  {$\widetilde{\mathscr H}(\mathbf 1,\mathbf 4)$};
\draw[<-,orange] (0.6,-7.1) -- (1.1,-7.1) node[right,orange]  {$\widetilde{\mathscr H}(\mathbf 1,\mathbf 5)$};
\fill[pattern=horizontal lines, pattern color=black] (-0.3,1.45)rectangle(0.3,6.55);
\fill[pattern=horizontal lines, pattern color=black] (-0.3,-1.45)rectangle(0.3,-6.55);
\fill[pattern=vertical lines, pattern color=black] (1.45,-0.3) rectangle (6.55,0.3);
\fill[pattern=vertical lines, pattern color=black] (-1.45,-0.3) rectangle (-6.55,0.3);

\draw[thick,red] (0,1) arc(90:60:1cm and 1cm) (0,1) arc(90:120:1cm and 1cm);
\draw[pattern=north west lines, pattern color=red]  (0,1) arc(90:60:1cm and 1cm) -- (0.385,0.6) arc(60:120:0.8cm);
\draw[pattern=north west lines, pattern color=red] (0,1) arc(90:120:1cm and 1cm) -- (-0.4,0.6);
\draw[thick,red] (0,1) arc(90:120:1cm and 1cm) -- (-0.4,0.6);
\draw[thick,red] (0,1) arc(90:60:1cm and 1cm) -- (0.385,0.6) arc(60:120:0.8cm);
\draw[pattern=north west lines, pattern color=orange]  (0,7) arc(-90:-60:1cm and 1cm) -- (0.385,7.4) arc(300:240:0.8cm);
\draw[pattern=north west lines, pattern color=orange] (0,7) arc(-90:-120:1cm and 1cm) -- (-0.4,7.4);
\draw[thick,orange] (0,7) arc(-90:-60:1cm and 1cm) -- (0.385,7.4) arc(300:240:0.8cm);
\draw[thick, orange] (0,7) arc(-90:-120:1cm and 1cm) -- (-0.4,7.4);
\draw[thick,yellow] (-1,0) arc(180:150:1cm and 1cm) (-1,0) arc(180:210:1cm and 1cm);
\draw[pattern=north west lines, pattern color=yellow]  (-1,0) arc(180:150:1cm and 1cm) -- (-0.6,0.385) arc(150:210:0.8cm);
\draw[pattern=north west lines, pattern color=yellow] (-1,0) arc(180:210:1cm and 1cm) -- (-0.6,-0.4);
\draw[thick,yellow] (-1,0) arc(180:150:1cm and 1cm)  (-1,0) arc(180:150:1cm and 1cm) -- (-0.6,0.385) arc(150:210:0.8cm);
\draw[thick,yellow] (-1,0) arc(180:150:1cm and 1cm) (-1,0) arc(180:210:1cm and 1cm) -- (-0.6,-0.4);
\draw[pattern=north west lines, pattern color=orange]  (-7,0) arc(0:-30:1cm and 1cm) -- (-7.4,-0.385) arc(-30:30:0.8cm);
\draw[pattern=north west lines, pattern color=orange] (-7,0) arc(0:30:1cm and 1cm) -- (-7.4,0.4);
\draw[thick,orange] (-7,0) arc(0:30:1cm and 1cm)  (-7,0) arc(0:-30:1cm and 1cm) -- (-7.4,-0.385) arc(-30:30:0.8cm);
\draw[thick,orange] (-7,0) arc(0:30:1cm and 1cm) (-7,0) arc(0:30:1cm and 1cm) -- (-7.4,0.4);
\draw[thick,green] (0,-1) arc(270:240:1cm and 1cm) (0,-1) arc(270:300:1cm and 1cm);
\draw[pattern=north west lines, pattern color=green]  (0,-1) arc(-90:-60:1cm and 1cm) -- (0.385,-0.6) arc(300:240:0.8cm);
\draw[pattern=north west lines, pattern color=green] (0,-1) arc(-90:-120:1cm and 1cm) -- (-0.4,-0.6);
\draw[thick,green] (0,-1) arc(-90:-60:1cm and 1cm) -- (0.385,-0.6) arc(300:240:0.8cm);
\draw[thick,green] (0,-1) arc(-90:-120:1cm and 1cm) -- (-0.4,-0.6);
\draw[thick,orange] (0,-7) arc(90:60:1cm and 1cm) (0,-7) arc(90:120:1cm and 1cm);
\draw[pattern=north west lines, pattern color=orange]  (0,-7) arc(90:60:1cm and 1cm) -- (0.385,-7.4) arc(60:120:0.8cm);
\draw[pattern=north west lines, pattern color=orange] (0,-7) arc(90:120:1cm and 1cm) -- (-0.4,-7.4);
\draw[thick,orange] (0,-7) arc(90:120:1cm and 1cm) -- (-0.4,-7.4);
\draw[thick,orange] (0,-7) arc(90:60:1cm and 1cm) -- (0.385,-7.4) arc(60:120:0.8cm);
\draw[thick,blue] (1,0) arc(0:30:1cm and 1cm) (1,0) arc(0:-30:1cm and 1cm);
\draw[pattern=north west lines, pattern color=blue]  (1,0) arc(0:-30:1cm and 1cm) -- (0.6,-0.385) arc(-30:30:0.8cm);
\draw[pattern=north west lines, pattern color=blue] (1,0) arc(0:30:1cm and 1cm) -- (0.6,0.4);
\draw[thick,blue] (1,0) arc(0:30:1cm and 1cm)  (1,0) arc(0:-30:1cm and 1cm) -- (0.6,-0.385) arc(-30:30:0.8cm);
\draw[thick,blue] (1,0) arc(0:30:1cm and 1cm) (1,0) arc(0:30:1cm and 1cm) -- (0.6,0.4);
\draw[pattern=north west lines, pattern color=orange]  (7,0) arc(180:150:1cm and 1cm) -- (7.4,0.385) arc(150:210:0.8cm);
\draw[pattern=north west lines, pattern color=orange] (7,0) arc(180:210:1cm and 1cm) -- (7.4,-0.4);
\draw[thick,orange] (7,0) arc(180:150:1cm and 1cm)  (7,0) arc(180:150:1cm and 1cm) -- (7.4,0.385) arc(150:210:0.8cm);
\draw[thick,orange] (7,0) arc(180:150:1cm and 1cm) (7,0) arc(180:210:1cm and 1cm) -- (7.4,-0.4);

\fill[pattern=north west lines, pattern color=red] (-9,7.5)rectangle(-8,7.9);
\draw[red,thick](-9,7.5)rectangle(-8,7.9) (-7.9,7.7) node[right,black] {$\overline{\mathscr H}(\mathbf 1,\mathbf 4)$};
\fill[pattern=north west lines, pattern color=blue] (-9,7)rectangle(-8,7.4);
\draw[blue,thick](-9,7)rectangle(-8,7.4) (-7.9,7.2) node[right,black] {$\overline{\mathscr H}(\mathbf 1,\mathbf 2)$};
\fill[pattern=north west lines, pattern color=green] (-9,6.5)rectangle(-8,6.9);
\draw[green,thick](-9,6.5)rectangle(-8,6.9) (-7.9,6.7) node[right,black] {$\overline{\mathscr H}(\mathbf 1,\mathbf 5)$};
\fill[pattern=north west lines, pattern color=yellow] (-9,6)rectangle(-8,6.4);
\draw[yellow,thick](-9,6)rectangle(-8,6.4) (-7.9,6.2) node[right,black] {$\overline{\mathscr H}(\mathbf 1,\mathbf 3)$};
\fill[pattern=vertical lines, pattern color=black] (-9,5.5) rectangle (-8,5.9) (-7.9,5.7) node[right,black] {$W_j(\mathbf 1, \cdot),\ j=2,\dots,L-3$};
\end{tikzpicture}%
    }\par
     \end{minipage}
     \caption{\label{figq5alto} Example of $5-$Potts model with $S=\{1,2,3,4,5\}$. Viewpoint from above on the set of minimal gates around the stable configuration $\mathbf 1$ at energy $2K+2+H(\mathbf 1)$. For any $\mathbf s\in\{\mathbf 2,\mathbf 3,\mathbf 4, \mathbf 5\}$, starting from $\mathbf 1$, the process hits $\mathcal X^s\backslash\{\mathbf 1\}$ for the first time in $\mathbf s$ with probability $\frac{1}{q-1}=\frac{1}{4}$.}
\end{figure}

\begin{theorem}[Minimal gates for the transition $\mathbf r\to\mathcal X^s\backslash\{\mathbf r\}$]\label{mingatesNOcond}
Consider $\mathbf r\in\mathcal{X}^s$. Then, the following sets are minimal gates for the transition $\mathbf r\to\mathcal X^s\backslash\{\mathbf r\}$:
\begin{itemize}
\item[\emph{(a)}]  $\bigcup_{\mathbf t\in\mathcal X^s\backslash \{\mathbf r\}}\overline{\mathscr P}(\mathbf r,\mathbf t)$ and $\bigcup_{\mathbf t\in\mathcal X^s\backslash \{\mathbf r\}}\widetilde{\mathscr P}(\mathbf r,\mathbf t)$;
\item[\emph{(b)}] $\bigcup_{\mathbf t\in\mathcal X^s\backslash \{\mathbf r\}}\overline{\mathcal Q}(\mathbf r,\mathbf t)$ and $\bigcup_{\mathbf t\in\mathcal X^s\backslash \{\mathbf r\}}\widetilde{\mathcal Q}(\mathbf r,\mathbf t)$;
\item[\emph{(c)}] $\bigcup_{\mathbf t\in\mathcal X^s\backslash \{\mathbf r\}}\overline{\mathscr H}_i(\mathbf r,\mathbf t)$ and $\bigcup_{\mathbf t\in\mathcal X^s\backslash \{\mathbf r\}}\widetilde{\mathscr H}_i(\mathbf r,\mathbf t)$ for any $i=1,\dots,K-3$;
\item[\emph{(d)}] $\bigcup_{\mathbf t\in\mathcal X^s\backslash \{\mathbf r\}}\mathcal W_j^{(h)}(\mathbf r,\mathbf t)$ for any $j=2,\dots,L-3$ and any $h=1,\dots,K-1$.
\end{itemize}
\end{theorem}
\begin{theorem}[Union of all minimal gates for the transition $\mathbf r\to\mathcal X^s\backslash\{\mathbf r\}$]\label{setmingatesNOcond}
Given $\mathbf r\in\mathcal{X}^s$, the union of all minimal gates for the transition $\mathbf r\to\mathcal X^s\backslash\{\mathbf r\}$ is given by
\begin{align}\label{mingatestransitiontype2}
\mathcal G(\mathbf r,\mathcal X^s\backslash\{\mathbf r\})=\bigcup_{\mathbf t\in\mathcal X^s\backslash \{\mathbf r\}}\mathcal F(\mathbf r,\mathbf t),
\end{align}
where
\begin{align}
\mathcal F(\mathbf r,\mathbf t)=\bigcup_{j=2}^{L-3}\mathcal W_j(\mathbf r,\mathbf t)\cup\overline{\mathscr H}(\mathbf r,\mathbf t)\cup \widetilde{\mathscr H}(\mathbf r,\mathbf t)\cup\overline{\mathcal Q}(\mathbf r,\mathbf t)\cup\widetilde{\mathcal Q}(\mathbf r,\mathbf t)\cup\overline{\mathscr P}(\mathbf r,\mathbf t)\cup\widetilde{\mathscr P}(\mathbf r,\mathbf t).
\end{align}
\end{theorem}
\begin{remark}\label{remarkkim}
Note that when $A=\{r\}$ and $B=S\backslash\{r\}$ the set of model-dependent gateway configurations given in \cite[Definition 8.1]{kim2021metastability} contains strictly $\mathcal G(\mathbf r,\mathcal X^s\backslash\{\mathbf r\})$, thus it is a gate but it is not minimal.
\end{remark}
We refer to Figure \ref{ciclieplateau} for an illustration of the energy landscape between two Potts stable states.
Moreover, in Figure \eqref{figq5alto}, we depict an example of restricted-gates for $5$-state Potts model in which
the set of minimal restricted-gate corresponds to one of the arms that collegues two different stable states.
\FloatBarrier
\begin{cor}[Crossing the gate]\label{corgateNOcond}
Consider any $\mathbf r\in\mathcal{X}^s$ and the transition from $\mathbf r$ to $\mathcal X^s\backslash\{\mathbf r\}$. Then, the following properties hold:
\begin{itemize}
\item[\emph{(a)}] $\lim_{\beta\to\infty} \mathbb P_\beta(\tau^{\mathbf r}_{\bigcup_{\mathbf t\in\mathcal X^s\backslash \{\mathbf r\}}\overline{\mathscr P}(\mathbf r,\mathbf t)}<\tau^{\mathbf r}_{\mathcal X^s\backslash\{\mathbf r\}})\hspace{-.5mm}=\hspace{-.5mm}\lim_{\beta\to\infty} \mathbb P_\beta(\tau^{\mathbf r}_{\bigcup_{\mathbf t\in\mathcal X^s\backslash \{\mathbf r\}}\widetilde{\mathscr P}(\mathbf r,\mathbf t)}<\tau^{\mathbf r}_{\mathcal X^s\backslash\{\mathbf r\}})=1$;
\item[\emph{(b)}] $\lim_{\beta\to\infty} \mathbb P_\beta(\tau^{\mathbf r}_{\bigcup_{\mathbf t\in\mathcal X^s\backslash \{\mathbf r\}}\overline{\mathcal Q}(\mathbf r,\mathbf t)}<\tau^{\mathbf r}_{\mathcal X^s\backslash\{\mathbf r\}})\hspace{-.5mm}=\hspace{-.5mm}\lim_{\beta\to\infty} \mathbb P_\beta(\tau^{\mathbf r}_{\bigcup_{\mathbf t\in\mathcal X^s\backslash \{\mathbf r\}}\widetilde{\mathcal Q}(\mathbf r,\mathbf t)}<\tau^{\mathbf r}_{\mathcal X^s\backslash\{\mathbf r\}})=1$;
\item[\emph{(c)}] $\lim_{\beta\to\infty} \mathbb P_\beta(\tau^{\mathbf r}_{\bigcup_{\mathbf t\in\mathcal X^s\backslash \{\mathbf r\}}\overline{\mathscr H}_i(\mathbf r,\mathbf t)}\hspace{-.5mm}<\hspace{-.5mm}\tau^{\mathbf r}_{\mathcal X^s\backslash\{\mathbf r\}})\hspace{-.5mm}=\hspace{-.5mm}\lim_{\beta\to\infty} \mathbb P_\beta(\tau^{\mathbf r}_{\bigcup_{\mathbf t\in\mathcal X^s\backslash \{\mathbf r\}}\widetilde{\mathscr H}_i(\mathbf r,\mathbf t)}\hspace{-.5mm}<\hspace{-.5mm}\tau^{\mathbf r}_{\mathcal X^s\backslash\{\mathbf r\}})=1$ for any $i=1,\dots,K-3$;
 \item[\emph{(d)}] $\lim_{\beta\to\infty} \mathbb P_\beta(\tau^{\mathbf r}_{\bigcup_{\mathbf t\in\mathcal X^s\backslash \{\mathbf r\}}\mathcal W_j^{(h)}(\mathbf r,\mathbf t)}<\tau^{\mathbf r}_{\mathcal X^s\backslash\{\mathbf r\}})=1$ for every $j=2,\dots,L-3,$ and $h\hspace{-.5mm}=\hspace{-.5mm}1,\dots,K-1$.
\end{itemize}
Moreover, \emph{(a)--(d)} imply
 \begin{align}
\lim_{\beta\to\infty}\hspace{-0.7mm}\mathbb P_\beta(\tau^{\mathbf r}_{\mathcal G(\mathbf r,\mathcal X^s\backslash\{\mathbf r\})}\hspace{-0.8mm}<\tau^{\mathbf r}_{\mathcal X^s\backslash\{\mathbf r\}})\hspace{-0.5mm}=1.
 \end{align}
\end{cor}
The above corollary implies that every geometrical gate and their union have to be crossed with probability tending to one in the asymptotic limit. In \cite{kim2021metastability}, the authors prove Corollary 8.9 that is similar to Corollary \ref{corgateNOcond}(d).

\subsubsection{Minimal gates for the transition from a stable state to another stable state}\label{subgatetarget}
Finally, we describe the union of the minimal gates for the transition from a stable configuration to another stable state. 
In Theorem \ref{mingatessingh0} we identify geometrically all the sets of minimal gates for the transition from a stable configuration to the another stable state. While in Theorem \ref{setmingatefin} we fully geometrically describe the union of all minimal gates for the same transition.
We assume $q>2$, otherwise when $q=2$, $|\mathcal X^s|=2$ and Theorems \ref{mingatessingh0}--\ref{setmingatefin} coincide with Theorems \ref{mingatescond}--\ref{mingatescondset}. We invite the reader to see Figure \ref{figq5alto} for a pictorial illustration of the structure of the minimal gates.
\begin{theorem}[Minimal gates for the transition $\mathbf r\to\mathbf s$]\label{mingatessingh0}
Consider $\mathbf r,\mathbf{s}\in\mathcal{X}^s$, $\mathbf r\neq \mathbf s$. Then, the following sets are minimal gates for the transition $\mathbf r\to\mathbf s$:
\begin{itemize}
\item[\emph{(a)}] $\bigcup_{\mathbf t\in\mathcal X^s\backslash\{\mathbf r\}} \overline{\mathscr P}(\mathbf r,\mathbf t)$, $\bigcup_{\mathbf t\in\mathcal X^s\backslash\{\mathbf s\}} \overline{\mathscr P}(\mathbf t,\mathbf s)$,\ $\bigcup_{\mathbf t\in\mathcal X^s\backslash\{\mathbf r\}} \widetilde{\mathscr P}(\mathbf r,\mathbf t)$, $\bigcup_{\mathbf t\in\mathcal X^s\backslash\{\mathbf s\}} \widetilde{\mathscr P}(\mathbf t,\mathbf s)$;
\item[\emph{(b)}] $\bigcup_{\mathbf t\in\mathcal X^s\backslash\{\mathbf r\}} \overline{\mathcal Q}(\mathbf r,\mathbf t)$, $\bigcup_{\mathbf t\in\mathcal X^s\backslash\{\mathbf s\}} \overline{\mathcal Q}(\mathbf t,\mathbf s)$,\ $\bigcup_{\mathbf t\in\mathcal X^s\backslash\{\mathbf r\}} \widetilde{\mathcal Q}(\mathbf r,\mathbf t)$, $\bigcup_{\mathbf t\in\mathcal X^s\backslash\{\mathbf s\}} \widetilde{\mathcal Q}(\mathbf t,\mathbf s)$;
\item[\emph{(c)}]  $\bigcup_{\mathbf t\in\mathcal X^s\backslash\{\mathbf r\}} \overline{\mathscr H}_i(\mathbf r,\mathbf t)$, $\bigcup_{\mathbf t\in\mathcal X^s\backslash\{\mathbf s\}} \overline{\mathscr H}_i(\mathbf t,\mathbf s)$,\ $\bigcup_{\mathbf t\in\mathcal X^s\backslash\{\mathbf r\}} \widetilde{\mathscr H}_i(\mathbf r,\mathbf t)$, $\bigcup_{\mathbf t\in\mathcal X^s\backslash\{\mathbf s\}} \widetilde{\mathscr H}_i(\mathbf t,\mathbf s)$ for any $i=1,\dots,K-3$;
\item[\emph{(d)}] $\bigcup_{\mathbf t\in\mathcal X^s\backslash\{\mathbf r\}} \mathcal W_j^{(h)}(\mathbf r,\mathbf t)$,\ $\bigcup_{\mathbf t\in\mathcal X^s\backslash\{\mathbf s\}}\mathcal W_j^{(h)}(\mathbf t,\mathbf s)$ for any $j=2,\dots,L-3$ and any $h=1,\dots,K-1$.
\end{itemize}
\end{theorem}
\begin{theorem}[Union of all minimal gates for the transition $\mathbf r\to\mathbf s$]\label{setmingatefin}
Consider $\mathbf r,\mathbf{s}\in\mathcal{X}^s$, $\mathbf r\neq \mathbf s$. Then, the union of all minimal gates for the transition $\mathbf r\to\mathbf s$ is given by
\begin{align}
\mathcal G(\mathbf r,\mathbf s)=\bigcup_{\mathbf t\in\mathcal X^s\backslash\{\mathbf r\}} \mathcal F(\mathbf r,\mathbf t)\cup\bigcup_{\mathbf t\in\mathcal X^s\backslash\{\mathbf s\}} {\mathcal F}(\mathbf t,\mathbf s),
\end{align}
where, for any $\mathbf t,\mathbf z\in\mathcal X^s$, $\mathbf t\neq\mathbf z$,
\begin{align}
\mathcal F(\mathbf t,\mathbf z)=\bigcup_{j=2}^{L-3} \mathcal W_j(\mathbf t,\mathbf z)\cup\overline{\mathscr H}(\mathbf t,\mathbf z)\cup \widetilde{\mathscr H}(\mathbf t,\mathbf z)\cup\overline{\mathcal Q}(\mathbf t,\mathbf z)\cup\widetilde{\mathcal Q}(\mathbf t,\mathbf z)\cup\overline{\mathscr P}(\mathbf t,\mathbf z)\cup \widetilde{\mathscr P}(\mathbf t,\mathbf z).
\end{align}
\end{theorem}
\begin{cor}[Crossing the gates]\label{corgateconv2}
Consider any $\mathbf r,\mathbf{s}\in\mathcal{X}^s$, with $\mathbf r\neq \mathbf s$, and the transition from $\mathbf r$ to $\mathbf s$. Then, the following properties hold:
\begin{itemize}
\item[\emph{(a)}] $\lim_{\beta\to\infty} \mathbb P_\beta(\tau^{\mathbf r}_{\bigcup_{\mathbf t\in\mathcal X^s\backslash\{\mathbf r\}} \overline{\mathscr P}(\mathbf r,\mathbf t)}<\tau^{\mathbf r}_{\mathbf s})=1$, $\lim_{\beta\to\infty} \mathbb P_\beta(\tau^{\mathbf r}_{\bigcup_{\mathbf t\in\mathcal X^s\backslash\{\mathbf s\}} \overline{\mathscr P}(\mathbf t,\mathbf s)}<\tau^{\mathbf r}_{\mathbf s})\hspace{-.5mm}=\hspace{-.5mm}1$\ and \ 
$\lim_{\beta\to\infty} \mathbb P_\beta(\tau^{\mathbf r}_{\bigcup_{\mathbf t\in\mathcal X^s\backslash\{\mathbf r\}} \widetilde{\mathscr P}(\mathbf r,\mathbf t)}<\tau^{\mathbf r}_{\mathbf s})=1$, $\lim_{\beta\to\infty} \mathbb P_\beta(\tau^{\mathbf r}_{\bigcup_{\mathbf t\in\mathcal X^s\backslash\{\mathbf s\}} \widetilde{\mathscr P}(\mathbf t,\mathbf s)}<\tau^{\mathbf r}_{\mathbf s})=1$;
\item[\emph{(b)}] $\lim_{\beta\to\infty} \mathbb P_\beta(\tau^{\mathbf r}_{\bigcup_{\mathbf t\in\mathcal X^s\backslash\{\mathbf r\}} \overline{\mathcal Q}(\mathbf r,\mathbf t)}<\tau^{\mathbf r}_{\mathbf s})=1$, $\lim_{\beta\to\infty} \mathbb P_\beta(\tau^{\mathbf r}_{\bigcup_{\mathbf t\in\mathcal X^s\backslash\{\mathbf s\}} \overline{\mathcal Q}(\mathbf t,\mathbf s)}<\tau^{\mathbf r}_{\mathbf s})\hspace{-.5mm}=\hspace{-.5mm}1$\ and \

$\lim_{\beta\to\infty} \mathbb P_\beta(\tau^{\mathbf r}_{\bigcup_{\mathbf t\in\mathcal X^s\backslash\{\mathbf r\}} \widetilde{\mathcal Q}(\mathbf r,\mathbf t)}<\tau^{\mathbf r}_{\mathbf s})=1$, $\lim_{\beta\to\infty} \mathbb P_\beta(\tau^{\mathbf r}_{\bigcup_{\mathbf t\in\mathcal X^s\backslash\{\mathbf s\}} \widetilde{\mathcal Q}(\mathbf t,\mathbf s)}<\tau^{\mathbf r}_{\mathbf s})=1$;
\item[\emph{(c)}] $\lim_{\beta\to\infty} \mathbb P_\beta(\tau^{\mathbf r}_{\bigcup_{\mathbf t\in\mathcal X^s\backslash\{\mathbf r\}} \overline{\mathscr H}_i(\mathbf r,\mathbf t)}<\tau^{\mathbf r}_{\mathbf s})=1$, $\lim_{\beta\to\infty} \mathbb P_\beta(\tau^{\mathbf r}_{\bigcup_{\mathbf t\in\mathcal X^s\backslash\{\mathbf s\}} \overline{\mathscr H}_i(\mathbf t,\mathbf s)}<\tau^{\mathbf r}_{\mathbf s})\hspace{-.5mm}=\hspace{-.5mm}1$ and  
$\lim_{\beta\to\infty} \mathbb P_\beta(\tau^{\mathbf r}_{\bigcup_{\mathbf t\in\mathcal X^s\backslash\{\mathbf r\}} \widetilde{\mathscr H}_i(\mathbf r,\mathbf t)}<\tau^{\mathbf r}_{\mathbf s})=1$, $\lim_{\beta\to\infty} \mathbb P_\beta(\tau^{\mathbf r}_{\bigcup_{\mathbf t\in\mathcal X^s\backslash\{\mathbf s\}} \widetilde{\mathscr H}_i(\mathbf t,\mathbf s)}<\tau^{\mathbf r}_{\mathbf s})=1$ for any $i=1,\dots,K-3$;
 \item[\emph{(d)}] for any  $j=2,\dots,L-3$, $h=1,\dots,K-1$, $\lim_{\beta\to\infty} \mathbb P_\beta(\tau^{\mathbf r}_{\bigcup_{\mathbf t\in\mathcal X^s\backslash\{\mathbf r\}} \mathcal W_j^{(h)}(\mathbf r,\mathbf t)}<\tau^{\mathbf r}_{\mathbf s})=1$\ and\ $\lim_{\beta\to\infty} \mathbb P_\beta(\tau^{\mathbf r}_{\bigcup_{\mathbf t\in\mathcal X^s\backslash\{\mathbf s\}}\mathcal W_j^{(h)}(\mathbf t,\mathbf s)}<\tau^{\mathbf r}_{\mathbf s})=1$.
\end{itemize}
 Moreover, \emph{(a)--(d)} imply
 \begin{align}
\lim_{\beta\to\infty}\hspace{-0.7mm}\mathbb P_\beta(\tau^{\mathbf r}_{\mathcal G(\mathbf r,\mathbf s)}\hspace{-0.8mm}<\tau^{\mathbf r}_{\mathbf s})\hspace{-0.5mm}=1.
 \end{align}
\end{cor}
This corollary implies that every geometrical gate for the transition $\mathbf r\to\mathbf s$ and their union have to be crossed with probability tending to one in the asymptotic limit.
\subsubsection{Minimal gates of the Ising model with zero external magnetic field}
When $q=2$, the Potts model corresponds to the Ising model with no external magnetic field, in which $S=\{-1,+1\}$ and $\mathcal X^s=\{\mathbf{-1}, \mathbf{+1}\}$. In this scenario, starting from $\mathbf{-1}$, the target is necessarily $\mathbf{+1}$ and the following corollary holds.
 
\begin{cor}[Union of all minimal gates for the Ising Model with zero external magnetic field]\label{corgate}
Consider the Ising model on a $K\times L$ grid $\Lambda$ with periodic boundary conditions and with zero external magnetic field. Then,  the following sets are minimal gates for the transition $\mathbf{-1}\to\mathbf{+1}$:
\begin{itemize}
\item[\emph{(a)}]  $\overline{\mathscr P}(\mathbf{-1},\mathbf{+1})$ and $\widetilde{\mathscr P}(\mathbf{-1},\mathbf{+1})$;
\item[\emph{(b)}] $\overline{\mathcal Q}(\mathbf{-1},\mathbf{+1})$ and $\widetilde{\mathcal Q}(\mathbf{-1},\mathbf{+1})$;
\item[\emph{(c)}] $\overline{\mathscr H}_i(\mathbf{-1},\mathbf{+1})$ and $\widetilde{\mathscr H}_i(\mathbf{-1},\mathbf{+1})$ for any $i=1,\dots,K-3$;
\item[\emph{(d)}] $\mathcal W_j^{(h)}(\mathbf{-1},\mathbf{+1})$ for any $j=2,\dots,L-3$ and any $h=1,\dots,K-1$.
\end{itemize}
Moreover
\begin{align}
\mathcal G(\mathbf{-1},\mathbf{+1})=\bigcup_{j=2}^{L-3}& \mathcal W_j(\mathbf{-1},\mathbf{+1})\cup\overline{\mathscr H}(\mathbf{-1},\mathbf{+1})\cup \widetilde{\mathscr H}(\mathbf{-1},\mathbf{+1})\cup\overline{\mathcal Q}(\mathbf{-1},\mathbf{+1})\notag\\&\cup\widetilde{\mathcal Q}(\mathbf{-1},\mathbf{+1})\cup\overline{\mathscr P}(\mathbf{-1},\mathbf{+1})\cup \widetilde{\mathscr P}(\mathbf{-1},\mathbf{+1}).
\end{align}
\end{cor}
\section{Restricted-tube of typical paths and tube of typical paths}\label{tube}
At the ending of this section we state our main results on the restricted-tube of typical paths and on the tube of typical paths for the transition from a Potts stable state to the other stable configurations and from a stable state to another one. 

\subsection{Definitions and notations}
In order to describe the tube of typical trajectories performing the transition between two Potts stable configurations, we recall the definitions given in Section \ref{defnotgates} and we add some new ones. 

\subsubsection{Model-independent definitions and notations}\label{deftubemodind}
In addition to the list of Section \ref{modinddefnotgates}, we give also the following model-independent definitions. In particular, these definitions are taken from \cite{nardi2016hitting}, \cite{cirillo2015metastability} and \cite{olivieri2005large}.
\begin{itemize}
\item[-] Given a non-empty subset $\mathcal A\subseteq\mathcal X$, it is said to be \textit{connected} if for any $\sigma,\eta\in\mathcal A$ there exists a path $\omega:\sigma\to\eta$ totally contained in $\mathcal A$. Moreover, we define $\partial A$ as the \textit{external boundary} of $\mathcal A$, i.e., the set
\begin{align}\label{boundaryset}
\partial\mathcal A:=\{\eta\not\in\mathcal A:\ P(\sigma,\eta)>0\ \text{for some}\ \sigma\in\mathcal A\}.
\end{align}
\item[-] A maximal connected set of equal energy states is called a \textit{plateau}.
\item[-] A non-empty subset $\mathcal C\subset\mathcal{X}$ is called \textit{cycle} if it is either a singleton or a connected set such that
\begin{align}\label{cycle}
\max_{\sigma\in\mathcal C} H(\sigma)<H(\mathscr{F}(\partial\mathcal C)).
\end{align}
When $\mathcal C$ is a singleton, it is said to be a \textit{trivial cycle}. We define \textit{extended cycle} a collection of connected equielevated cycles, i.e., cycles of equal energy which belong to the same plateau. It is easy to see that an example of extended cycle is a plateau that may be depicted as union of equielevated trivial cycles. 
Let $\mathscr C(\mathcal{X})$ be the set of cycles and extended cycles of $\mathcal{X}$. 
\item[-] For any $\mathcal C\in\mathscr C(\mathcal X)$, we define external boundary $\partial\mathcal C$ by \eqref{boundaryset}, i.e., as the set can be reached from $\mathcal C$ in one step of the dynamics.
\item[-] For any $\mathcal C\in\mathscr C(\mathcal X)$, we define 
\begin{align}\label{principalboundary}
\mathcal B(\mathcal C):=
\begin{cases}
\mathscr F(\partial\mathcal C) &\text{if}\ \mathcal C\ \text{is non-trivial cycle},\\
\{\eta\in\partial\mathcal C: H(\eta)< H(\sigma)\} &\text{if}\ \mathcal C=\{\sigma\}\ \text{is trivial cycle},\\
\{\eta\in\partial\mathcal C: \exists\sigma\in\mathcal C\ \text{s.t.}\ H(\eta)< H(\sigma)\} &\text{if}\ \mathcal C\ \text{is  extended cycle},
\end{cases}
\end{align}
as the \textit{principal boundary} of $\mathcal C$. Furthermore, let $\partial^{np}\mathcal C$ be the \textit{non-principal boundary} of $\mathcal C$, i.e., $\partial^{np}\mathcal C:=\partial\mathcal C\backslash\mathcal B(\mathcal C).$
\item[-] Given a non-empty set $\mathcal{A}$ and $\sigma\in\mathcal{X}$, we define the \textit{initial cycle} $\mathcal C_{\mathcal{A}}(\sigma)$ by 
\begin{align}\label{initialcycle}
\mathcal C_{\mathcal{A}}(\sigma):=\{\sigma\}\cup\{\eta\in\mathcal{X}:\ \Phi(\sigma,\eta)<\Phi(\sigma,\mathcal A)\}.
\end{align}
If $\sigma\in\mathcal A$, then $\mathcal C_{\mathcal{A}}(\sigma)=\{\sigma\}$ and thus is a trivial cycle. Otherwise, $\mathcal C_{\mathcal{A}}(\sigma)$ is either a trivial cycle (when $\Phi(\sigma, \mathcal{A})=H(\sigma))$ or a non-trivial cycle containing $\sigma$, if  $\Phi(\sigma, \mathcal{A})>H(\sigma)$. In any case, if $\sigma\notin\mathcal{A}$, then $\mathcal C_{\mathcal{A}}(\sigma)\cap {\mathcal{A}}=\varnothing$.
\item[-] The \textit{relevant cycle} $\mathcal C_{\mathcal{A}}^+(\sigma)$ is 
\begin{align}\label{relevantcycle}
\mathcal C_{\mathcal{A}}^+(\sigma):=\{\eta\in\mathcal{X}:\ \Phi(\sigma,\eta)<\Phi(\sigma,\mathcal A)+\delta/2\},
\end{align} 
where $\delta$ is the minimum energy gap between an optimal and a non-optimal path from $\sigma$ to $\mathcal A$.
\item[-] Given a non-empty set $\mathcal A\subset\mathcal X$, we denote by $\mathcal M(\mathcal A)$ the \textit{collection of maximal cycles and extended cycles that partitions} $\mathcal A$, i.e., 
\begin{align}
\mathcal M(\mathcal A):=\{\mathcal C\in\mathscr C(\mathcal X)|\ \mathcal C\ \text{maximal by inclusion under constraint}\ \mathcal C\subseteq\mathcal A\}.\notag
\end{align}
\item[-] We call \textit{cycle-path} a finite sequence $(\mathcal C_1,\dots,\mathcal C_m)$ of trivial, non-trivial and extended cycles $\mathcal C_1,\dots,\mathcal C_m\in\mathscr C(\mathcal X)$, such that 
$\mathcal C_i\cap\mathcal C_{i+1}=\varnothing$ and $\partial\mathcal C_i\cap\mathcal C_{i+1}\neq\varnothing$, for every $i=1,\dots,m-1$. We denote the set of cycle-paths that lead from $\sigma$ to $\mathcal A$ and consist of maximal cycles  in $\mathcal X\backslash\mathcal A$ by 
\begin{align}
\hspace{0pt}\mathcal {P}_{\sigma,\mathcal A}:=\{\text{cycle-path}\ (\mathcal C_1,\dots,\mathcal C_m)|\ \mathcal C_1,\dots,\mathcal C_m&\in\mathcal M(\mathcal C_{\mathcal{A}}^+(\sigma)\backslash A),\sigma\in\mathcal C_1, \partial\mathcal C_m\cap\mathcal A\neq\varnothing\}.\notag
\end{align}
\item[-]  Given a non-empty set $\mathcal{A}\subset\mathcal X$ and $\sigma\in\mathcal{X}$, we constructively define a mapping $G: \Omega_{\sigma,A}\to\mathcal P_{\sigma,\mathcal A}$. More precisely, given $\omega=(\omega_1,\dots,\omega_n)\in\Omega_{\sigma,A}$, we set $m_0=1$, $\mathcal C_1=\mathcal C_{\mathcal A}(\sigma)$ and define recursively $m_i:=\min\{k>m_{i-1}|\ \omega_k\notin\mathcal C_i\}$ and $\mathcal C_{i+1}:=\mathcal C_{\mathcal A}(\omega_{m_i})$. We note that $\omega$ is a finite sequence and $\omega_n\in\mathcal A$, so there exists an index $n(\omega)\in\mathbb N$ such that $\omega_{m_{n(\omega)}}=\omega_n\in\mathcal A$ and there the procedure stops. The way the sequence $(\mathcal C_1,\dots,\mathcal C_{m_{n(\omega)}})$ is constructed shows that it is a cycle-path with $\mathcal C_1,\dots,\mathcal C_{m_{n(\omega)}}\subset\mathcal M(\mathcal X\backslash\mathcal A)$. Moreover, the fact that $ \omega\in\Omega_{\sigma,A}$ implies that $\sigma\in\mathcal C_1$ and that $\partial\mathcal C_{n(\omega)}\cap\mathcal A\neq\varnothing$, hence $G(\omega)\in\mathcal P_{\sigma,\mathcal A}$ and the mapping is well-defined.

\item[-] We say that a cycle-path $(\mathcal C_1,\dots,\mathcal C_m)$ is \textit{connected via typical jumps} to $\mathcal A\subset\mathcal X$ or simply $vtj-$\textit{connected} to $\mathcal A$ if
\begin{align}\label{cyclepathvtj}
\mathcal B(\mathcal C_i)\cap\mathcal C_{i+1}\neq\varnothing,\ \ \forall i=1,\dots,m-1,\ \ \text{and}\ \ \mathcal B(\mathcal C_m)\cap\mathcal A\neq\varnothing.\end{align}
Let $J_{\mathcal C,\mathcal A}$ be the collection of all cycle-path $(\mathcal C_1,\dots,\mathcal C_m)$ vtj-connected to $\mathcal A$ such that $\mathcal C_1=\mathcal C$. 
\item[-] Given a non-empty set $\mathcal{A}$ and $\sigma\in\mathcal{X}$, we define $\omega\in\Omega_{\sigma,A}$ as a \textit{typical path} from $\sigma$ to $\mathcal A$ if its corresponding cycle-path $G(\omega)$ is vtj-connected to $\mathcal A$ and we denote by $\Omega_{\sigma,A}^{\text{vtj}}$ the collection of all typical paths from $\sigma$ to $\mathcal A$, i.e., 
\begin{align}
\Omega_{\sigma,A}^{\text{vtj}}:=\{\omega\in\Omega_{\sigma,A}|\ G(\omega)\in J_{\mathcal C_{\mathcal A}(\sigma),\mathcal A}\}.
\end{align}
It is useful to recall the following \cite[Lemma 3.12]{nardi2016hitting} in which the authors give the following equivalent characterization of a typical path from $\sigma\notin\mathcal A$ and $\mathcal A\subset\mathcal X$.
\begin{lemma}
Consider a non empty subset $\mathcal A\subset\mathcal X$ and $\sigma\notin\mathcal A$. Then
\begin{align}\label{chartyppath}
\omega\in\Omega_{\sigma,A}^{\text{vtj}}\iff\omega\in\Omega_{\sigma,A}\ \text{and}\ \Phi(\omega_{i+1},\mathcal A)\le\Phi(\omega_i,\mathcal A)\ \forall i=1,\dots,|\omega|.
\end{align}\end{lemma}
\item[-] We define the \textit{tube of typical paths} from $\sigma$ to $\mathcal A$, $T_{\mathcal A}(\sigma)$, as the subset of states $\eta\in\mathcal X$ that can be reached from $\sigma$ by means of a typical path which does not enter $\mathcal A$ before visiting $\eta$, i.e.,
\begin{align}
T_{\mathcal A}(\sigma):=\{\eta\in\mathcal X|\ \exists\omega\in\Omega_{\sigma,A}^{\text{vtj}}:\ \eta\in\omega\}.
\end{align}
Moreover, we define $\mathfrak T_{\mathcal A}(\sigma)$ as the set of all maxiaml cycles and maximal extended cycles that belong to at least one vtj-connected path from $\mathcal C_{\mathcal A}(\sigma)$ to $\mathcal A$, i.e.,
\begin{align}\label{tubostorto}
\mathfrak T_{\mathcal A}(\sigma):=\{\mathcal C\in\mathcal M(\mathcal C_{\mathcal A}^+(\sigma)\backslash\mathcal A)|\ \exists(\mathcal C_1,\dots,\mathcal C_n)\in J_{\mathcal C_{\mathcal A}(\sigma),\mathcal A}\ \text{and}\exists j\in\{1,\dots,n\}:\ \mathcal C_j=\mathcal C\}.
\end{align}
Note that $\mathfrak T_{\mathcal A}(\sigma)=\mathcal M(T_{\mathcal A}(\sigma)\backslash\mathcal A)$ and that the boundary of $T_{\mathcal A}(\sigma)$ consists of states either in $\mathcal A$ or in the non-principal part of the boundary of some $\mathcal C\in\mathfrak T_{\mathcal A}(\sigma)$:
\begin{align}
\partial T_{\mathcal A}(\sigma)\backslash\mathcal A\subseteq \bigcup_{\mathcal C\in\mathfrak T_{\mathcal A}(\sigma)} (\partial\mathcal C\backslash\mathcal B(\mathcal C))=:\partial^{np}\mathfrak T_{\mathcal A}(\sigma).
\end{align}
For sake of semplicity, we will also refer to $\mathfrak T_{\mathcal A}(\sigma)$ as tube of typical paths from $\sigma$ to $\mathcal A$.
\item[-] Given $|\mathcal X^s|>2$ and $\sigma,\sigma'\in\mathcal X^s, \sigma\neq\sigma'$, we define the \textit{restricted-tube of typical paths} between $\sigma$ and $\sigma'$, $\mathcal U_{\sigma'}(\sigma)$, as the subset of states $\eta\in\mathcal X$ that can be reached from $\sigma$ by means of a typical path which does not intersect $\mathcal X^s\backslash\{\sigma,\sigma'\}$ and does not visit $\sigma'$ before visiting $\eta$, i.e,
\begin{align}
\mathcal U_{\sigma'}(\sigma):=\{\eta\in\mathcal X|\ \exists\omega\in\Omega_{\sigma,\sigma'}^{\text{vtj}}\ \text{s.t.}\ \omega\cap\mathcal X^s\backslash\{\sigma,\sigma'\}=\varnothing\ \text{and}\ \eta\in\omega\}. 
\end{align}
Moreover, we set $\mathscr U_{\sigma'}(\sigma)$ as the set of maximal cycles and maximal extended cycles that belong to at last one vtj-connected path from $\mathcal C_{\sigma'}(\sigma)$ to $\sigma'$ such that does not intersect $\mathcal X^s\backslash\{\sigma,\sigma'\}$, i.e.,
\begin{align}\label{restrictedtube}
\hspace{-10pt}\mathscr U_{\sigma'}(\sigma):=\{\mathcal C\in&\mathcal M(\mathcal C_{\{\sigma'\}}^+(\sigma)\backslash\{\sigma'\})|\exists(\mathcal C_1,\dots,\mathcal C_m)\in J_{\mathcal C_{\sigma'}(\sigma),\{\sigma'\}}\ \text{such that}\notag \\ &\bigcup_{i=1}^m \mathcal C_i\cap\mathcal X^s\backslash\{\sigma,\sigma'\}=\varnothing\ \text{and}\ \exists j\in\{1,\dots,n\}:\ \mathcal C_j=\mathcal C\}.
\end{align}
Note that $\mathscr U_{\sigma'}(\sigma)=\mathcal M(\mathcal U_{\sigma'}(\sigma)\backslash(\mathcal X^s\backslash\{\sigma\}))$ and that the boundary of $\mathcal U_{\sigma'}(\sigma)$ consists of $\sigma'$ and of states in the non-principal part of the boundary of some $\mathcal C\in\mathscr U_{\sigma'}(\sigma)$:
\begin{align}
\partial\mathcal U_{\sigma'}(\sigma)\backslash\{\sigma'\}\subseteq\bigcup_{\mathcal C\in\mathscr U_{\sigma'}(\sigma)} (\partial\mathcal C\backslash\mathcal B(\mathcal C)) =:\partial^{np}\mathscr U_{\sigma'}(\sigma).
\end{align}
For sake of semplicity, we will also refer to $\mathscr U_{\sigma'}(\sigma)$ as restricted-tube of typical paths from $\sigma$ to $\sigma'$.
\end{itemize}
\begin{remark}
Note that the notion of extended cyles is taken from \cite{olivieri2005large}. In particular, using also the extended cycles for defining a cycle-path vtj-connected, we get that this object is the so-called \textit{standard cascade} in \cite{olivieri2005large}.
\end{remark}

\subsubsection{Model-dependent definitions and notations}\label{moddeptubo}
In this section we add some model-dependent definitions to the list given in Section \ref{moddepdefnot}.\\
For any $s,r\in\{1,\dots,q\}$, let $\mathcal X(r,s)$ be the set defined by
\begin{align}
\mathcal X(r,s):=\{\sigma\in\mathcal X: \sigma(v)\in\{r,s\}\ \forall v\in V\}.
\end{align}
Let $R_{\ell_1\times\ell_2}$ be a rectangle in $\mathbb R^2$ with horizontal side of length $\ell_1$ and vertical side of length $\ell_2$. 
\begin{itemize}
\item[-] We define the set $\text C^s(\sigma)\subseteq\mathbb R^2$ as the union of unit closed squares centered at the vertices $v\in V$ such that $\sigma(v)=s$. We define $s$-\textit{clusters} the maximal connected components $C^s_1,\dots,C^s_m,\ m\in\mathbb N$, of $\text C^s(\sigma)$ and we consider separately two $s$-clusters which share only one point. In particular, two clusters $C^s_1, C^s_2$ of spins are said to be \textit{interacting} if either $C_1^s$ and $C_2^s$ intersect or are disjoint but there exists a site $v\notin C_1^s\cup C_2^s$ such that $\sigma(v)\neq s$ with two distinct nearest-neighbor sites $u, w$ lying inside $C_1^s, C_2^s$ respectively. In particular, we say that a $q$-Potts configuration has $s$-\textit{interacting clusters} when all its $s$-clusters are interacting.
\item[-] We set $R(\text C^s(\sigma))$ as the smallest surrounding rectangle to the union of the clusters of spins $s$ in $\sigma$.
\end{itemize} 

Now we define some extended cycles that are crucial to describe the tube of typical paths.
\begin{itemize}
 \item[-] We define 
\begin{align}\label{defKbar}
\overline{\mathcal K}(r,s)&:=\{\sigma\in\mathcal X(r,s): H(\sigma)=2K+2+H(\mathbf r),\  \sigma\ \text{has either a $s$-cluster or more }\notag\\ &\text{$s$-interacting clusters and}\ R(\text C^s(\sigma))=R_{2\times(K-1)}\}\cup\overline{\mathcal Q}(\mathbf r,\mathbf s)\cup\overline{\mathscr P}(\mathbf r,\mathbf s).
\end{align} 
Note that $\overline{\mathscr H}(\mathbf r,\mathbf s)\subset\overline{\mathcal K}(r,s)$.
\item[-] We set
\begin{align}\label{defD1}
\overline{\mathcal D}_1(r,s)&:=\{\sigma\in\mathcal X(r,s):\ H(\sigma)=2K+H(\mathbf r),\ \sigma\ \text{has either a $s$-cluster }\notag\\& \text{or more  $s$-interacting clusters}\ \text{such that}\ R(\text C^s(\sigma))=R_{2\times(K-2)}\},
\\ &\hspace{-58pt}\text{and}&\notag \\ \label{defE1}
\overline{\mathcal E}_1(r,s)&:=\{\sigma\in\mathcal X(r,s): H(\sigma)=2K+H(\mathbf r),\ \sigma\ \text{has either a $s$-cluster or more  }\notag\\ &\text{$s$-interacting clusters}\ \text{such that}\ R(\text C^s(\sigma))=R_{1\times(K-1)}\}\cup\bar R_{1,K}(r, s).
\end{align}
\item[-] For any $i=2,\dots,K-2$, we define
\begin{align}\label{defDi}
\overline{\mathcal D}_i(r,s):=&\{\sigma\in\mathcal X(r,s):\ H(\sigma)=2K-2i+2+H(\mathbf r),\sigma\ \text{has either a $s$-cluster}\notag\\& \text{ or more $s$-interacting clusters such that}\ R(\text C^s(\sigma))=R_{2\times(K-(i+1))}\},
\\ &\hspace{-54pt}\text{and}&\notag \\
\label{defEi}\overline{\mathcal E}_i(r,s):=&\{\sigma\in\mathcal X(r,s):\ H(\sigma)=2K-2i+2+H(\mathbf r),\ \sigma\ \text{has either a $s$-cluster}\notag\\ &\text{or more $s$-interacting clusters}\text{such that}\  R(\text C^s(\sigma))=R_{1\times(K-i)}\}.
\end{align}
\item[-] Similarly, for any $i=1,\dots,K-2$ we set $\widetilde{\mathcal K}(r,s)$, $\widetilde{\mathcal D}_i(r,s)$, $\widetilde{\mathcal E}_i(r,s)$ by interchanging the role of the spins $r$ and $s$, i.e., they are defined as the collection of those configurations which have either a $r$-cluster or more than one $r$-interacting clusters such that $R(\text C^r(\sigma))$ satifies the same conditions given in \eqref{defKbar}--\eqref{defEi} given for $s$.
\end{itemize}
We refer to Figure \ref{tubecycle} and to Figures \ref{examplesprincbound}--\ref{thirdexamplesprincbound} for an illustration of an example of the extended cycles defined above. 

\subsection{Main results}\label{mainrestube}
We are now ready to state the main results on the restricted-tube of typical paths between two Potts stable configurations and on the tube of typical trajectories from a stable state to the other stable states and from a stable state to another stable configuration. In particular, we prove Theorems \ref{exittuberes}, \ref{exittubestable} and \ref{exittubetarget} using \cite[Lemma 3.13]{nardi2016hitting} and \cite[Proposition 2.7]{cirillo2015metastability}.\\

In order to understand Theorems \ref{exittuberes}, \ref{exittubestable}, \ref{exittubetarget}, one may think about the following analogy. If Figure \ref{figq5alto} represents valleys having actual physical depth, and one imagines pouring a liquid in the valley corresponding to the configuration \textbf1, then the tubes of typical trajectories are all the different ways in which the liquid might flow out of that valley.

\subsubsection{Restricted-tube of typical paths between two Potts stable configurations}\label{mainrestuberes}
We briefly recall that given $\mathbf r,\mathbf s\in\mathcal X^s, \mathbf r\neq\mathbf s$, the restricted-tube of typical paths $\mathcal U_{\mathbf s}(\mathbf r)$ is the set of those configurations belonging to at least a typical path $\omega\in\Omega_{\mathbf r,\mathbf s}^{\text{vtj}}$ such that $\omega\cap\mathcal X^s\backslash\{\mathbf r,\mathbf s\}=\varnothing$. 
Since in absence of external magnetic field the energy landscape between two Potts stable configurations is characterized by many extended-cycles, we describe the restricted-tube of typical paths defined in general in  \eqref{restrictedtube}. For our model, let
\begin{align}\label{aligntuberes}
\mathscr U_{\mathbf s}&(\mathbf r):=\bar R_{1,1}(r,s)\cup\bigcup_{i=1}^{K-2} (\overline{\mathcal D}_i(r,s)\cup\overline{\mathcal E}_i(r,s))\cup\overline{\mathcal K}(r,s)\cup\bigcup_{h=2}^{K-2} \bar B_{1,K-1}^h(r,s)\cup\bigcup_{j=2}^{L-2} \bigcup_{h=1}^{K-1} \bar B_{j,K}^h(r,s)\notag\\&\cup\bigcup_{j=2}^{L-2} \bar R_{j,K}(r,s)\cup\bigcup_{h=2}^{K-2} \tilde B_{1,K-1}^h(r,s)\cup\widetilde{\mathcal K}(r,s)\cup\bigcup_{i=1}^{K-2}(\widetilde{\mathcal D}_i(r,s)\cup\widetilde{\mathcal E}_i(r,s))\cup\tilde R_{1,1}(r,s).
\end{align}

As illustrated in the next result, $\mathscr U_{\mathbf s}(\mathbf r)$ includes those configurations with no-vanishing probability of being visited by the Markov chain $\{X_t\}^\beta_{t\in\mathbb N}$ started in $\mathbf r$ before hitting $\mathbf s$ in the limit $\beta\to\infty$.

\begin{theorem}[Restricted-tube of typical paths]\label{exittuberes}
For every $\mathbf{r}, \mathbf s\in\mathcal{X}^s$, $\mathbf s\neq\mathbf r$, we have that $\mathscr U_{\mathbf s}(\mathbf r)$ is the restricted-tube of typical paths for the transition $\mathbf r\to\mathbf s$.
Moreover, there exists $k>0$ such that for $\beta$ sufficiently large 
\begin{align}\label{timetuberes}
\mathbbm P_\beta(\tau_{\partial^{np}\mathscr U_{\mathbf s}(\mathbf r)}^\mathbf r\le\tau_\mathbf s^\mathbf r|\tau^\mathbf r_{\mathbf s}<\tau^\mathbf r_{\mathcal X^s\backslash\{\mathbf r,\mathbf s\}})\le e^{-k\beta}.
\end{align}
\end{theorem}
\begin{figure}[h!]
\begin{minipage}[c]{0.95\textwidth}
    \centering
    \makebox[0pt]{%
\centering
\begin{tikzpicture}[scale=0.8, transform shape]
\draw[gray] (0,0)--(9.4,0)--(9.4,0.5)--(0,0.5)--(0,0);

\draw[dotted] (-0.6,0)--(10,0) (-0.6,-0.8)--(10,-0.8) (-0.6,-1.6)--(10,-1.6) (-0.6,-2.4)--(10,-2.4) (-0.6,-3.2)--(10,-3.2) (-0.6,-4.8)--(10,-4.8) (-0.6,-5.6)--(10,-5.6) (-0.6,-6.4)--(10,-6.4) (-0.6,-7.2)--(10,-7.2);
\draw (0,0.25) node[right] {\footnotesize{$\overline{\mathcal K}(r,s)$}};

\draw[gray] (0,-0.8)--(4.2,-0.8)--(4.2,-0.3)--(0,-0.3)--(0,-0.8);\draw (0,-0.55) node[right] {\footnotesize{$\overline{\mathcal E}_1(r,s)$}};
\draw[gray] (9.4,-0.8)--(5.2,-0.8)--(5.2,-0.3)--(9.4,-0.3)--(9.4,-0.8); \draw (9.4,-0.55) node[left] {\footnotesize{$\overline{\mathcal D}_1(r,s)$}};
 
\draw[gray] (0.3,-1.6)--(4.2,-1.6)--(4.2,-1.1)--(0.3,-1.1)--(0.3,-1.6);\draw (0.3,-1.35) node[right] {\footnotesize{$\overline{\mathcal E}_2(r,s)$}};
\draw[gray] (9.1,-1.6)--(5.2,-1.6)--(5.2,-1.1)--(9.1,-1.1)--(9.1,-1.6); \draw (9.1,-1.35) node[left] {\footnotesize{$\overline{\mathcal D}_2(r,s)$}};
 
\draw[gray] (0.6,-2.4)--(4.2,-2.4)--(4.2,-1.9)--(0.6,-1.9)--(0.6,-2.4);\draw (0.6,-2.15) node[right] {\footnotesize{$\overline{\mathcal E}_3(r,s)$}};
\draw[gray] (8.8,-2.4)--(5.2,-2.4)--(5.2,-1.9)--(8.8,-1.9)--(8.8,-2.4); \draw (8.8,-2.15) node[left] {\footnotesize{$\overline{\mathcal D}_3(r,s)$}};
 
\draw[gray] (0.9,-3.2)--(4.2,-3.2)--(4.2,-2.7)--(0.9,-2.7)--(0.9,-3.2);\draw (0.9,-2.95) node[right] {\footnotesize{$\overline{\mathcal E}_4(r,s)$}};
\draw[gray] (8.5,-3.2)--(5.2,-3.2)--(5.2,-2.7)--(8.5,-2.7)--(8.5,-3.2); \draw (8.5,-2.95) node[left] {\footnotesize{$\overline{\mathcal D}_4(r,s)$}};
 
 \draw (3.25,-3.75) node[thick] {{$\vdots$}}; \draw(6.15,-3.75) node[thick] {{$\vdots$}};
 
\draw[gray] (1.5,-4.8)--(4.2,-4.8)--(4.2,-4.3)--(1.5,-4.3)--(1.5,-4.8);\draw (1.5,-4.55) node[right] {\footnotesize{$\overline{\mathcal E}_{K-3}(r,s)$}};
\draw[gray] (7.9,-4.8)--(5.2,-4.8)--(5.2,-4.3)--(7.9,-4.3)--(7.9,-4.8); \draw (7.9,-4.55) node[left] {\footnotesize{$\overline{\mathcal D}_{K-3}(r,s)$}}; 

\draw[gray] (1.8,-5.6)--(4.2,-5.6)--(4.2,-5.1)--(1.8,-5.1)--(1.8,-5.6);\draw (1.8,-5.35) node[right] {\footnotesize{$\overline{\mathcal E}_{K-2}(r,s)$}};
\draw[gray] (7.6,-5.6)--(5.2,-5.6)--(5.2,-5.1)--(7.6,-5.1)--(7.6,-5.6); \draw (7.6,-5.35) node[left] {\footnotesize{$\overline{\mathcal D}_{K-2}(r,s)$}};

\draw[gray] (4,-6.4)--(5.4,-6.4)--(5.4,-5.9)--(4,-5.9)--(4,-6.4);\draw (4.7,-6.15) node{\footnotesize{$\bar R_{1,1}(\mathbf r,\mathbf s)$}};

\draw (4.7,-6.95) node{\footnotesize{$\mathbf r$}};

\draw[->,thick,black!60!white] (2.2,0)--(2.2,-0.25);\draw[->,thick,black!90!white] (7.2,0)--(7.2,-0.25);
\draw[->,thick,black!60!white] (2.4,-0.8)--(2.4,-1.05);\draw[->,thick,black!90!white] (7,-0.8)--(7,-1.05);
\draw[->,thick,black!60!white] (2.6,-1.6)--(2.6,-1.85);\draw[->,thick,black!90!white] (6.8,-1.6)--(6.8,-1.85);
\draw[->,thick,black!60!white] (2.8,-2.4)--(2.8,-2.65);\draw[->,thick,black!90!white] (6.6,-2.4)--(6.6,-2.65);
\draw[->,thick,black!60!white] (3,-3.2)--(3,-3.45);\draw[->,thick,black!90!white] (6.4,-3.2)--(6.4,-3.45);
\draw[->,thick,black!60!white] (3.4,-4.8)--(3.4,-5.05);\draw[->,thick,black!90!white] (6,-4.8)--(6,-5.05);
\draw[->,thick,black!60!white] (3.6,-5.6)--(4.1,-5.8);\draw[->,thick,black!90!white] (5.8,-5.6)--(5.3,-5.8);
\draw[->,thick] (4.7,-6.4)-- (4.7,-6.8);

\draw[->,thick,black!60!white] (0,0.25)--(-0.3,0.25)--(-0.3,-1.3)--(0.2,-1.3); \draw[->,thick,black!90!white] (9.4,0.25)--(9.7,0.25)--(9.7,-1.3)--(9.2,-1.3);
\draw[->,thick,black!60!white] (0,-0.6)--(-0.2,-0.6)--(-0.2,-2.1)--(0.5,-2.1);\draw[->,thick,black!90!white] (9.4,-0.6)--(9.6,-0.6)--(9.6,-2.1)--(8.9,-2.1);
\draw[->,thick,black!60!white] (0.3,-1.4)--(0,-1.4)--(0,-2.9)--(0.8,-2.9);\draw[->,thick,black!90!white] (9.1,-1.4)--(9.4,-1.4)--(9.4,-2.9)--(8.6,-2.9);
\draw[->,thick,black!60!white] (0.6,-2.2)--(0.3,-2.2)--(0.3,-3.7)--(1.1,-3.7);\draw[->,thick,black!90!white] (8.8,-2.2)--(9.1,-2.2)--(9.1,-3.7)--(7.9,-3.7);
\draw[->,thick,black!60!white] (1.5,-4.6)--(0.9,-4.6)--(0.9,-6.15)--(3.9,-6.15);
\draw[->,thick,black!60!white] (1.2,-4.4)--(1.2,-5.3)--(1.7,-5.3);\draw[->,thick,black!90!white] (8.2,-4.4)--(8.2,-5.3)--(7.7,-5.3); \draw[dotted,thick,black!60!white] (1.2,-4.4)--(1.2,-4.15); \draw[dotted,thick,black!90!white] (8.2,-4.4)--(8.2,-4.15);
\draw[->,thick,black!60!white] (1.8,-5.4)--(1.2,-5.4)--(1.2,-6.95)--(4.3,-6.95);\draw[->,thick,black!90!white] (7.6,-5.4)--(8.2,-5.4)--(8.2,-6.95)--(5.1,-6.95);

\draw[->,thick,black!90!white] (5.2,-0.7)--(4.7,-1.2)--(4.3,-1.5);\draw[->,thick,black!90!white](4.7,-1.2)--(4.3,-2);
\draw[->,thick,black!90!white] (5.2,-1.5)--(4.7,-2)--(4.3,-2.3);\draw[->,thick,black!90!white](4.7,-2)--(4.3,-2.8);
\draw[->,thick,black!90!white] (5.2,-2.3)--(4.7,-2.8)--(4.3,-3.1);\draw[->,thick,black!90!white](4.7,-2.8)--(4.3,-3.6);
\draw[->,thick,black!90!white] (5.2,-4.7)--(4.7,-5.2)--(4.3,-5.5);\draw[->,thick,black!90!white](4.7,-5.2)--(4.7,-5.8);
\draw[thick,black!90!white] (5.2,-3.1)--(4.7,-3.6);\draw[dotted,thick,black!90!white] (4.7,-3.6)--(4.55,-3.75);

\draw (-0.8,0.2) node[left] {\tiny{$2K+2+H(\mathbf r)$}};
\draw (-0.8,-0.6) node[left] {\tiny{$2K+H(\mathbf r)$}};
\draw (-0.8,-1.4) node[left] {\tiny{$2K-2+H(\mathbf r)$}};
\draw (-0.8,-2.2) node[left] {\tiny{$2K-4+H(\mathbf r)$}};
\draw (-0.8,-3) node[left] {\tiny{$2K-6+H(\mathbf r)$}};
\draw (-1.4,-3.65) node {\tiny{$\vdots$}};
\draw (-0.8,-4.6) node[left] {\tiny{$8+H(\mathbf r)$}};
\draw (-0.8,-5.4) node[left] {\tiny{$6+H(\mathbf r)$}};
\draw (-0.8,-6.2) node[left] {\tiny{$4+H(\mathbf r)$}};
\draw (-0.8,-7) node[left] {\tiny{$H(\mathbf r)$}};
\end{tikzpicture}%
    }\par
     \end{minipage}
\caption{\label{tubecycle} Illustration of the first descent from $\overline{\mathcal K}(r,s)$ to $\mathbf r$. The rectangles denote extended cycles, i.e., the sets of trivial equielevated cycles. The arrows denote the connection between each extended cycle and the sets which belong to its principal boundary.}
\end{figure}
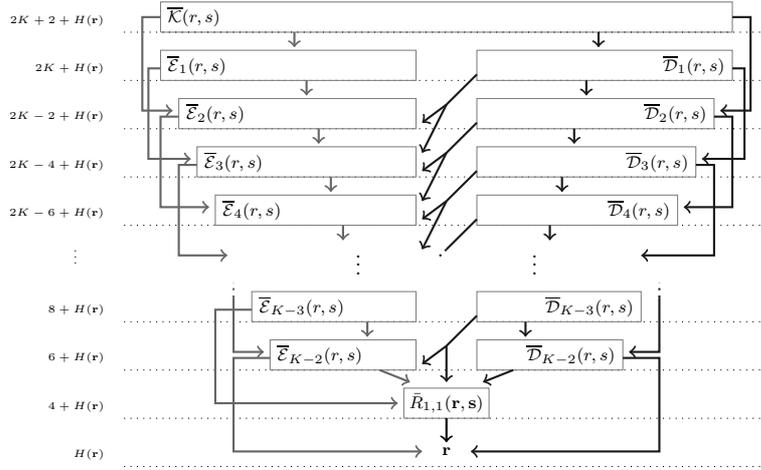\FloatBarrier
\subsubsection{Tube of typical paths between a stable state and the other stable states}\label{mainrestubeset}
Using the results about the restricted-tube of typical paths, we prove the following results on the the tube of typical trajectories from a stable configuration to all the other stable states.  We assume $q>2$, since in the case $q=2$ the Hamiltonian has only two global minima, $|\mathcal X^s|=2$, and the result is given by Theorem \ref{exittuberes}.\\

We recall the model-independent definition \eqref{tubostorto} and 
define $\mathfrak T_{\mathcal X^s\backslash\{\mathbf r\}}$ for Potts model with $q>2$: 
\begin{align}\label{aligntubestable}
\mathfrak T_{\mathcal X^s\backslash\{\mathbf r\}}(\mathbf r):=\bigcup_{\mathbf t\in\mathcal X^s\backslash\{\mathbf r\}}\mathscr U_{\mathbf t}(\mathbf r).
\end{align} 

\begin{theorem}[Tube of typical paths for the transition $\mathbf r\to\mathcal X^s\backslash\{\mathbf r\}$]\label{exittubestable}
For any $\mathbf r\in\mathcal{X}^s$, we have that $\mathfrak T_{\mathcal X^s\backslash\{\mathbf r\}}(\mathbf r)$ is the tube of typical trajectories for the transition $\mathbf r\to\mathcal X^s\backslash\{\mathbf r\}$ and there exists $k>0$ such that for $\beta$ sufficiently large
\begin{align}\label{timetubestable}
\mathbbm P_\beta(\tau_{\partial^{np}\mathfrak T_{\mathcal X^s\backslash\{\mathbf r\}}(\mathbf r)}^\mathbf r\le\tau_{\mathcal X^s\backslash\{\mathbf r\}}^\mathbf r)\le e^{-k\beta}.
\end{align}
\end{theorem}

\subsubsection{Tube of typical paths between a stable state and another stable state}\label{mainrestubetarget} 
Finally, using the results about the tube of typical paths from a stable state to the other stable configurations, we describe the tube of typical paths from a stable configuration to another stable state. Arguing like in Section \ref{subgatetarget}, we describe the typical trajectories for the transition $\mathbf r\to\mathbf s$ in terms of a sequence of transitions between two stable states such that the path followed by the process does not intersect other stable configurations.  We assume $q>2$, otherwise when $q=2$ the Hamiltonian has only two global minima and the results on the minimal gates coincide with Theorem \ref{exittuberes}.\\
We recall the model-independent definition \eqref{tubostorto} and 
define $\mathfrak T_{\mathbf s}(\mathbf r)$ for Potts model with $q>2$: 
\begin{align}\label{aligntubetarget}
\mathfrak T_{\mathbf s}(\mathbf r):=\bigcup_{\mathbf t\in\mathcal X^s\backslash\{\mathbf r\}} \mathscr U_{\mathbf t}(\mathbf r)\cup\bigcup_{\mathbf t,\mathbf t' \in\mathcal X^s\backslash\{\mathbf r,\mathbf s\},\mathbf t\neq\mathbf t'}\mathscr U_{\mathbf t}(\mathbf t')\cup\bigcup_{\mathbf t'\in\mathcal X^s\backslash\{\mathbf s\}}\mathscr U_{\mathbf s}(\mathbf t').
\end{align}

\begin{theorem}[Tube of typical paths for the transition $\mathbf r\to\mathbf s$]\label{exittubetarget}
For any $\mathbf r,\mathbf s\in\mathcal{X}^s$, $\mathbf r\neq \mathbf s$ we have that $\mathfrak T_{\mathbf s}(\mathbf r)$ is the tube of typical trajectories for the transition $\mathbf r\to\mathbf s$ and there exists $k>0$ such that for $\beta$ sufficiently large
\begin{align}\label{timetubetarget}
\mathbbm P_\beta(\tau_{\partial^{np}\mathfrak T_{\mathbf s}(\mathbf r)}^\mathbf r\le\tau_\mathbf s^\mathbf r)\le e^{-k\beta}.
\end{align} 
\end{theorem}
\subsubsection{Tube of typical paths for the Ising model with zero magnetic field}\label{mainresisingtube}
For sake of completeness, we give the following result on the tube of typical paths for the Ising model with zero magnetic field.\\
We recall the tube of typical paths defined in general in \eqref{tubostorto} and 
define $\mathfrak T_{\mathbf{+1}}(\mathbf{-1})$ for Ising model:
\begin{align}&
\mathfrak T_{\mathbf{+1}}(\mathbf{-1}):=\bar R_{1,1}(-1,+1)\cup\bigcup_{i=1}^{K-2} (\overline{\mathcal D}_i(-1,+1)\cup\overline{\mathcal E}_i(-1,+1))\cup\overline{\mathcal K}(-1,+1)\notag \\ 
&\cup\bigcup_{h=2}^{K-2} \bar B_{1,K-1}^h(-1,+1)\cup\bigcup_{j=2}^{L-2} \bigcup_{h=1}^{K-1} \bar B_{j,K}^h(-1,+1)\cup\bigcup_{j=2}^{L-2} \bar R_{j,K}(-1,+1)\cup\bigcup_{h=2}^{K-2} \tilde B_{1,K-1}^h(-1,+1)\notag\\&\cup\widetilde{\mathcal K}(-1,+1)\cup\bigcup_{i=1}^{K-2}(\widetilde{\mathcal D}_i(-1,+1)\cup\widetilde{\mathcal E}_i(-1,+1))\cup\tilde R_{1,1}(-1,+1).
\end{align}
\begin{cor}[Tube of typical paths for the Ising Model with zero external magnetic field]\label{cortube}
Consider the Ising model on a $K\times L$ grid $\Lambda$ with periodic boundary conditions and with zero external magnetic field. Then, we have that $\mathfrak T_{\mathbf{+1}}(\mathbf{-1})$ is the tube of typical trajectories for the transition $\mathbf{+1}\to\mathbf{-1}$ and there exists $k>0$ such that for $\beta$ sufficiently large
\begin{align}
\mathbbm P_\beta(\tau_{\partial^{np}\mathfrak T_{\{\mathbf{+1}\}}(\mathbf{-1})}^\mathbf r\le\tau_{\mathbf{+1}}^\mathbf{-1})\le e^{-k\beta}, 
\end{align}
\end{cor}
\section{Minimal restricted-gates}\label{minimalresgatessection}
In order to prove our main results on the set of minimal gates, we first describe the set of all \textit{minimal restricted-gates} for the transition from $\mathbf r\in\mathcal X^s$ to $\mathbf s\in\mathcal X^s$, $\mathbf r\neq\mathbf s$. To do this, we first collect some relevant properties and results concerning the energy landscape between $\mathbf r$ and $\mathbf s$ by studying those optimal paths $\omega\in\Omega_{\mathbf r,\mathbf s}^{opt}$ such that 
\begin{align}\label{conditionopt}
\omega\cap\mathcal X^s\backslash\{\mathbf r,\mathbf s\}=\varnothing.
\end{align}
Note that condition \eqref{conditionopt} implies that the system does not visit $\mathbf r$ and $\mathbf s$ but it does not exclude the presence of some $\sigma\in\omega$ such that exists $v\in V$ with $\sigma(v)\in S\backslash\{r,s\}$.   

\subsection{Energy landscape between two Potts stable configurations}

In \cite[Theorem 2.1]{nardi2019tunneling}, the authors prove that the communication energy \eqref{comheight} between any $\mathbf r,\mathbf{s}\in\mathcal{X}^s$, with $\mathbf r\neq \mathbf s$, is given by
\begin{align}\label{comheightstable}
\Phi(\mathbf r, \mathbf s)=2\min\{K,L\}+2+H(\mathbf r)=2K+2+H(\mathbf r).
\end{align}
where the last equality follows by our assumption $K<L$. Hence, in view of \eqref{comheightstable} we have that any optimal path $\omega\in\Omega_{\mathbf r,\mathbf s}^{opt}$ does not pass through configurations whose energy is strictly larger than $2K+2+H(\mathbf r)$.  
\begin{remark}\label{remarkdaric}
For any $\sigma\in\mathcal X$ and for every $\mathbf r\in \mathcal{X}^s$,
\begin{align}\label{disen}
H(\sigma)-H(\mathbf r)&=H(\sigma)+|E|=|E|-\sum_{(v,w)\in E} \mathbbm{1}_{\{\sigma(v)=\sigma(w)\}}=\sum_{(v,w)\in E} \mathbbm{1}_{\{\sigma(v)\neq\sigma(w)\}}. 
\end{align}
Note that the total number of edges that connect two vertices with different spins, say $r\in S$ and $s\in S\backslash\{r\}$, in a configuration $\sigma\in\mathcal X$ is equal to the perimeter of the same-spin clusters in $C^r_1,C^r_2,\dots$. Thus, thanks to \eqref{comheightstable} and \eqref{disen}, it follows that for any $\sigma$, that belongs to an optimal path $\omega\in\Omega_{\mathbf r,\mathbf s}^{opt}$, the total perimeter of its clusters with the same spin value cannot be larger than $2K+2$.
\end{remark}
The following lemma is an immediate consequence of \eqref{disen}.
\begin{lemma}\label{lemmaenBR}
Consider $\mathbf r,\mathbf s\in\mathcal X^s$, $\mathbf r\neq\mathbf s$. Then, for every $j=1,\dots,L-1$,
\begin{itemize}
\item[\emph{(a)}] any $\sigma\in\bar B_{j,K}^h(r,s)=\tilde B_{L-j-1,K}^{K-h}(r,s)$, $h=1,\dots,K-1$, is such that $H(\sigma)=H(\mathbf s)+2K+2=\Phi(\mathbf r, \mathbf s)$;
\item[\emph{(b)}] any $\sigma\in\bar B_{1,K-1}^h(r,s)$ and any $\sigma\in\tilde B_{1,K-1}^{h}(r,s)$, $h=2,\dots,K-2$, are such that $H(\sigma)=H(\mathbf s)+2K+2=\Phi(\mathbf r, \mathbf s)$;
\item[\emph{(c)}] any $\sigma\in\bar R_{j,K}(r,s)=\tilde R_{L-j,K}(r,s)$ is such that $H(\sigma)=H(\mathbf s)+2K$;
\item[\emph{(d)}] any $\sigma\in\bar R_{j,K-1}(r,s)\cup\tilde R_{j,K-1}(r,s)$ is such that \begin{align}\begin{cases}H(\sigma)=H(\mathbf s)+2K,\ &\text{if}\ j=1;\\ H(\sigma)=H(\mathbf s)+2K+2=\Phi(\mathbf r, \mathbf s),\ &\text{if}\ j=2;\\ H(\sigma)>\Phi(\mathbf r, \mathbf s),\ &\text{if}\ j=3,\dots,L-1.\end{cases}\end{align}
\end{itemize}
\end{lemma}
In the next Lemma we point out which configurations communicate by one step of the dynamics along an optimal path with those states belonging to the sets defined at the beginning of Section \ref{moddepdefnot}. Given $s\in S$, let 
\begin{align}\label{ennekappa}
N_s(\sigma):=|\{v\in V:\ \sigma(v)=s\}|
\end{align}
be the number of vertices with spin $s$ in the configuration $\sigma$. \\
\begin{lemma}\label{strisciastriscia}
Consider $\mathbf r,\mathbf s\in\mathcal X^s$, $\mathbf r\neq\mathbf s$. Given a configuration $\sigma$, let $\bar\sigma$ be a configuration which communicates with $\sigma$ along an optimal path from $\mathbf r$ to $\mathbf s$ that does not intersect $\mathcal X^s\backslash\{\mathbf r,\mathbf s\}$. For any $j=2,\dots,L-2$, the following properties hold:
\begin{itemize}
\item[\emph{(a)}] if $\sigma\in\bar B_{j,K}^h(r,s)$ and $N_s(\sigma)>N_s(\bar\sigma)$, then \[\begin{cases}\bar\sigma\in\bar R_{j,K}(r,s),\ &\text{if}\ h=1;\\\bar\sigma\in\bar B_{j,K}^{h-1}(r,s),\ &\text{if}\ h=2,\dots,K-1;\end{cases}\]
\item[\emph{(b)}] if $\sigma\in\bar B_{j,K}^h(r,s)$ and $N_s(\sigma)<N_s(\bar\sigma)$, then \[\begin{cases}\bar\sigma\in\bar B_{j,K}^{h+1}(r,s),\ &\text{if}\ h=1,\dots,K-2;\\\bar\sigma\in\bar R_{j+1,K}(r,s),\ &\text{if}\ h=K-1;\end{cases}\]
\item[\emph{(c)}] if $\sigma\in\bar R_{j,K}(r,s)$ and $N_s(\sigma)>N_s(\bar\sigma)$, then $\bar\sigma\in\bar B_{j-1,K}^{K-1}(r,s)$;
\item[\emph{(d)}] if $\sigma\in\bar R_{j,K}(r,s)$ and $N_s(\sigma)<N_s(\bar\sigma)$, then $\bar\sigma\in\bar B_{j,K}^{1}(r,s)$.
\end{itemize}
\end{lemma}
\textit{Proof.} Consider $\sigma\in\bar B_{j,K}^h(r,s)$ for some $j=2,\dots,L-2$, $h=1,\dots,K-1$. Let $\bar\sigma\in\mathcal X$ and $v\in V$ be a configuration and a vertex such that flipping the spin in $v$ we move from $\sigma$ to $\bar\sigma$. Thanks to \eqref{energydiff} the following implications hold:
\begin{itemize}
\item[(i)] if $\bar\sigma(v)=t\in S\backslash\{r,s\}$, then $H(\bar\sigma)-H(\sigma)\ge2$;
\item[(ii)] if $\sigma(v)=s$, $v$ has four nearest-neighbor spins $s$ and $\bar\sigma(v)=r$, then $H(\bar\sigma)-H(\sigma)=4$;
\item[(iii)] if $\sigma(v)=r$, $v$ has four nearest-neighbor spins $r$ and $\bar\sigma(v)=s$, then $H(\bar\sigma)-H(\sigma)=4$;
\item[(iv)] if $\sigma(v)=s$, $v$ has three nearest-neighbor spins $s$ and $\bar\sigma(v)=r$ , then $H(\bar\sigma)-H(\sigma)=2$;
\item[(v)]  if $\sigma(v)=r$, $v$ has three nearest-neighbor spins $r$ and $\bar\sigma(v)=s$, then $H(\bar\sigma)-H(\sigma)=2$.
\end{itemize}

Since Lemma \ref{lemmaenBR} holds, it follows that in all the above five cases, we have $H(\bar\sigma)>\Phi(\mathbf r, \mathbf s)$ which is not admissible. Hence, the only configurations $\bar\sigma$ that communicate with $\sigma\in\bar B_{j,K}^h(r,s)$, such that $H(\bar\sigma)\le\Phi(\mathbf r,\mathbf s)$, are those which are obtained by flipping either a spin from $s$ to $r$ or a spin from $r$ to $s$ among the spins with two nearest-neighbor spins $s$ and  two nearest-neighbor spins $r$ in $\sigma$. In particular, following an optimal path from $\sigma\in\bar B_{j,K}^h(r,s)$ to $\mathbf r$ we have \[\begin{cases}\bar\sigma\in\bar R_{j,K}(r,s),\ &\text{if}\ h=1;\\\bar\sigma\in\bar B_{j,K}^{h-1}(r,s),\ &\text{if}\ h=2,\dots,K-1.\end{cases}\] Otherwise, following an optimal path from $\sigma\in\bar B_{j,K}^h(r,s)$ to $\mathbf s$, we have \[\begin{cases}\bar\sigma\in\bar B_{j,K}^{h+1}(r,s),\ &\text{if}\ h=1,\dots,K-2;\\ \bar\sigma\in\bar R_{j+1,K}(r,s),\ &\text{if}\ h=K-1.\end{cases}\]

When $\sigma\in\bar R_{j,K}(r,s)$ for some $j=2,\dots,L-2$, the proof is similar. Indeed, in view of \eqref{energydiff} we have  
\begin{itemize}
\item[(i)] if $\bar\sigma(v)=t\in S\backslash\{r,s\}$, then $H(\bar\sigma)-H(\sigma)=4$;
\item[(ii)] if $\sigma(v)=s$ and $v$ has four nearest-neighbor spins $s$ and $\bar\sigma(v)=r$, then $H(\bar\sigma)-H(\sigma)=4$;
\item[(iii)] if $\sigma(v)=r$ and $v$ has four nearest-neighbor spins $r$ and $\bar\sigma(v)=s$, then $H(\bar\sigma)-H(\sigma)=4$.
\end{itemize}
Moreover, by Lemma \ref{lemmaenBR} we know that $H(\sigma)=2K+H(\mathbf r)$. Thus (i), (ii) and (iii) imply $H(\bar\sigma)>\Phi(\mathbf r,\mathbf s)$ which is not admissible. Hence, for any  $j=2,\dots,L-2$, $\sigma\in\bar R_{j,K}(r,s)$ communicates only with configurations obtained either by flipping a spin from $s$ to $r$ among those spins with three nearest-neighbor spins $s$ and one nearest-neighbor spin $r$ in $\sigma$ or by flipping a spin from $r$ to $s$ among those spins with three nearest-neighbor spins $r$ and one nearest-neighbor spin $s$ in $\sigma$. Following an optimal path from $\sigma$ to $\mathbf r$, these configurations belong to $\bar B_{j-1,K}^{K-1}(r,s)$, while along an optimal path from $\sigma$ to $\mathbf s$ to $\bar B_{j,K}^1(r,s)$. 
$\qed$\\

We remark that thanks to \eqref{equivalenzarbk}, Lemma \ref{strisciastriscia} may be also used for describing the transition between those configurations which belong to either some $\tilde B_{j,K}^h(r,s)$ or some $\tilde R_{j,K}(r,s)$, for $j=2,\dots,L-2$ and $h=1,\dots,K-1$.

\subsection{Geometric properties of the Potts model with zero external magnetic field}
A \textit{two dimensional polyomino} on $\mathbb{Z}^2$ is a finite union of unit squares. The area of a polyomino is the number of its unit squares, while its perimeter is the cardinality of its boundary, namely, the number of unit edges of the dual lattice which intersect only one of the unit squares of the polyomino itself. Thus, the perimeter is the number of interfaces on $\mathbb Z^2$ between the sites inside the polyomino and those outside. We define $M_n$ as the set of all the polyominoes with minimal perimeter among those with area $n$. We call \textit{minimal polyominoes} the elements of $M_n$. \\
Let
\begin{align}\label{nufol}
\mathcal V_n^s:=\{\sigma\in\mathcal X:\ N_s(\sigma)=n\},
\end{align}
be the set of configurations with $n$ spins equal to $s$, with $s\in\{1,\dots,q\}$ and $n=0,\dots,KL$. Note that, given two different stable configurations $\mathbf r, \mathbf s\in\mathcal X^s$, every optimal path from $\mathbf r$ to $\mathbf s$ has to intersect at least one time $\mathcal V_n^s$ for any $n=0,\dots,KL$. We are now able to prove some useful lemmas.

\begin{lemma}\label{lemmaavvolg}
Let $\mathbf r,\mathbf s\in\mathcal{X}^s$ be two different stable configurations and let $\omega$ be an optimal restricted-path
for the transition from $\mathbf r$ to $\mathbf s$.
There exists $K^*\in \{0,1,\ldots,(K-1)^2\}$ such that in any $\sigma\in\omega$ with $K^*+1\le N_s(\sigma)\le (K-1)^2$ at least a cluster of spins $s$ belongs to either $\bar R_{j,K}(r,s)$ or $\bar B_{j,K}^h(r,s)$, for some $h=1,\dots,K-1$ and $j=2,\dots,L$. In other words, at least a cluster of spin $s$ wraps around $\Lambda$. 
\end{lemma}
\textit{Proof.} The strategy for the proof is to construct a path $\tilde\omega:\mathbf r\to\tilde\sigma$ for some $\tilde\sigma\in\bar R_{K-1,K-1}(r,s)$ as a sequence of configurations in which the unique cluster of spins $s$ is a polyomino (cluster that does not wrap around the torus) with minimal perimeter among those with the same area. After this construction we compare $\tilde\omega$ with any optimal path $\omega:\mathbf r\to\mathbf s$ up to $\tilde\eta\in\mathcal V^s_{(K-1)^2}$ with the aim to understand what is the last polyomino with perimeter smaller or equal than $2K+2$, see Remark \ref{remarkdaric}. 
In \cite[Theorem 2.2]{alonso1996three} the authors show that the set of minimal polyominoes of area $n$, $M_n$, includes squares or quasi-squares with possibly a bar on one of the longest sides. Thus, we define $\tilde\omega:=(\tilde\omega_0, \tilde\omega_1, \tilde\omega_2, \dots,\tilde\sigma)$ as the sequence of configurations having these shapes in which the number of spins $s$ increases following the clockwise direction, see Figure \ref{defomegatilde} for an example of this construction. 
\begin{figure}[!ht]
\centering
\begin{tikzpicture}[scale=0.6, transform shape]
\draw [fill=gray,black!15!white] (1.6,1.6) rectangle (2,2);
\fill[black!18!white] (1.6,2) rectangle (2,2.4);
\fill[black!24!white] (2,2) rectangle (2.4,2.4);
\fill[black!27!white] (2,1.6) rectangle (2.4,2);
\fill[black!30!white] (2,1.2) rectangle (2.4,1.6);
\fill[black!33!white] (1.6,1.2) rectangle (2,1.6);
\fill[black!36!white] (1.2,1.2) rectangle (1.6,1.6);
\fill[black!39!white] (1.2,1.6) rectangle (1.6,2);
\fill[black!42!white] (1.2,2) rectangle (1.6,2.4);
\fill[black!45!white] (1.2,2.4) rectangle (1.6,2.8);
\fill[black!48!white] (1.6,2.4) rectangle (2,2.8);
\fill[black!51!white] (2,2.4) rectangle (2.4,2.8);
\fill[black!54!white] (2.4,2.4) rectangle (2.8,2.8);
\fill[black!57!white] (2.4,2) rectangle (2.8,2.4);
\fill[black!60!white] (2.4,1.6) rectangle (2.8,2);
\fill[black!63!white] (2.4,1.2) rectangle (2.8,1.6);
\fill[black!66!white] (2.4,0.8) rectangle (2.8,1.2);
\fill[black!69!white] (2,0.8) rectangle (2.4,1.2);
\fill[black!72!white] (1.6,0.8) rectangle (2,1.2);

\draw[<-,thick] (1,1.8) -- (1,1); \draw[thick] (1,1) -- (1.4,1);
\draw[step=0.4cm,color=black] (0,0) grid (6,4.4);
\end{tikzpicture}%

\caption{\label{defomegatilde} First steps of path $\tilde\omega$ on a $11\times 15$ grid $\Lambda$. The white squares have spin $r$, the other colors denote spin $s$. The arrow indicates the order in which the spins are flipped from $r$ to $s$. The colors of the squares indicate when they have been flipped, with darker squares having been flipped later.}
\end{figure}
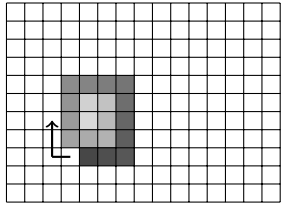\FloatBarrier
Since $N_s(\tilde\omega_{j+1})=N_s(\tilde\omega_j)+1$ for all $j=0,1,\dots$, there exists some $K^*>0$ such that the perimeter of the $s$-cluster of $\tilde\omega_{K^*}$ is equal to $2K+2$ and the perimeter of the $s$-cluster of $\tilde\omega_{K^*+1}$ is strictly larger than $2K+2$. It follows that $H(\tilde\omega_{K^*})=2K + 2 + H(\mathbf r)$ and $H(\tilde\omega_{K^*+1})>2K + 2 + H(\mathbf r)$, thus $\tilde\omega_{K^*+1}$ does not belong to any $\omega\in\Omega_{\mathbf r,\mathbf s}^{opt}$. Explicit computations show that $K^*$ is given by
\begin{align}
K^*:=
\begin{cases}
\frac{K^2+2K+1}{4},\ &\text{if}\ K\ \text{is odd},\\
\frac{K^2+2K}{4},\ &\text{if}\ K\ \text{is even}.
\end{cases}
\end{align}
If $K$ is odd (resp.~even), the $s$-cluster of $\tilde\omega_{K^*}$  is a $\frac{K+1}{2}\times\frac{K+1}{2}$ square (resp.~$\frac{K}{2}\times(\frac{K}{2}+1)$ rectangle).
It follows that the configuration(s) in the intersection $\omega\cap\mathcal V_n^s$ for any $\omega\in\Omega_{\mathbf r,\mathbf s}^{opt}$ and $K^*+1\le n\le (K-1)^2$, do not contain any $\tilde\omega_n$ for any $n\ge K^*+1$.

Hence, any optimal path intersects $\mathcal V_n^s$ with $n>K^*$ only in configurations with $s$-clusters wrapping around $\Lambda$. Moreover, the intersection belongs to either $\bar R_{j,K}(r,s)$ or $\bar B_{j,K}^h(r,s)$, for some $h=1,\dots,K-1$ and $j=2,\dots,L$. $\qed$.

\begin{lemma}\label{intPQH}
Consider $\mathbf r, \mathbf s\in\mathcal X^s$, $\mathbf r\neq\mathbf s$. Then, for any $\omega\in\Omega_{\mathbf r,\mathbf s}^{opt}$ such that $\omega$ is an optimal restricted-path between $\mathbf r$ and $\mathbf s$, 
we have
\begin{itemize}
\item[\emph{(a)}] $\omega\cap\overline{\mathscr P}(\mathbf r,\mathbf s)\neq\varnothing$, $\omega\cap\widetilde{\mathscr P}(\mathbf r,\mathbf s)\neq\varnothing$;
\item[\emph{(b)}] $\omega\cap\overline{\mathcal Q}(\mathbf r,\mathbf s)\neq\varnothing$, $\omega\cap\widetilde{\mathcal Q}(\mathbf r,\mathbf s)\neq\varnothing$;
\item[\emph{(c)}] $\omega\cap\overline{\mathscr H}_i(\mathbf r,\mathbf s)\neq\varnothing$, $\omega\cap\widetilde{\mathscr H}_i(\mathbf r,\mathbf s)\neq\varnothing$ for any $i=1,\dots,K-3$.
\end{itemize} 
\end{lemma}

\textit{Proof.} We begin by proving (a). We prove the statement only for $\overline{\mathscr P}(\mathbf r,\mathbf s)$ since the proof for $\widetilde{\mathscr P}(\mathbf r,\mathbf s)$ follows by switching the roles of $r$ and $s$ and using the symmetry of the model.\\
Let $\omega\in\Omega_{\mathbf r,\mathbf s}^{opt}$ be any optimal restricted-path between $\mathbf r$ and $\mathbf s$.
Thanks to Lemma \ref{lemmaavvolg}, there exists $K^*\in\mathbb{N}$ such that, when $n>K^*$, every $\omega$ intersects $\mathcal{V}_n^s$ in configurations which belong to either $\bar B_{j,K}^h(r,s)$ or $\bar R_{j,K}(r,s)$ for some $j=2,\dots,L-2$ and $h=1,\dots,K-1$. Moreover, in view of Lemma \ref{strisciastriscia}, it follows that $\omega$ reaches these configurations only moving among configurations belonging to either $\bar B_{j,K}^h(r,s)$ or $\bar R_{j,K}(r,s)$ with $j=2,\dots,L-2$, $h=1,\dots,K-1$. 
Note that, given an optimal path $\omega$ from either $\bar B_{j,K}^h(r,s)$ or $\bar R_{j,K}(r,s)$ to $\mathbf r$ and given $\sigma\in\omega\cap\bar R_{2,K}(r,s)$, the only $\bar\sigma\in\mathcal V_{2K-1}^s$ which communicates with $\sigma$ belongs to $\overline{\mathscr P}(\mathbf r,\mathbf s)=\bar B_{1,K}^{K-1}(r,s)$, see Lemma \ref{strisciastriscia}(c). Indeed, we move from $\sigma$ to $\bar\sigma$ by flipping a spin from $s$ to $r$ and the only possibility to not overcome $\Phi(\mathbf r, \mathbf s)$ is to flip a spin $s$ with three nearest-neighbor spins $s$.\\
From now on, using the reversibility of the dynamics, we prefer to study 
$\omega\in\Omega_{\mathbf r,\mathbf s}^{opt}$ an optimal restricted-path between $\mathbf r$ and $\mathbf s$
by analyzing instead its time reversal $\omega^T=(\omega_n,\dots,\omega_0)$. Indeed, a path $\omega=(\omega_0,\dots,\omega_n)$ from $\mathbf r$ to $\mathbf s$ is optimal if and only if the path $\omega^T=(\omega_n,\dots,\omega_0)$ is optimal.\\
Let us move to the proof of (b). To aid the understanding, we suggest to use Figure \ref{figdaPaQ} as a reference for this part of the proof.  We prove the statement only for $\overline{\mathcal Q}(\mathbf r,\mathbf s)$ since the proof for $\widetilde{\mathcal Q}(\mathbf r,\mathbf s)$ again follows from symmetry considerations.\\

Consider $\omega\in\Omega_{\mathbf r,\mathbf s}^{opt}$ such that it is an optimal restricted-path between $\mathbf r$ and $\mathbf s$  
and take $\sigma\in\overline{\mathscr P}(\mathbf r,\mathbf s)\cap\omega$ This exists in view of (a). Note that from Lemma \ref{lemmaenBR} we have $H(\sigma)=\Phi(\mathbf r, \mathbf s)$. Since $\sigma\in\mathcal V_{2K-1}^s$, we have to move from $\sigma$ to $\bar\sigma$ by removing a spin $s$ and the only possibility to not overcome $\Phi(\mathbf r, \mathbf s)$ is to change from $s$ to $r$ a spin $s$ with two nearest-neighbor spins $r$. Indeed, in a such a way the perimeter of the cluster with spins $s$ does not increase and $H(\bar\sigma)$ does not exceed $\Phi(\mathbf r, \mathbf s)$, see Remark \ref{remarkdaric}.
Hence, given $\sigma\in\overline{\mathscr P}(\mathbf r,\mathbf s)$, the only configurations $\bar\sigma\in\mathcal V_{2K-2}^s$ which communicate with $\sigma$, along an optimal path from $\sigma$ to $\mathbf r$, belong to either $\bar R_{2,K-1}(r,s)$ or $\bar B_{1,K}^{K-2}(r,s)$, i.e., to $\overline{\mathcal Q}(\mathbf r,\mathbf s)$.\\

Finally, we carry out the proof of (c). We prove the statement only for $\overline{\mathscr H}_i(\mathbf r,\mathbf s)$, $i=1,\dots,K-3$, since the proof for $\widetilde{\mathscr H}(\mathbf r,\mathbf s)$, $i=1,\dots,K-3$, again follows from symmetry considerations. For simplicity, we split the proof in several steps. \\

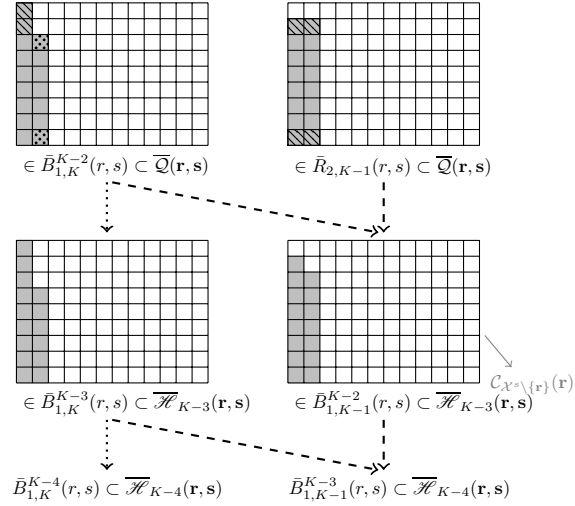
\begin{figure}[h!]
\centering
\begin{tikzpicture}[scale=0.7, transform shape]
\draw [fill=gray,lightgray] (0,0) rectangle (0.3,2.7);
\draw [fill=gray,lightgray] (0.3,0) rectangle (0.6,1.8);
\draw [fill=gray,lightgray] (5.1,0) rectangle (5.4,2.4);
\draw [fill=gray,lightgray] (5.4,0) rectangle (5.7,2.1);
\draw [fill=gray,lightgray] (0,4.5) rectangle (0.3,7.2);
\draw [fill=gray,lightgray] (0.3,4.5) rectangle (0.6,6.6);
\draw [fill=gray,lightgray] (5.1,4.5) rectangle (5.4,6.9);
\draw [fill=gray,lightgray] (5.4,4.5) rectangle (5.7,6.9);
\draw[pattern=north west lines, pattern color=black] (0,6.6) rectangle (0.3,7.2);
\draw[pattern=crosshatch dots, pattern color=black] (0.3,6.3) rectangle (0.6,6.6);
\draw[pattern=crosshatch dots, pattern color=black] (0.3,4.5) rectangle (0.6,4.8);
\draw[pattern=north west lines, pattern color=black] (5.1,6.6) rectangle (5.7,6.9);
\draw[pattern=north west lines, pattern color=black] (5.1,4.5) rectangle (5.7,4.8);
\draw[step=0.3cm,color=black] (0,4.5) grid (3.6,7.2) node[right] at (0,4.1){$\in\bar B_{1,K}^{K-2}(r,s)\subset\overline{\mathcal Q}(\mathbf r,\mathbf s)$};
\draw[step=0.3cm,color=black] (5.1,4.5) grid (8.7,7.2) node[right] at (5.1,4.1){$\in\bar R_{2,K-1}(r,s)\subset\overline{\mathcal Q}(\mathbf r,\mathbf s)$};
\draw[step=0.3cm,color=black] (0,0) grid (3.6,2.7) node[right] at (0,-0.4){$\in\bar B_{1,K}^{K-3}(r,s)\subset\overline{\mathscr H}_{K-3}(\mathbf r,\mathbf s)$};
\draw[step=0.3cm,color=black] (5.1,0) grid (8.7,2.7) node[right] at (5.1,-0.4){$\in\bar B_{1,K-1}^{K-2}(r,s)\subset\overline{\mathscr H}_{K-3}(\mathbf r,\mathbf s)$};
\draw [->,dotted,thick] (1.7, 3.8) -- (1.7,2.85);
\draw [->,dashed,thick] (1.8, 3.8) -- (6.8,2.85);
\draw [->,dashed,thick] (6.9, 3.8) -- (6.9, 2.85);
\draw [->,gray] (8.8,0.9) -- (9.3,0.3) node[below] at (9.7,0.3){$\mathcal C_{\mathcal X^s\backslash\{\mathbf r\}}(\mathbf r)$};
\draw [->,dotted,thick] (1.7, -0.7) -- (1.7,-1.65) node[below] at (1.9,-1.65){$\bar B_{1,K}^{K-4}(r,s)\subset\overline{\mathscr H}_{K-4}(\mathbf r,\mathbf s)$};
\draw [->,dashed,thick] (1.8, -0.7) -- (6.8,-1.65) node[below] at (7.2,-1.65){$\bar B_{1,K-1}^{K-3}(r,s)\subset\overline{\mathscr H}_{K-4}(\mathbf r,\mathbf s)$};
\draw [->,dashed,thick] (6.9, -0.7) -- (6.9, -1.65);
\end{tikzpicture}
\caption{\label{passdaQaQH} Example on a $9\times 12$ grid $\Lambda$ of Step 1. White vertices have spin $r$, gray vertices have spin $s$. Starting from $\bar B_{1,K}^{K-2}(r,s)$, the path can enter into $\overline{\mathscr H}_{K-3}(\mathbf r,\mathbf s)$ by flipping to $r$ either a spin $s$ among those with dots or a spin $s$ among those with lines. Note that from $\bar B_{1,K-1}^{K-3}(r,s)$ the path can enter into $C_{\mathcal X^s\backslash\{\mathbf r\}}(\mathbf r)$ in one step by flipping from $s$ to $r$ the spin $s$ with three nearest-neighbor $r$.}
\end{figure}

\textbf{Step 1.}\ \ \ We claim that, given $\bar\sigma\in\bar R_{2,K-1}(r,s)\cup\bar B_{1,K}^{K-2}(r,s)$,  the only configurations $\hat\sigma\in\mathcal V_{2K-3}^s$ which communicate with $\bar\sigma$, along an optimal path from $\bar\sigma$ to $\mathbf r$, belong to either $\bar B_{1,K-1}^{K-2}(r,s)$ or $\bar B_{1,K}^{K-3}(r,s)$, see Figure \ref{passdaQaQH}.\\
We remark that $\bar\sigma\in\mathcal V_{2K-2}^s$ and, thanks to Lemma \ref{lemmaenBR}, that $H(\bar\sigma)=\Phi(\mathbf r, \mathbf s)$. Hence, we have to move from $\bar\sigma$ to $\hat\sigma$ by removing a spin $s$ without increasing the energy and the only possibility is flipping from $s$ to $r$ a spin $s$ among those with two nearest-neighbor spins $s$. This can happen in many ways. Assume first that $\bar\sigma\in\bar R_{2,K-1}(r,s)$, then $\hat\sigma\in\bar B_{1,K-1}^{K-2}(r,s)\subset\overline{\mathscr H}_i(\mathbf r,\mathbf s)$ for any $i=1,\dots,K-3$ and both the proof of (c) and the proof of the claim are concluded. Otherwise, when $\bar\sigma\in\bar B_{1,K}^{K-2}(r,s)$, we can flip from $s$ to $r$ either
\begin{itemize}
\item[(i)] a spin $s$ with two nearest-neighbor spins $s$ which lies on the column full of $s$ or
\item[(ii)] a spin $s$ among those with two nearest-neighbor spins $s$ on the incomplete column of spins $s$.
\end{itemize}
In particular, in case (i), $\hat\sigma\in\bar B_{1,K-1}^{K-2}(r,s)\subset\overline{\mathscr H}_i(\mathbf r,\mathbf s)$ for any $i=1,\dots,K-3$ and both the proof of (c) and of the claim are completed. Otherwise, in case (ii), $\hat\sigma\in\bar B_{1,K}^{K-3}(r,s)\subset\overline{\mathscr H}_{K-3}(\mathbf r,\mathbf s)$. Thus claim is verified. However, in this case we conclude that only $\overline{\mathscr H}_{K-3}(\mathbf r,\mathbf s)$ is a gate and it is necessary to consider other steps to prove that also $\overline{\mathscr H}_i(\mathbf r,\mathbf s)$ for any $i=1,\dots,K-4$ is a gate. 

\FloatBarrier

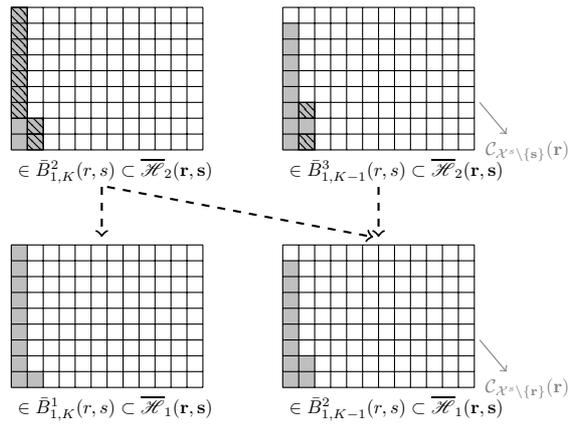
\begin{figure}[h!]
\centering
\begin{tikzpicture}[scale=0.7, transform shape]
\draw [fill=gray,lightgray] (0,0) rectangle (0.3,2.7);
\draw [fill=gray,lightgray] (0.3,0) rectangle (0.6,0.3);
\draw [fill=gray,lightgray] (5.1,0) rectangle (5.4,2.4);
\draw [fill=gray,lightgray] (5.4,0) rectangle (5.7,0.6);
\draw [fill=gray,lightgray] (0,4.5) rectangle (0.3,7.2);
\draw [fill=gray,lightgray] (0.3,4.5) rectangle (0.6,5.1);
\draw [fill=gray,lightgray] (5.1,4.5) rectangle (5.4,6.9);
\draw [fill=gray,lightgray] (5.4,4.5) rectangle (5.7,5.4);
\draw[pattern=north west lines, pattern color=black] (0,5.1) rectangle (0.3,7.2);
\draw[pattern=north west lines, pattern color=black] (0.3,4.8) rectangle (0.6,5.1);
\draw[pattern=north west lines, pattern color=black] (0.3,4.5) rectangle (0.6,4.8);
\draw[pattern=north west lines, pattern color=black] (5.4,4.5) rectangle (5.7,4.8);
\draw[pattern=north west lines, pattern color=black] (5.4,5.1) rectangle (5.7,5.4);
\draw[step=0.3cm,color=black] (0,4.5) grid (3.6,7.2) node[right] at (0,4.1){$\in\bar B_{1,K}^{2}(r,s)\subset\overline{\mathscr H}_2(\mathbf r,\mathbf s)$};
\draw[step=0.3cm,color=black] (5.1,4.5) grid (8.7,7.2) node[right] at (5.1,4.1){$\in\bar B_{1,K-1}^{3}(r,s)\subset\overline{\mathscr H}_2(\mathbf r,\mathbf s)$};
\draw[step=0.3cm,color=black] (0,0) grid (3.6,2.7) node[right] at (0,-0.4){$\in\bar B_{1,K}^{1}(r,s)\subset\overline{\mathscr H}_1(\mathbf r,\mathbf s)$};
\draw[step=0.3cm,color=black] (5.1,0) grid (8.7,2.7) node[right] at (5.1,-0.4){$\in\bar B_{1,K-1}^{2}(r,s)\subset\overline{\mathscr H}_1(\mathbf r,\mathbf s)$};
\draw [->,dashed,thick] (1.7, 3.8) -- (1.7,2.85);
\draw [->,dashed,thick] (1.8, 3.8) -- (6.8,2.85);
\draw [->,dashed,thick] (6.9, 3.8) -- (6.9, 2.85);
\draw [->,gray] (8.8,0.9) -- (9.3,0.3) node[below] at (9.7,0.3){$\mathcal C_{\mathcal X^s\backslash\{\mathbf r\}}(\mathbf r)$};
\draw [->,gray] (8.8,5.4) -- (9.3,4.8) node[below] at (9.7,4.8){$\mathcal C_{\mathcal X^s\backslash\{\mathbf s\}}(\mathbf r)$};
\end{tikzpicture}
\caption{\label{stepfinale} Example on a $9\times 12$ grid $\Lambda$ of the final step of the proof of Lemma \ref{intPQH}(c). White vertices have spin $r$, gray vertices have spin $s$. The lines denote the spins $s$ that can be flipped to $r$ without increasing the energy.} 
\end{figure}\FloatBarrier
\textbf{Step 2.}\ \ \ We claim that, given $\hat\sigma\in\bar  B_{1,K}^{K-3}(r,s)$,  the only configurations of $\mathcal V_{2K-4}^s$ which communicate with $\hat\sigma$, along an optimal path between $\hat\sigma$ and $\mathbf r$, belong to either $\bar B_{1,K-1}^{K-3}(r,s)$ or $\bar B_{1,K}^{K-4}(r,s)$, see Figure \eqref{step2}. \\

Since $H(\hat\sigma)=\Phi(\mathbf r, \mathbf s)$ (see Lemma \ref{lemmaenBR}), then in order to not increase the energy and to reduce the number of spins $s$, the moves (i) and (ii) of Step 1 are the only possibilities. It follows that $\omega$ can pass through $\bar B_{1,K}^{K-3}(r,s)$ coming from a configuration that lies in $\bar B_{1,K-1}^{K-3}(r,s)$ or in $\bar B_{1,K}^{K-4}(r,s)$. In any case, the claim is verified. However, we conclude the proof of (c) only in the first case. Indeed, $\bar B_{1,K-1}^{K-3}(r,s)\subset\overline{\mathscr H}_i(\mathbf r,\mathbf s)$ for any $i=1,\dots,K-4$.
Otherwise, if $\omega$ visits $\bar B_{1,K}^{K-4}(r,s)$, then we conclude that $\overline{\mathscr H}_{K-4}(\mathbf r,\mathbf s)$ is a gate and we have to move on to complete the proof of (c). Iterating the above construction, if at a certain step $\omega$ intersects $\mathcal V_{m}^s$, for $m=K+2,\dots,2K-2$, in a configuration of $\bar B_{1,K-1}^{K-3}(r,s)\subset\overline{\mathscr H}(\mathbf r,\mathbf s)$, then (c) is satisfied and the proof is completed at that step. Otherwise, if $\omega$ intersects every $\mathcal V_{m}^s$, for $m=K+2,\dots,2K-2$, in configurations belonging to $\bar B^i_{1,K}(r,s)$ for $i=1,\dots,K-5$, then the above construction leads to a configuration that lies in $\bar B_{1,K}^2(r,s)\subset\overline{\mathscr H}_2(\mathbf r,\mathbf s)$ and item (c) is satisfied because any $\eta\in\bar B_{1,K}^2(r,s)$ communicates with $\mathcal V_{K+1}^s$ only through configurations belonging to either $\bar B_{1,K}^1(r,s)\subset\overline{\mathscr H}_1(\mathbf r,\mathbf s)$ or to $\bar B_{1,K-1}^2(r,s)\subset\overline{\mathscr H}_1(\mathbf r,\mathbf s)$, see Figure \ref{stepfinale}. Indeed, (i) and (ii) of Step 1 are the only admissible options to move from $\bar B_{1,K}^2(r,s)$ to $\mathcal V_{K+1}^s$ following an optimal path. $\qed$\\

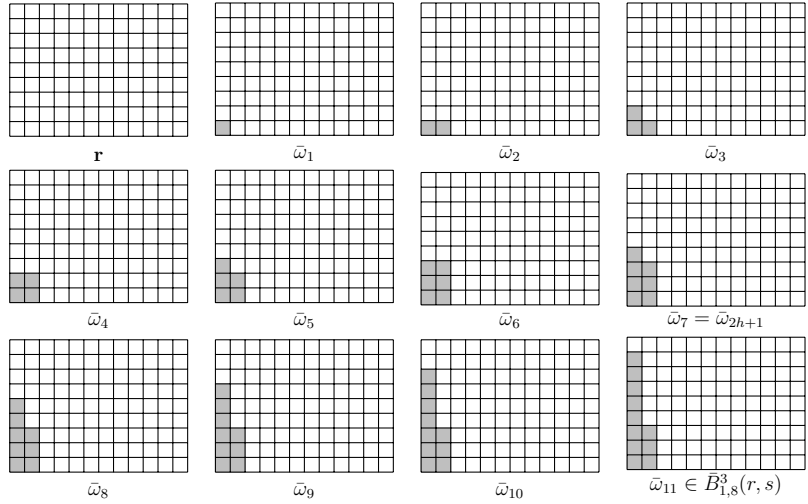
\begin{figure}[h!]
    \centering
\begin{tikzpicture}[scale=0.65, transform shape]
\draw[step=0.3cm,color=black] (0,0) grid (3.6,2.7);
\draw (1.8,-0.2) node[below] {\large$\mathbf r$};
\end{tikzpicture}\ \ \
\begin{tikzpicture}[scale=0.65, transform shape]
\draw [fill=gray,lightgray] (0,0) rectangle (0.3,0.3);
\draw[step=0.3cm,color=black] (0,0) grid (3.6,2.7);
\draw (1.8,-0.1) node[below] {\large$\bar\omega_1$};
\end{tikzpicture}\ \ \
\begin{tikzpicture}[scale=0.65, transform shape]
\draw [fill=gray,lightgray] (0,0) rectangle (0.6,0.3);
\draw[step=0.3cm,color=black] (0,0) grid (3.6,2.7);
\draw (1.8,-0.1) node[below] {\large$\bar\omega_2$};
\end{tikzpicture}\ \ \
\begin{tikzpicture}[scale=0.65, transform shape]
\draw [fill=gray,lightgray] (0,0) rectangle (0.6,0.3);
\draw [fill=gray,lightgray] (0,0.3) rectangle (0.3,0.6);
\draw[step=0.3cm,color=black] (0,0) grid (3.6,2.7);
\draw (1.8,-0.1) node[below] {\large$\bar\omega_3$};
\end{tikzpicture}\\
\begin{tikzpicture}[scale=0.65, transform shape]
\draw [fill=gray,lightgray] (0,0) rectangle (0.6,0.6);
\draw[step=0.3cm,color=black] (0,0) grid (3.6,2.7);
\draw (1.8,-0.1) node[below] {\large$\bar\omega_4$};
\end{tikzpicture}\ \ \
\begin{tikzpicture}[scale=0.65, transform shape]
\draw [fill=gray,lightgray] (0,0) rectangle (0.6,0.6);
\draw [fill=gray,lightgray] (0,0.6) rectangle (0.3,0.9);
\draw[step=0.3cm,color=black] (0,0) grid (3.6,2.7);
\draw (1.8,-0.1) node[below] {\large$\bar\omega_5$};
\end{tikzpicture}\ \ \
\begin{tikzpicture}[scale=0.65, transform shape]
\draw [fill=gray,lightgray] (0,0) rectangle (0.6,0.9);
\draw[step=0.3cm,color=black] (0,0) grid (3.6,2.7);
\draw (1.8,-0.05) node[below] {\large$\bar\omega_6$};
\end{tikzpicture}\ \ \
\begin{tikzpicture}[scale=0.65, transform shape]
\draw [fill=gray,lightgray] (0,0) rectangle (0.6,0.9);
\draw [fill=gray,lightgray] (0,0.9) rectangle (0.3,1.2);
\draw[step=0.3cm,color=black] (0,0) grid (3.6,2.7);
\draw (1.8,0) node[below] {\large$\bar\omega_7=\bar\omega_{2h+1}$};
\end{tikzpicture}\\
\begin{tikzpicture}[scale=0.65, transform shape]
\draw [fill=gray,lightgray] (0,0) rectangle (0.6,0.9);
\draw [fill=gray,lightgray] (0,0.9) rectangle (0.3,1.5);
\draw[step=0.3cm,color=black] (0,0) grid (3.6,2.7);
\draw (1.8,-0.1) node[below] {\large$\bar\omega_8$};
\end{tikzpicture}\ \ \
\begin{tikzpicture}[scale=0.65, transform shape]
\draw [fill=gray,lightgray] (0,0) rectangle (0.6,0.9);
\draw [fill=gray,lightgray] (0,0.9) rectangle (0.3,1.8);
\draw[step=0.3cm,color=black] (0,0) grid (3.6,2.7);
\draw (1.8,-0.1) node[below] {\large$\bar\omega_9$};
\end{tikzpicture}\ \ \
\begin{tikzpicture}[scale=0.65, transform shape]
\draw [fill=gray,lightgray] (0,0) rectangle (0.6,0.9);
\draw [fill=gray,lightgray] (0,0.9) rectangle (0.3,2.1);
\draw[step=0.3cm,color=black] (0,0) grid (3.6,2.7);
\draw (1.8,-0.1) node[below] {\large$\bar\omega_{10}$};
\end{tikzpicture}\ \ \
\begin{tikzpicture}[scale=0.65, transform shape]
\draw [fill=gray,lightgray] (0,0) rectangle (0.6,0.9);
\draw [fill=gray,lightgray] (0,0.9) rectangle (0.3,2.4);
\draw[step=0.3cm,color=black] (0,0) grid (3.6,2.7);
\draw (1.8,0.05) node[below] {\large$\bar\omega_{11}\in\bar B_{1,8}^3(r,s)$};
\end{tikzpicture}
\caption{\label{pathconfH} Example of $\bar\omega: \mathbf r\to\sigma$ of Lemma \ref{pathH} where $\sigma\in\bar B_{1,K-1}^h(r,s)$ with $K=9$ and $h=3$.  White vertices have spin $r$, gray vertices have spin $s$.}
\end{figure}

In the proof of \cite[Proposition 2.5]{nardi2019tunneling}, the authors define a reference path $\omega^*$ between any pair of different stable configurations of a $q$-state Potts model on a $K\times L$ grid $\Lambda$.  Before stating the last lemma of the section, we briefly introduce this path. We say that a path $\omega:\sigma\to\sigma'$ is the \textit{concatenation} of the $L$ paths $\omega^{(i)}=(\omega^{(i)}_0,\dots,\omega^{(i)}_{m_i}),\ \text{for some}\ m_i\in\mathbb{N},\ i=1,\dots,L,$ if $\omega=(\omega^{(1)}_0=\sigma,\dots,\omega^{(1)}_{m_1},\omega^{(2)}_0,\dots,\omega^{(2)}_{m_2},\dots,\omega^{(L)}_0,\dots,\omega^{(L)}_{m_L}=\sigma').$
\begin{definition}\label{remarkrepathNZ}
Given any $\mathbf r, \mathbf s\in\mathcal X^s$, $\mathbf r\neq\mathbf s$, the \textit{reference path} $\omega^*$ is an optimal path from $\mathbf r$ to $\mathbf s$ that is formed by a sequence of configurations in which the cluster of spins $s$ grows gradually column by column. During the first $K$ steps, $\omega^*$ passes through configurations in which the spins on a particular column, say $c_j$ for some $j=0,\dots,L-1$, become $s$, then it crosses those configurations in which the spins on either $c_{j+1}$ or $c_{j-1}$ become $s$ and so on. More precisely, without loss of generality we can start to flip the spins on the first column $c_0$ and define $\omega^*$ as the concatenation of $L$ paths ${\omega^*}^{(1)},\dots,{\omega^*}^{(L)}$ such that ${\omega^*}^{(i)}: \eta_{i-1}\to\eta_{i}$, where $\eta_0:=\mathbf r$, $\eta_L:=\mathbf s$ and for any $i=1,\dots,L$, $\eta_{i}$ is defined as
\begin{align}
\eta_i(v):=\begin{cases}
s,\ &\text{if}\ v\in\bigcup_{j=0}^{i -1} c_j,\\
r,\ &\text{otherwise}.
\end{cases}
\end{align}
In particular, for any $i=1,\dots,L$, we define ${\omega^*}^{(i)}=({\omega^*}^{(i)}_0,\dots,{\omega^*}^{(i)}_K)$ as
\begin{itemize}
\item[-] ${\omega^*}^{(i)}_0=\eta_{i-1}$; 
\item[-] ${\omega^*}^{(i)}_h=({\omega^*}^{(i)}_{h-1})^{(h-1,i),s}$, for $h=1,\dots,K-1$;
\item[-] ${\omega^*}^{(i)}_K=\eta_{i}$. 
\end{itemize}
Note that for any $i=1,\dots,L-1$, $h=1,\dots,K-1$, we have
\begin{itemize}
\item[-] ${\omega^*}^{(1)}_h\in\bar R_{1,h}(r,s)$;
\item[-] $\eta_i\in\bar R_{i,K}(r,s)$;
\item[-] ${\omega^*}^{(i+1)}_h\in\bar B_{i,K}^h(r,s).$
\end{itemize}
Using Lemma \ref{lemmaenBR} and the fact that $\Phi(\mathbf r,\mathbf s)=2K+2+H(\mathbf r)$, we see, indeed, that $\omega^*$ is an optimal path.
\end{definition}
\begin{lemma}\label{pathH}
Consider the $q$-state Potts model on a $K\times L$ grid $\Lambda$ with periodic boundary conditions.  Let $\mathbf r, \mathbf s\in\mathcal X^s$, $\mathbf r\neq\mathbf s$. For any $\sigma\in\bar B^1_{1,K}(r,s)\cup\bigcup_{h=2}^{K-2}\bar B^h_{1,K-1}(r,s)$ there exists a path $\bar\omega=(\bar\omega_0,\dots,\bar\omega_{n})$ from $\mathbf r$ to $\sigma$ such that 
\begin{align}\label{itemlemma}
H(\bar\omega_i)<2K+2+H(\mathbf r),
\end{align}
for any $i=0,\dots,n-1$. Similarly, there exists $\tilde\omega$ from $\mathbf s$ to any $\sigma\in\tilde B^1_{1,K}(r,s)\cup\bigcup_{h=2}^{K-2}\tilde B^h_{1,K-1}(r,s)$ with the same properties of $\bar\omega$.
\end{lemma}
\textit{Proof.} We prove that there exists $\bar\omega: \mathbf r\to\sigma$ which satisfies \eqref{itemlemma} for any $\sigma\in\bar B^1_{1,K}(r,s)\cup\bigcup_{h=2}^{K-2}\bar B^h_{1,K-1}(r,s)$; by reversing the roles of $r$ and $s$, the proof of the existence of $\tilde\omega$ from $\mathbf s$ to any $\sigma\in\tilde B^1_{1,K}(r,s)\cup\bigcup_{h=2}^{K-2}\tilde B^h_{1,K-1}(r,s)$ is analogous. 
The definition of $\bar B^1_{1,K}(r,s)\cup\bigcup_{h=2}^{K-2}\bar B^h_{1,K-1}(r,s)$ gives rise to the two following scenarios, see \eqref{definsHbar}.
If $\sigma\in\bar B_{1,K}^1(r,s)$, then $\bar\omega$ is given by the first steps of the path $\omega^*$ depicted in Definition \ref{remarkrepathNZ}, i.e., $\bar\omega=(\mathbf r, {\omega^*}^{(1)}_1,\dots, {\omega^*}^{(1)}_{K},{\omega^*}^{(2)}_1=\sigma).$ 
Indeed, without loss of generality, we consider $c_0$ as the column in which $\sigma$ has all spins $s$ and the construction of Definition \ref{remarkrepathNZ} holds.\\
Otherwise, we have $\sigma\in\bar B_{1,K-1}^h(r,s)$ for $h=2,\dots,K-2$. Possibly relabeling the columns, we build $\bar\omega$ taking into account the columns $c_0, c_1$ of the grid. In Figure \ref{pathconfH} we depict an example of the path $\bar\omega$.
For every $h=2,\dots,K-2$ and for any odd value $i$ from $1$ to $2h-1$, we set
$\bar\omega_i=\bar\omega_{i-1}^{(\frac{i-1}{2},0),s},\ \ \bar\omega_{i+1}=\bar\omega_i^{(\frac{i-1}{2},1),s}.$
Then, we set
$\bar\omega_j=\bar\omega_{j-1}^{(j-h,0)},$
for any $j=2h+1,\dots,K-1+h.$ 
As we see in Figure \ref{pathconfH}, after $2h$ steps $\bar\omega$ arrives in $\bar\omega_{2h}\in\bar R_{2,h}(r,s)$ and its next configurations belong to $\bar B_{1,j-h}^h(r,s)$ for $j=2h+1,\dots,K-1+h$. Finally, \eqref{itemlemma} is satisfied in view of \eqref{disen}.
$\qed$\\
\FloatBarrier

\subsection{Study of the set of all minimal restricted-gates between two different stable states}
We are now able to prove the following results concerning the set of minimal restricted-gates from $\mathbf r\in\mathcal X^s$ to $\mathbf s\in\mathcal X^s, \mathbf s\neq\mathbf r$.\\

\textit{Proof of Theorem \ref{mingatescond}.} 
In order to prove that a set $\mathcal W_{\text{\tiny{RES}}} \subset\mathcal S(\mathbf r, \mathbf s)$ is a minimal restricted-gate for the transition from $\mathbf r$ to $\mathbf s$ we show that 
\begin{itemize}
\item[(i)] $\mathcal W_{\text{\tiny{RES}}}$ is a restricted-gate, i.e., every $\omega\in\Omega_{\mathbf r,\mathbf s}^{opt}$ such that  it is an optimal restricted-path between $\mathbf r$ and $\mathbf s$
intersects $\mathcal W_{\text{\tiny{RES}}}$,
\item[(ii)] for any $\eta\in\mathcal W_{\text{\tiny{RES}}}$ there exists an optimal path $\omega' \in\Omega_{\mathbf r,\mathbf s}^{opt}$ such that $\omega'\cap(\mathcal W_{\text{\tiny{RES}}}\backslash\{\eta\})=\varnothing$.
\end{itemize} 
Hence, we now show that the sets defined in (a), (b), (c) and (d) of Theorem \ref{mingatescond} satisfy the conditions above. Using Lemma \ref{intPQH}(a), $\overline{\mathscr P}(\mathbf r,\mathbf s)$ and $\widetilde{\mathscr P}(\mathbf r,\mathbf s)$ are gates for the transition from $\mathbf r$ to $\mathbf s$. Next let us show that for any $\eta\in\overline{\mathscr P}(\mathbf r,\mathbf s)$ there exists an optimal path $\omega' \in\Omega_{\mathbf r,\mathbf s}^{opt}$ such that $\omega'\cap({\overline{\mathscr P}(\mathbf r,\mathbf s)\backslash\{\eta\})}=\varnothing$. It is enough to consider $\omega'$ as the path $\omega^*$ of Definition \ref{remarkrepathNZ} and to rewrite it in order to have $\omega'\cap\overline{\mathscr P}(\mathbf r,\mathbf s)=\{\eta\}$, i.e., ${\omega^*}^{(2)}_{K-1}=\eta$.
By the symmetry of the model, we prove similarly that there exists a such a path also for any $\eta\in\widetilde{\mathscr P}(\mathbf r,\mathbf s)$.\\
Using Lemma \ref{intPQH}(b), $\overline{\mathcal Q}(\mathbf r,\mathbf s)$ and $\widetilde{\mathcal Q}(\mathbf r,\mathbf s)$ are gates for the transition from $\mathbf r$ to $\mathbf s$. Next let us show that for any $\eta\in\overline{\mathcal Q}(\mathbf r,\mathbf s)$ there exists an optimal path $\omega' \in\Omega_{\mathbf r,\mathbf s}^{opt}$ such that $\omega'\cap\left(\overline{\mathcal Q}(\mathbf r,\mathbf s)\backslash\{\eta\}\right)=\varnothing$. We distinguish two cases:
\begin{itemize}
\item[(i)] if $\eta\in\bar R_{2,K-1}(r,s)$, given $\bar\eta\in\bar B_{1,K-1}^{K-2}(r,s)$ and  $\hat\eta\in\overline{\mathscr P}(\mathbf r,\mathbf s)$ which communicate with $\eta$, then $\omega'$ is the optimal path given by the concatenation of
\begin{itemize}
\item[-] the path $\bar\omega: \mathbf r\to\bar\eta$ of Lemma \ref{pathH};
\item[-] the path $(\bar\eta,\eta,\hat\eta)$; 
\item[-] the portion of the path $\omega^*$ in Definition \ref{remarkrepathNZ} from  ${\omega^*}^{(2)}_{K-1}=\hat\eta$ to $\mathbf s$, 
\end{itemize} 
so that $\omega'\cap\overline{\mathcal Q}(\mathbf r,\mathbf s)=\{\eta\}$. 
\item[(ii)] if $\eta\in\bar B_{1,K}^{K-2}(r,s)$, then to define $\omega'$ it is enough to consider the path $\omega^*$ of Definition \ref{remarkrepathNZ} and to construct it in order to have $\omega'\cap\overline{\mathcal Q}(\mathbf r,\mathbf s)=\{\eta\}$, i.e., ${\omega^*}^{(2)}_{K-2}=\eta$.
\end{itemize}
Thanks to the symmetry of the model, we define $\omega'$ in an analogous way for any $\eta\in\widetilde{\mathcal Q}(\mathbf r,\mathbf s)$.\\
Using Lemma \ref{intPQH}(c), for any $i=1,\dots,K-3$, $\overline{\mathscr H}_i(\mathbf r,\mathbf s)$ and $\widetilde{\mathscr H}_i(\mathbf r,\mathbf s)$ are gates for the transition from $\mathbf r$ to $\mathbf s$. Next let us show that for any $\eta\in\overline{\mathscr H}_i(\mathbf r,\mathbf s)$ and any $i=1,\dots,K-3$ there exists an optimal path $\omega' \in\Omega_{\mathbf r,\mathbf s}^{opt}$ such that $\omega'\cap\left(\overline{\mathscr H}_i(\mathbf r,\mathbf s)\backslash\{\eta\}\right)=\varnothing$. We have to distinguish two cases:
\begin{itemize}
\item[(i)] if $\eta\in\bar B_{1,K}^i(r,s)$, then $\omega'$ is given by the path $\omega^*$ of Definition \ref{remarkrepathNZ} defined in order to have $\omega'\cap\overline{\mathscr H}_i(\mathbf r,\mathbf s)=\{\eta\}$, i.e., ${\omega^*}^{(2)}_{i}=\eta$;
\item[(ii)] if $\eta\in\bar B_{1,K-1}^h(r,s)$, for some $h=2,\dots,K-2$, then $\omega'$ corresponds to the optimal path given by the concatenation of
\begin{itemize}
\item[-] the path $\bar\omega: \mathbf r\to\eta$ of Lemma \ref{pathH};
\item[-] the path $(\eta,\bar\eta)$ with $\bar\eta\in\bar B_{1,K}^h(r,s)$, such that the bar of length $h$ is in the same position as in $\eta$;
\item[-] the portion of the path $\omega^*$ in Definition \ref{remarkrepathNZ} from ${\omega^*}^{(2)}_h=\bar\eta$ to $\mathbf s$,
\end{itemize} 
so that $\omega'\cap\overline{\mathscr H}_i(\mathbf r,\mathbf s)=\{\eta\}$.
\end{itemize}
Thanks to the symmetry of the model, we define $\omega'$ for any $\eta\in\widetilde{\mathscr H}_i(\mathbf r,\mathbf s)$ and any $i=1,\dots,K-3$ following the same strategy.\\
From Lemma \ref{lemmaavvolg} and Lemma \ref{strisciastriscia}, we conclude that $\mathcal W_j^{(h)}(\mathbf r,\mathbf s)$ are gates for the transition $\mathbf r\to\mathbf s$ for any $j=2,\dots,L-3$ and any $h=1,\dots,K-1$. Indeed, by Lemma \ref{lemmaavvolg}, there exists $K^*\in\mathbb{N}$ such that when $n>K^*$ every $\omega\in\Omega_{\mathbf r,\mathbf s}^{opt}$ intersects $\mathcal{V}_n^s$ in configurations which belong to either $\bar R_{j,K}(r,s)$ or $\mathcal W_j^{(h)}(\mathbf r,\mathbf s)=\bar B_{j,K}^h(r,s)$  for some $j=2,\dots,L-3$, $h=1,\dots,K-1$. Moreover, by Lemma \ref{strisciastriscia}, we know that $\omega$ reaches these configurations only moving among configurations lying either in $\bar R_{j,K}(r,s)$ or in $\mathcal W_j^{(h)}(\mathbf r,\mathbf s)=\bar B_{j,K}^h(r,s)$. Hence, between its last visit to $\overline{\mathscr P}(\mathbf r,\mathbf s)$ and its first visit to $\widetilde{\mathscr P}(\mathbf r,\mathbf s)$, $\omega$ passes at least once through each $\mathcal W_j^{(h)}(\mathbf r,\mathbf s)$, $j=2,\dots,L-3$. Thus, to conclude the proof we have to show that for every $\eta\in W_j^{(h)}(\mathbf r,\mathbf s)$, there exists a path $\omega'\in\Omega_{\mathbf r,\mathbf s}^{opt}$ such that $\omega'\cap(\mathcal W_j^{(h)}(\mathbf r,\mathbf s)\backslash\{\eta\})=\varnothing$. For any $j=2,\dots,L-3$ and any $h=1,\dots,K-1$, we can define this path $\omega'$ as the path $\omega^*$ of Definition \ref{remarkrepathNZ}, which we rewrite in order to have $\omega'\cap\mathcal W_j^{(h)}(\mathbf r,\mathbf s)=\{\eta\}$, i.e., ${\omega^*}^{(j+1)}_{h}=\eta$. $\qed$\\

\begin{remark}\label{remuness}
A saddle $\eta\in\mathcal S(\sigma,\sigma')$ is unessential if for any $\omega\in\Omega_{\sigma,\sigma'}^{opt}$ such that $\omega\cap\eta\neq\varnothing$ the following conditions are both satisfied:
\begin{itemize}
\item[(a)] $\{\text{argmax}_\omega H\}\backslash \{\eta\}\neq\varnothing$,
\item[(b)] there exists $\omega'\in\Omega_{\sigma,\sigma'}^{opt}$ such that $\{\text{argmax}_{\omega'} H\}\subseteq\{\text{argmax}_\omega H\}\backslash \{\eta\}$.
\end{itemize}
\end{remark}
\textit{Proof of Theorem \ref{mingatescondset}.} In view of Theorem \ref{mingatescond}, we have \[\bigcup_{j=2}^{L-3} \mathcal W_j(\mathbf r,\mathbf s) \cup \overline{\mathscr H}(\mathbf r,\mathbf s)\cup \widetilde{\mathscr H}(\mathbf r,\mathbf s)\cup\overline{\mathcal Q}(\mathbf r,\mathbf s)\cup\widetilde{\mathcal Q}(\mathbf r,\mathbf s)\cup\overline{\mathscr P}(\mathbf r,\mathbf s)\cup\widetilde{\mathscr P}(\mathbf r,\mathbf s)\subseteq\mathcal F(\mathbf r,\mathbf s).\] Hence, we have only to prove the opposite inclusion. In order to do this, we use the characterization of minimal gates as essential saddles given in \cite[Theorem 5.1]{manzo2004essential}. Thus, if we prove that any 
\begin{align}\label{sellanonessenziale}
\eta\in\mathcal S(\mathbf r, \mathbf s)\backslash\biggl[\bigcup_{j=2}^{L-3} \mathcal W_j(\mathbf r,\mathbf s) \cup \overline{\mathscr H}(\mathbf r&,\mathbf s)\cup\widetilde{\mathscr H}(\mathbf r,\mathbf s)\cup\overline{\mathcal Q}(\mathbf r,\mathbf s)\cup\widetilde{\mathcal Q}(\mathbf r,\mathbf s)\cup\overline{\mathscr P}(\mathbf r,\mathbf s)\cup\widetilde{\mathscr P}(\mathbf r,\mathbf s)\biggr]\end{align} is an unessential saddle for the restricted transition from $\mathbf r$ to $\mathbf s$, the proof is completed.\\
Consider $\eta$ as in \eqref{sellanonessenziale} and some $\omega\in\Omega_{\mathbf r,\mathbf s}^{opt}$, $\omega=(\omega_0,\dots,\omega_n)$, such that $\omega\cap\mathcal X^s\backslash\{\mathbf r,\mathbf s\}$ and such that $\eta\in\omega$. See Figure \ref{figdimteoG} for an example of $\omega$.  By Lemma \ref{strisciastriscia} and Lemma \ref{lemmaenBR} the condition (i) of Remark \ref{remuness} is satisfied. Indeed, any $\omega\in\Omega_{\mathbf r,\mathbf s}^{opt}$ passes through many configurations with energy value equal to $\Phi(\mathbf r, \mathbf s)$.\ Hence, to conclude that $\eta$ is an unessential saddle we have to prove that condition (b) of Remark \ref{remuness} is verified.
By Lemma \ref{intPQH}(c), there exist $\bar\eta\in\omega\cap\overline{\mathscr H}(\mathbf r,\mathbf s)$ and $\tilde\eta\in\omega\cap\widetilde{\mathscr H}(\mathbf r,\mathbf s)$, where $\bar\eta$ is the last configuration visited by $\omega$ in $\overline{\mathscr H}(\mathbf r,\mathbf s)$ and $\tilde\eta$ is the first configuration visited by $\omega$ in $\widetilde{\mathscr H}(\mathbf r,\mathbf s)$. Moreover, Lemma \ref{intPQH}(a) implies that there exist $\bar\eta^*\in\omega\cap\overline{\mathscr P}(\mathbf r,\mathbf s)$ and $\tilde\eta^*\in\omega\cap\widetilde{\mathscr P}(\mathbf r,\mathbf s)$, where $\bar\eta^*$ is the last configuration visited by $\omega$ in $\overline{\mathscr P}(\mathbf r,\mathbf s)$ and $\tilde\eta^*$ is the first configuration visited by $\omega$ in $\widetilde{\mathscr P}(\mathbf r,\mathbf s)$.

\begin{figure}[h!]
       \begin{minipage}[c]{1\textwidth}
    \centering
    \makebox[0pt]{%
\begin{tikzpicture}[scale=0.8, transform shape]
\fill[color=black] (0.2,1.8)  circle (1.5pt) node[left] {$\mathbf r$};
\fill[color=black] (13.8,0.3)  circle (1.5pt) node[right] {$\mathbf s$};
\draw[dotted] (2.1,-2) -- (2.1,3);\draw[<-] (2,2.8) -- (1.6,3.1) node[left]{$\mathcal V_{K+1}^s$};
\draw[dotted] (11.9,-2) -- (11.9,3);\draw[<-] (12,2.8) -- (12.4,3.1) node[right]{$\mathcal V_{LK-K-1}^s$};
\draw[dotted] (2.5,-2) -- (2.5,3);\draw[<-] (2.5,3.1) -- (2.5,3.5) node[above] {$\mathcal V_{K+2}^s$};
\draw[dotted] (11.5,-2) -- (11.5,3);\draw[<-] (11.5,3.1) -- (11.5,3.5) node[above] {$\mathcal V_{LK-K-2}^s$};
\draw[dotted] (3.1,-2) -- (3.1,3);
\draw[dotted] (3.5,-2) -- (3.5,3);
\draw[dotted] (3.9,-2) -- (3.9,3);\draw[<-] (3.9,3.1) -- (4.3,3.5) node[right] {$\mathcal V_{2K-1}^s$};
\draw[dotted] (10.9,-2) -- (10.9,3);
\draw[dotted] (10.5,-2) -- (10.5,3);
\draw[<-] (10.1,3.1) -- (9.6,3.5) node[left] {$\mathcal V_{LK-2K+1}^s$};
\draw[dotted] (10.1,-2) -- (10.1,3);

\draw[dashed, thick] (3.1,1.6) -- (2.1,-0.6); 
\draw[dashed,thick] (10.9,1.6) -- (11.9,-0.6); 
\draw[ultra thick] (3.9,1.5) -- (3.9,0.1); 
\draw[ultra thick] (10.1,1.5) -- (10.1,0.1); 
\draw[ultra thick,black!40!white] (3.5,1.5) -- (3.5,0.1); 
\draw[ultra thick,black!40!white] (10.5,1.5) -- (10.5,0.1); 

\draw[dotted,ultra thick,black!15!white] (0.2,1.8) -- (2.75,0.8) -- (5.8,0.8);
\draw[dotted,ultra thick,black!15!white]  (8.2,0.8) -- (11.25,0.8) -- (13.8,0.3);
\draw (0.2,1.8) -- (4,1.8) arc (90:-90:0.08) -- (2.5,1.6)  arc (90:270:0.1) -- (3.9,1.4) arc (90:-90:0.1) -- (2.7,1.2) arc (90:270:0.2) -- (5.8,0.8);
\draw (13.8,0.3) -- (10.5,0.3) arc (270:90:0.15) -- (12.1,0.6) arc (-90:90:0.1) -- (8.2,0.8);
\draw[dashed] (6.2,0.8) -- (7.8,0.8);

\fill[color=black] (3.9,0.8) circle (1.5pt) node[below right] {$\bar\eta^*$};
\fill[color=black] (2.75,0.8) circle (1.5pt); \draw (2.65,0.8) node[below right] {$\bar\eta$};
\fill[color=black] (10.1,0.8) circle (1.5pt) node[above left] {$\tilde\eta^*$};
\fill[color=black] (11.25,0.8) circle (1.5pt); \draw (11.15,0.8) node[above right] {$\tilde\eta$};
\fill[color=black] (3.9,1.8) circle (1.5pt) node[above right] {$\eta$};

\draw[black!40!white,ultra thick] (0.5,-2.5) -- (1.26,-2.5) node[right,black] {$\bar B_{1,K}^{K-2}(r,s)\cup\bar R_{2,K-1}(r,s)\subset\overline{\mathcal Q}(\mathbf r,\mathbf s)$}; 
\draw[dashed,thick] (0.5,-3) -- (1.25,-3) node[right] {$\overline{\mathscr H}(\mathbf r,\mathbf s)$}; 
\draw[black, ultra thick] (0.5,-3.5) -- (1.25,-3.5) node[right] {$\overline{\mathscr P}(\mathbf r,\mathbf s)$}; 
\draw[ultra thick,black!40!white] (8.5,-2.5) -- (9.25,-2.5) node[right,black] {$\tilde B_{1,K}^{K-2}(r,s)\cup\tilde R_{2,K-1}(r,s)\subset\widetilde{\mathcal Q}(\mathbf r,\mathbf s)$}; 
\draw[dashed,thick] (8.5,-3) -- (9.25,-3) node[right,black] {$\widetilde{\mathscr H}(\mathbf r,\mathbf s)$}; 
\draw[black, ultra thick] (8.5,-3.5) -- (9.25,-3.5) node[right] {$\widetilde{\mathscr P}(\mathbf r,\mathbf s)$}; 
\end{tikzpicture}%
    }\par
     \end{minipage}
     \caption{\label{figdimteoG} Example of the paths $\omega$ (solid black path) and $\omega'$ (dotted gray path) of the proof of Theorem \ref{setmingatefin}. }
\end{figure}
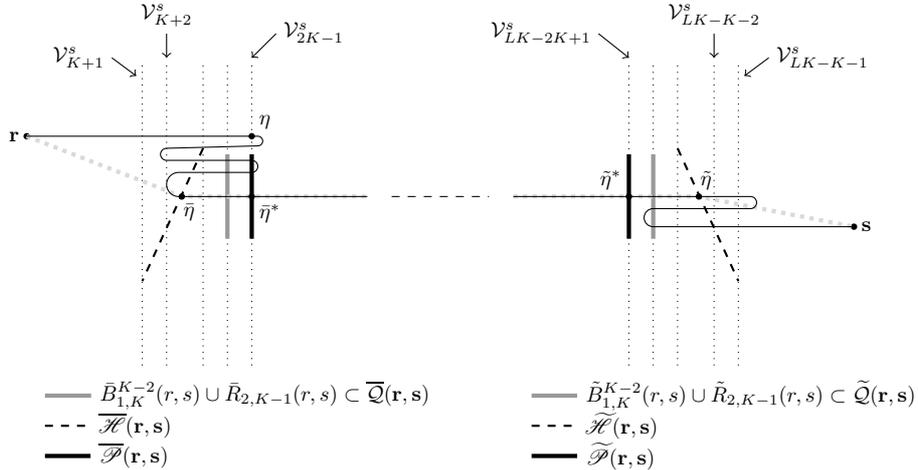\FloatBarrier
 In view of the proof of Lemma \ref{intPQH}, after visiting $\bar\eta$, $\omega$ interstects $\mathcal S(\mathbf r, \mathbf s)$ only in saddles belonging to either $\overline{\mathcal Q}(\mathbf r,\mathbf s)$ or $\overline{\mathscr P}(\mathbf r,\mathbf s)$ or $ \mathcal W_j(\mathbf r,\mathbf s)$, for some $j=2,\dots,L-3$, until it intersects $\widetilde{\mathscr P}(\mathbf r,\mathbf s)$ in $\tilde\eta^*$.  Similarly, after the visit in $\bar\eta^*$ and before the arrival in $\tilde\eta$, $\omega$ passes only through saddles belonging to either $\mathcal W_j(\mathbf r,\mathbf s)$, for some $j=2,\dots,L-3$, or $\widetilde{\mathscr P}(\mathbf r,\mathbf s)$ or $\widetilde{\mathcal Q}(\mathbf r,\mathbf s)$.
It follows that after $\bar\eta$ and before $\tilde\eta$, $\omega$ intersects $\mathcal S(\mathbf r, \mathbf s)$ only in those saddles which belong to $\bigcup_{j=2}^{L-3} \mathcal W_j(\mathbf r,\mathbf s) \cup \overline{\mathscr H}(\mathbf r,\mathbf s)\cup\widetilde{\mathscr H}(\mathbf r,\mathbf s)\cup\overline{\mathcal Q}(\mathbf r,\mathbf s)\cup\widetilde{\mathcal Q}(\mathbf r,\mathbf s)\cup\overline{\mathscr P}(\mathbf r,\mathbf s)\cup\widetilde{\mathscr P}(\mathbf r,\mathbf s)$.\\
Now consider the paths $\bar\omega: \mathbf r\to\bar\eta$ and $\tilde\omega: \mathbf s\to\tilde\eta$, which exist in view of Lemma \ref{pathH}, and take the time reversal of $\tilde\omega$, i.e., $\tilde\omega^T=(\omega_n=\tilde\eta,\tilde\omega_{n-1},\dots,\omega_1,\omega_0=\mathbf s)$.
Thus, if $\omega=(\omega_0=\mathbf r,\dots,\omega_i=\bar\eta,\dots,\omega_j=\tilde\eta,\dots,\omega_n=\mathbf s)$, the path $\omega'\in\Omega_{\mathbf r,\mathbf s}^{opt}$, may be defined as
\begin{itemize}
\item[-] $\omega'\equiv\bar\omega$ from $\mathbf r$ to $\bar\eta$;
\item[-] $\omega'\equiv(\omega_i=\bar\eta,\dots,\omega_j=\tilde\eta)$ from $\bar\eta$ to $\tilde\eta$;
\item[-] $\omega'\equiv\tilde\omega^T$ from $\tilde\eta$ to $\mathbf s$.
\end{itemize} 
Thus, (b) of Remark \ref{remuness} is verified.
$\qed$\\

\textit{Proof of Corollary \ref{corgatecond}.} Since Theorem \ref{mingatescond} holds, the corollary follows by \cite[Theorem 5.4]{manzo2004essential}. $\qed$\\

\section{Minimal gates}\label{proofmaingates}
We are now able to carry out the proof of the main results on the minimal gates for the transitions from a stable state to the other stable configurations and from a stable state to another stable configuration.

\subsection{The minimal gates from a stable state to the other stable states}\label{proofappendixsub}

\ \ \ \ \textit{Proof of Theorem \ref{mingatesNOcond}.} We recall that a subset $\mathcal W\subset\mathcal S(\mathbf r,\mathcal X^s\backslash\{\mathbf r\})$ is a minimal gate for the transition $\mathbf r\to\mathcal X^s\backslash\{\mathbf r\}$ if  
\begin{itemize}
\item[(i)] $\mathcal W$ is a gate for the transition $\mathbf r\to\mathcal X^s\backslash\{\mathbf r\}$, i.e., every $\omega\in\Omega_{\mathbf r,\mathcal X^s\backslash\{\mathbf r\}}^{opt}$ intersects $\mathcal W$, 
\item[(ii)] for any $\eta\in\mathcal W$ there exists an optimal path $\omega' \in\Omega_{\mathbf r,\mathcal X^s\backslash\{\mathbf r\}}^{opt}$ such that $\omega'\cap(\mathcal W\backslash\{\eta\})=\varnothing$.
\end{itemize}  
We begin to prove that the sets depicted in (a) of Theorem \ref{mingatesNOcond} are minimal gates for the transition $\mathbf r\to\mathcal X^s\backslash\{\mathbf r\}$.
Consider any $\omega\in\Omega_{\mathbf r,\mathcal X^s\backslash\{\mathbf r\}}^{opt}$ and let $\mathbf s\in\mathcal X^s\backslash\{\mathbf r\}$ be the first configuration visited by $\omega$ in $\mathcal X^s\backslash\{\mathbf r\}$.   From Theorem \ref{mingatescond}(a) we have $\omega\cap\overline{\mathscr P}(\mathbf r,\mathbf s)\neq\varnothing$ and $\omega\cap\widetilde{\mathscr P}(\mathbf r,\mathbf s)\neq\varnothing$. Thus, \[\omega\cap\biggl(\bigcup_{\mathbf t\in\mathcal X^s\backslash \{\mathbf r\}}\overline{\mathscr P}(\mathbf r,\mathbf t)\biggr)\neq\varnothing\ \ \ \text{and}\ \ \ \omega\cap\biggl(\bigcup_{\mathbf t\in\mathcal X^s\backslash \{\mathbf r\}}\widetilde{\mathscr P}(\mathbf r,\mathbf t)\biggr)\neq\varnothing\]
and (i) is verified. \\
Now consider $\eta\in\bigcup_{\mathbf t\in\mathcal X^s\backslash \{\mathbf r\}}\overline{\mathscr P}(\mathbf r,\mathbf t)$. There exists $\mathbf s\in\mathcal X^s\backslash\{\mathbf r\}$ such that $\eta\in\overline{\mathscr P}(\mathbf r,\mathbf s)$. Let $\omega' \in\Omega_{\mathbf r,\mathcal X^s\backslash\{\mathbf r\}}^{opt}$ be the optimal path from $\mathbf r$ to $\mathbf s\in\mathcal X^s\backslash\{\mathbf r\}$ constructed in the proof of Theorem \ref{mingatescond}(a), such that $\omega'\cap(\bigcup_{\mathbf t\in\mathcal X^s\backslash \{\mathbf r\}}\overline{\mathscr P}(\mathbf r,\mathbf t)\backslash\{\eta\})=\varnothing$ and $\omega'\cap\mathcal X^s\backslash\{\mathbf r,\mathbf s\}=\varnothing$. Hence, (ii) is verified for $\bigcup_{\mathbf t\in\mathcal X^s\backslash \{\mathbf r\}}\overline{\mathscr P}(\mathbf r,\mathbf t)$. By the symmetry of the model, we can argue similarly to prove that (ii) holds also for $\bigcup_{\mathbf t\in\mathcal X^s\backslash \{\mathbf r\}}\widetilde{\mathscr P}(\mathbf r,\mathbf t)$.\\
Following the same strategy, we prove that the sets in (b), (c) and (d) are minimal gates for the transition from $\mathbf r$ to $\mathcal X^s\backslash\{\mathbf r\}$ by replacing Theorem \ref{mingatescond}(a) with Theorem \ref{mingatescond}(b), (c), (d), respectively. See Appendix \ref{appendixproof1} for the details. $\qed$\\

\textit{Proof of Theorem \ref{setmingatesNOcond}.} In view of Theorem \ref{mingatesNOcond} we have
\begin{align}
\bigcup_{\mathbf t\in\mathcal X^s\backslash \{\mathbf r\}}\biggl[\bigcup_{j=2}^{L-3} \mathcal W_j(\mathbf r,\mathbf t)\cup\overline{\mathscr H}&(\mathbf r,\mathbf t)\cup \widetilde{\mathscr H}(\mathbf r,\mathbf t)\cup\overline{\mathcal Q}(\mathbf r,\mathbf t)\notag\\&\cup\widetilde{\mathcal Q}(\mathbf r,\mathbf t)\cup\overline{\mathscr P}(\mathbf r,\mathbf t)\cup \widetilde{\mathscr P}(\mathbf r,\mathbf t)\biggr]\subseteq\mathcal G(\mathbf r,\mathcal X^s\backslash\{\mathbf r\}).\notag
\end{align}
Hence, we only have to prove the opposite inclusion. In order to do this, we again use the characterization of minimal gates as essential saddles given in \cite[Theorem 5.1]{manzo2004essential}. Thus, our strategy is to prove that any
\begin{align}\label{saddleunessentialbis}
\eta\in\mathcal S(\mathbf r,\mathcal X^s\backslash\{\mathbf r\})\backslash\bigcup_{\mathbf t\in\mathcal X^s\backslash \{\mathbf r\}}\biggl[\bigcup_{j=2}^{L-3}& \mathcal W_j(\mathbf r,\mathbf t)\cup\overline{\mathscr H}(\mathbf r,\mathbf t)\cup \widetilde{\mathscr H}(\mathbf r,\mathbf t)\cup\overline{\mathcal Q}(\mathbf r,\mathbf t)\notag\\&\cup\widetilde{\mathcal Q}(\mathbf r,\mathbf t)\cup\overline{\mathscr P}(\mathbf r,\mathbf t)\cup \widetilde{\mathscr P}(\mathbf r,\mathbf t)\biggr]
\end{align}
is an unessential saddle. In particular, for any saddle $\eta$ as in \eqref{saddleunessentialbis} and for any $\omega\in\Omega_{\mathbf r,\mathcal X^s\backslash\{\mathbf r\}}^{opt}$ such that $\eta\in\omega$, we have to prove that the conditions of Remark \ref{remuness} are satisfied. Let $\mathbf s\in\mathcal X^s\backslash\{\mathbf r\}$ be the first stable state visited by $\omega$. By Lemma \ref{strisciastriscia} and Lemma \ref{lemmaenBR} the condition (a) of Remark \ref{remuness} is satisfied. Moreover, let $\omega'\in\Omega_{\mathbf r,\mathcal X^s\backslash\{\mathbf r\}}^{opt}$ be the optimal path from $\mathbf r$ to $\mathbf s\in\mathcal X^s\backslash\{\mathbf r\}$ constructed in the proof of Theorem \ref{mingatescondset} such that \[\{\text{argmax}_{\omega'} H\}\subseteq\{\text{argmax}_\omega H\}\backslash \{\eta\}\] and $\omega'\cap\mathcal X^s\backslash\{\mathbf r,\mathbf s\}=\varnothing$. Thus, condition (b) of Remark \ref{remuness} is satisfied and $\eta$ is an unessential saddle. $\qed$\\

\textit{Proof of Corollary \ref{corgateNOcond}.} Since Theorem \ref{mingatesNOcond} holds, the corollary follows by \cite[Theorem 5.4]{manzo2004essential}. $\qed$

\subsection{The minimal gates from a stable state to another stable state}\label{proofmingateappendix}

\textit{Proof of Theorem \ref{mingatessingh0}.} We recall that a subset $\mathcal W\subset\mathcal S(\mathbf r,\mathbf s)$ is a minimal gate for the transition from $\mathbf r$ to $\mathbf s$ if  
\begin{itemize}
\item[(i)] $\mathcal W$ is a gate for the transition from $\mathbf r$ to $\mathbf s$, i.e., every $\omega\in\Omega_{\mathbf r,\mathbf s}^{opt}$ intersects $\mathcal W$,
\item[(ii)] for any $\eta\in\mathcal W$ there exists an optimal path $\omega' \in\Omega_{\mathbf r,\mathbf s}^{opt}$ such that $\omega'\cap(\mathcal W\backslash\{\eta\})=\varnothing$.
\end{itemize}  
We prove that all the sets given in the items (a)--(d) satisfy (i) and (ii). Let us begin to consider the sets of item (a). First, we focus on the sets
\begin{align}\label{Pgatersfixed}
\bigcup_{\mathbf t\in\mathcal X^s\backslash\{\mathbf r\}} \overline{\mathscr P}(\mathbf r,\mathbf t)\ \ \ \text{and}\ \bigcup_{\mathbf t\in\mathcal X^s\backslash\{\mathbf r\}} \widetilde{\mathscr P}(\mathbf r,\mathbf t).
\end{align} Let $\omega\in\Omega_{\mathbf r,\mathbf s}^{opt}$ and let $\mathbf z\in\mathcal X^s\backslash\{\mathbf r\}$ be the first stable configuration visited by $\omega$ after $\mathbf r$. By Theorem \ref{mingatescond}(a), we get $\omega\cap \overline{\mathscr P}(\mathbf r,\mathbf z)\neq\varnothing$ and $\omega\cap\widetilde{\mathscr P}(\mathbf r,\mathbf z)\neq\varnothing$. Thus,
$$\omega\cap(\bigcup_{\mathbf t\in\mathcal X^s\backslash\{\mathbf r\}} \overline{\mathscr P}(\mathbf r,\mathbf t))\neq\varnothing\ \ \ \text{and}\ \ \ \omega\cap(\bigcup_{\mathbf t\in\mathcal X^s\backslash\{\mathbf r\}} \widetilde{\mathscr P}(\mathbf r,\mathbf t))\neq\varnothing$$
and (i) is satisfied. Let us now prove that for any $\eta\in\bigcup_{\mathbf t\in\mathcal X^s\backslash\{\mathbf r\}} \overline{\mathscr P}(\mathbf r,\mathbf t)$ there exists $\omega'\in\Omega_{\mathbf r,\mathbf s}^{opt}$ that satisfies (ii). There exists $\mathbf z\in\mathcal X^s\backslash\{\mathbf r\}$ such that $\eta\in \overline{\mathscr P}(\mathbf r,\mathbf z)$. Let $\omega^{(1)}:\mathbf r\to\mathbf z$ be the path constructed in the proof of Theorem \ref{mingatescond}(a) such that $\omega^{(1)}\cap\mathcal X^s\backslash\{\mathbf r,\mathbf z\}=\varnothing$. Let $\omega^{(2)}:\mathbf z\to\mathbf s$ be the reference path given in Definition \ref{remarkrepathNZ}. If $\mathbf z=\mathbf s$, then we define the optimal path $\omega'$ as the path $\omega^{(1)}$. Otherwise, we define $\omega'$ as the concatenation of the paths $\omega^{(1)}$ and $\omega^{(2)}$. Concerning $\bigcup_{\mathbf t\in\mathcal X^s\backslash\{\mathbf r\}} \widetilde{\mathscr P}(\mathbf r,\mathbf t)$, we argue similarly to prove (ii) using the symmetry of the model on $\Lambda$. In both cases, (ii) is verified and the sets in \eqref{Pgatersfixed} are minimal gates for the transition $\mathbf r\to\mathbf s$.\\
Second, let us focus on the sets 
\begin{align}\label{secondgatePrsfix}
\bigcup_{\mathbf t\in\mathcal X^s\backslash\{\mathbf s\}} \overline{\mathscr P}(\mathbf t,\mathbf s)\ \ \ \text{and}\ \bigcup_{\mathbf t\in\mathcal X^s\backslash\{\mathbf s\}} \widetilde{\mathscr P}(\mathbf t,\mathbf s).
\end{align}
 Let $\omega\in\Omega_{\mathbf r,\mathbf s}^{opt}$ and let $\mathbf z\in\mathcal X^s\backslash\{\mathbf s\}$ be the last stable configuration visited by $\omega$ before hitting $\mathbf s$. By Theorem \ref{mingatescond}(a), we get $\omega\cap\overline{\mathscr P}(\mathbf z,\mathbf s)\neq\varnothing$ and $\omega\cap\widetilde{\mathscr P}(\mathbf z,\mathbf s)\neq\varnothing$. Thus,
$$\omega\cap(\bigcup_{\mathbf t\in\mathcal X^s\backslash\{\mathbf s\}} \overline{\mathscr P}(\mathbf t,\mathbf s))\neq\varnothing\ \ \ \text{and}\ \ \ \omega\cap(\bigcup_{\mathbf t\in\mathcal X^s\backslash\{\mathbf s\}} \widetilde{\mathscr P}(\mathbf t,\mathbf s))\neq\varnothing$$
and (i) is satisfied. Let us now prove that for any $\eta\in\bigcup_{\mathbf t\in\mathcal X^s\backslash\{\mathbf s\}} \overline{\mathscr P}(\mathbf t,\mathbf s)$ there exists $\omega'\in\Omega_{\mathbf r,\mathbf s}^{opt}$ that satisfies (ii). There exists $\mathbf z\in\mathcal X^s\backslash\{\mathbf s\}$ such that $\eta\in\overline{\mathscr P}(\mathbf z,\mathbf s)$. Let $\omega^{(1)}:\mathbf r\to\mathbf z$ be the reference path given in Definition \ref{remarkrepathNZ}. Let $\omega^{(2)}:\mathbf z\to\mathbf s$ be
the path constructed in the proof of Theorem \ref{mingatescond}(a) such that $\omega^{(2)}\cap\mathcal X^s\backslash\{\mathbf z,\mathbf s\}=\varnothing$. If $\mathbf z=\mathbf r$, then we define the optimal path $\omega'$ as the path $\omega^{(2)}$. Otherwise, we define $\omega'$ as the concatenation of the paths $\omega^{(1)}$ and $\omega^{(2)}$. Concerning $\bigcup_{\mathbf t\in\mathcal X^s\backslash\{\mathbf s\}} \widetilde{\mathscr P}(\mathbf t,\mathbf s)$, we argue similarly to prove (ii) using the symmetry of the model on $\Lambda$. In both cases, (ii) is verified and the sets in \eqref{secondgatePrsfix} are minimal gates for the transition $\mathbf r\to\mathbf s$.\\ 
Following the same strategy, we prove that the sets in (b), (c) and (d) are minimal gates for the transition from $\mathbf r$ to $\mathbf s$ by using Theorem \ref{mingatescond}(b), (c), (d), respectively, instead of Theorem \ref{mingatescond}(a). See Appendix \ref{appendixproof2} for the details. $\qed$\\ 

\textit{Proof of Theorem \ref{setmingatefin}.} Our aim is to prove that $\mathcal G(\mathbf r,\mathbf s)$ only contains the minimal gates of Theorem \ref{mingatessingh0}. Thus, our strategy is to prove that any 
\begin{align}\label{lastunesssaddlefisso}
\eta\in\mathcal S(\mathbf r,\mathbf s)\backslash\biggl(\bigcup_{\mathbf t\in\mathcal X^s\backslash\{\mathbf r\}} \mathcal F(\mathbf r,\mathbf t)\cup\bigcup_{\mathbf t\in\mathcal X^s\backslash\{\mathbf s\}}{\mathcal F}(\mathbf t,\mathbf s)\biggr)
\end{align} 
does not belong to a minimal gate. To do this, note that 
\begin{align}
\mathcal S(\mathbf r,\mathbf s)&\backslash\biggl(\bigcup_{\mathbf t\in\mathcal X^s\backslash\{\mathbf r\}} \mathcal F(\mathbf r,\mathbf t)\cup\bigcup_{\mathbf t\in\mathcal X^s\backslash\{\mathbf s\}}{\mathcal F}(\mathbf t,\mathbf s)\biggr)=\notag\\&
\bigcup_{\mathbf t,\mathbf z\in\mathcal X^s\backslash\{\mathbf r,\mathbf s\}, \mathbf t\neq\mathbf z} \mathcal F(\mathbf t,\mathbf z)\cup\biggl(\mathcal S(\mathbf r,\mathbf s)\backslash\bigcup_{\mathbf t,\mathbf z\in\mathcal X^s, \mathbf t\neq\mathbf z} \mathcal F(\mathbf t,\mathbf z)\biggr).
\end{align}
The proof of the theorem is given by the following claims
\begin{enumerate}
\item $\bigcup_{\mathbf t,\mathbf z\in\mathcal X^s\backslash\{\mathbf r,\mathbf s\}, \mathbf t\neq\mathbf z} \mathcal F(\mathbf t,\mathbf z)$ is not a gate;
\item any $\eta\in\mathcal S(\mathbf r,\mathbf s)\backslash\bigcup_{\mathbf t,\mathbf z\in\mathcal X^s, \mathbf t\neq\mathbf z} \mathcal F(\mathbf t,\mathbf z)$ is an ussential saddle.
\end{enumerate}
Indeed, by the first claim we obtain that $\bigcup_{\mathbf t,\mathbf z\in\mathcal X^s\backslash\{\mathbf r,\mathbf s\}, \mathbf t\neq\mathbf z} \mathcal F(\mathbf t,\mathbf z)$ is not a minimal gate since it does not satisfy the condition (i) given at the beginning of the proof of Theorem \ref{mingatessingh0}. While by the second claim we obtain that any $\eta\in\mathcal S(\mathbf r,\mathbf s)\backslash\bigcup_{\mathbf t,\mathbf z\in\mathcal X^s, \mathbf t\neq\mathbf z} \mathcal F(\mathbf t,\mathbf z)$ is not a minimal gate by \cite[Theorem 5.1]{manzo2004essential}. \\
In order to prove the claim 1, it is enough to prove that there exists an optimal path $\omega\in\Omega_{\mathbf r,\mathbf s}^{opt}$ such that $\omega\cap\mathcal F(\mathbf t,\mathbf z)=\varnothing$ for any $\mathbf t,\mathbf z\in\mathcal X^s\backslash\{\mathbf r,\mathbf s\}, \mathbf t\neq\mathbf z$. We define this path as the reference path $\omega:\mathbf r\to\mathbf s$ given in Definition \ref{remarkrepathNZ}. \\
Let us now prove the claim 2. For any $\eta\in\mathcal S(\mathbf r,\mathbf s)\backslash\bigcup_{\mathbf t,\mathbf z\in\mathcal X^s, \mathbf t\neq\mathbf z} \mathcal F(\mathbf t,\mathbf z)$ and for any $\omega\in\Omega_{\mathbf r,\mathbf s}^{opt}$ such that $\eta\in\omega$, we have to show that both the conditions of Remark \ref{remuness} are verified. By Lemma \ref{strisciastriscia} and Lemma \ref{lemmaenBR} the condition (a) of Remark \ref{remuness} is satisfied. Next we move to prove condition (b). Let $\mathbf t_1,\dots,\mathbf t_{m-1}\in\mathcal X^s$ be the stable configurations visited by $\omega$ in $\mathcal X^s$ before hitting $\mathbf s$. If we set $\mathbf t_0=\mathbf r$, $\mathbf t_m=\mathbf s$ and $\mathbf t_i\neq\mathbf t_{i+1}$ for all $i=0,\dots,m-1$, $m\in\mathbb N$, we rewrite $\omega$ as the concatenation of the $m$ paths $\omega^{(i)}:\ \mathbf t_i\to\mathbf t_{i+1}$. Let us assume that $\eta\in\omega^{(j)}$. Let $\omega'^{(j)}\in\Omega_{\mathbf t_j,\mathbf t_{j+1}}^{opt}$ be the optimal path constructed in the proof of Theorem \ref{mingatescondset} such that \[\{\text{argmax}_{\omega'^{(j)}} H\}\subseteq\{\text{argmax}_{\omega^{(j)}} H\}\backslash\{\eta\}\] and $\omega'^{(j)}\cap\mathcal X^s\backslash\{\mathbf t_j,\mathbf t_{j+1}\}=\varnothing$. Thus we can define a path $\omega'$ such that \[\{\text{argmax}_{\omega'} H\}\subseteq\{\text{argmax}_{\omega} H\}\backslash \{\eta\}\] as the concatenation of the $m$ paths $\omega^{(1)},\dots,\omega^{(j-1)},\dots,\omega'^{(j)},\omega^{(j+1)},\dots,\omega^{(m)}$. Hence, both the conditions of Remark \ref{remuness} are satisfied and $\eta$ is an unessential saddle.
$\qed$\\

\textit{Proof of Corollary \ref{corgateconv2}.} Since Theorem \ref{mingatessingh0} holds, the corollary follows by \cite[Theorem 5.4]{manzo2004essential}. $\qed$

\section{Restricted-tube and tube of typical paths}\label{proofmaintube}
In this section we prove the main results on the restricted-tube of typical paths and on the tube of typical paths stated in Section \ref{mainrestube}. 
\begin{remark}\label{remprinbound}
Given a $q$-Potts configuration $\sigma\in\mathcal X$ on a grid-graph $\Lambda$, a vertex $v\in V$ and a spin value $s\in\{1,\dots,q\}$ such that $\sigma(v)\neq s$, using \eqref{energydiff} we have \[H(\sigma)-H(\sigma^{v,s})\in\{-4,-2,0,2,4\}.\]
It follows that we can depict  the principal boundary of an extended cycles $\mathcal C$ in \eqref{defKbar}--\eqref{defEi} by the union of those configurations $\bar\sigma\in\partial\mathcal C$ such that either
\begin{itemize}
\item[(i)] $H(\bar\sigma)-H(\sigma)=-2,$ or
\item[(ii)]  $H(\bar\sigma)-H(\sigma)=-4.$
\end{itemize}
\end{remark}
Let us now state some useful lemmas which we prove at the end this Section by using Remark \ref{remprinbound}. 

\begin{lemma}\label{lemmaprinboundHK}
Let $\mathbf r,\mathbf s\in\mathcal{X}^s$, $\mathbf r\neq\mathbf s$. Then,
\begin{align}\label{prinlabelal}
\mathcal B(\overline{\mathcal K}(r,s))=\overline{\mathcal D}_1(r,s)\cup\overline{\mathcal D}_2(r,s)\cup\overline{\mathcal E}_1(r,s)\cup\overline{\mathcal E}_2(r,s).
\end{align}
\end{lemma}
The principal boundary \eqref{prinlabelal} is illustrated in the first row of Figure \ref{tubecycle}. 
\begin{lemma}\label{lemmaprinboundDE}
Let $\mathbf r,\mathbf s\in\mathcal{X}^s$, $\mathbf r\neq\mathbf s$. Then, for any $i=1,\dots,K-4$, 
\begin{align}\label{prinlabelal3}
&\mathcal B(\overline{\mathcal D}_i(r,s))=\overline{\mathcal D}_{i+1}(r,s)\cup\overline{\mathcal D}_{i+2}(r,s)\cup\overline{\mathcal E}_{i+1}(r,s)\cup\overline{\mathcal E}_{i+2}(r,s),\\
\label{prinlabelal2}
&\mathcal B(\overline{\mathcal E}_i(r,s))=\overline{\mathcal E}_{i+1}(r,s)\cup\overline{\mathcal E}_{i+2}(r,s).
\end{align}
\end{lemma}
The principal boundaries \eqref{prinlabelal3} and \ref{prinlabelal2} are illustrated in the middle rows of Figure \ref{tubecycle}. 
\begin{lemma}\label{lemmaprinboundfinal}
Let $\mathbf r,\mathbf s\in\mathcal{X}^s$, $\mathbf r\neq\mathbf s$. Then,
\begin{align}
&\mathcal B(\overline{\mathcal D}_{K-3}(r,s))=\overline{\mathcal D}_{K-2}(r,s)\cup\overline{\mathcal E}_{K-2}(r,s)\cup\bar R_{1,1}(\mathbf r,\mathbf s),\\
&\mathcal B(\overline{\mathcal E}_{K-3}(r,s))=\overline{\mathcal E}_{K-2}(r,s)\cup\bar R_{1,1}(r, s),\\
\label{claimlastlemma}
&\mathcal B(\overline{\mathcal D}_{K-2}(r,s))=\mathcal B(\overline{\mathcal E}_{K-2}(r,s))=\bar R_{1,1}(r, s).
\end{align}
\end{lemma}
 The principal boundaries described in Lemma \ref{lemmaprinboundfinal} are illustrated in the last rows of Figure \ref{tubecycle}. \\

First we describe the ``restricted-tube'' of typical paths, then we exploit this result to describe the tube of typical paths for the transition from a stable state to the other stable configurations and for the transition from a stable state to another stable configuration.  \\

\textit{Proof of Theorem \ref{exittuberes}.} We want to prove that $\mathscr U_{\mathbf s}({\mathbf r})$ in \eqref{aligntuberes} satisfies the following properties: it includes $\mathcal C\in\mathcal M(\mathcal C_{\mathbf s}^+(\mathbf r)\backslash\mathbf s)$ that belong to at least a cycle-path $(\mathcal C_1,\dots,\mathcal C_n)\in J_{\mathbf r,\mathbf s}$, $n\in\mathbb N$, such that $\bigcup_{i=1}^n \mathcal C_i\cap\mathcal X^s\backslash\{\mathbf r,\mathbf s\}=\varnothing$, $\mathcal C_1=\mathcal C_{\mathbf s}(\mathbf r)$ and $\mathbf s\in\partial\mathcal C_n$, see \eqref{deftubemodind}. \\

We start by studying the first descent from a trivial cycle $\{\xi^*\}$ for some $\xi^*\in\bar R_{\lfloor\frac{L}{2}\rfloor,K}(r,s)$ to $\mathbf r$, where $\lfloor n \rfloor:=\max\{m\in\mathbb Z:\ m\le n\}$. Using the symmetry of the model on $\Lambda$, we describe similarly the first descent from the same configuration $\xi^*$ to $\mathbf s$. Finally using reversibility, we will obtain a complete description of $\mathscr U_{\mathbf s}({\mathbf r})$ by joining the time reversal of the first descent from $\{\xi^*\}$ to $\mathbf r$ with the first discent from $\{\xi^*\}$ to $\mathbf s$.\\

For sake of semplicity, we separate the description of the first descent from $\{\xi^*\}$ for some $\xi^*\in\bar R_{\lfloor\frac{L}{2}\rfloor,K}(r,s)$ to $\mathbf r$ in more parts. We start by studying the typical trajectories followed by the process during the transition from $\{\xi^*\}$ to $\bar R_{2,K}(\mathbf r,\mathbf s)\subset\partial\overline{\mathcal K}(r,s)$, see \eqref{defKbar}, and then we study the typical paths followed for the first descent from $\overline{\mathcal K}(r,s)$ to $\mathbf r$. It is useful to remark that $\partial\mathcal C_{\mathbf s}(\mathbf r)\cap\overline{\mathcal K}(r,s)\neq\varnothing$. 
Using Lemma \ref{strisciastriscia}(a) and (c), for any $i=\lfloor\frac{L}{2}\rfloor-1,\dots,2\ $ we define a cycle-path $(\mathcal C_i^0,\mathcal C_i^1,\mathcal C_i^2)$ such that  
\begin{itemize}
\item[-] $\mathcal C_{i}^{0}=\{\eta_K\}$ for $\eta_K\in\bar R_{i+1,K}$, 
\item[-] $\mathcal C_{i}^{1}=\bigcup_{j=1}^{K-1} \{\eta_j\}$ for $\eta_j\in\bar B_{i,K}^j$,
\item[-]$\mathcal C_{i}^{2}=\{\eta_0\}$ for $\eta_0\in\bar R_{i,K}$,
\end{itemize}
where $\eta_{K-1},\dots,\eta_0$ are chosen in such a way that there exists $v\in V$ such that $\eta_i:=\eta_{i+1}^{v,r}$. Note that $\mathcal C_{i}^{0}$, $\mathcal C_{i}^{2}$ are non trivial cycles. For any $i=\lfloor\frac{L}{2}\rfloor-1,\dots,2$, using Lemma \ref{strisciastriscia}, we also remark that  for any $\sigma\in\bar R_{i,K}$, $\mathscr F(\partial\{\sigma\})\subset\bar B_{i-1,K}^{K-1}\cup\bar B_{i,K}^{1}$.
Then, using Lemma \ref{lemmaenBR} we have that $\{\sigma\}$ satisfies \eqref{cycle} and that
\begin{align}
\mathcal B(\{\sigma\})=\mathscr F(\partial\{\sigma\}).
\end{align}
Moreover, using Lemma \ref{lemmaenBR} and Lemma \ref{strisciastriscia} we remark that for any $i=\lfloor\frac{L}{2}\rfloor-1,\dots,2$, $\mathcal C_i^1$ is a plateau and its principal boundary is given by
\begin{align}
\mathcal B(\mathcal C_i^1)=\mathcal C_i^0\cup\mathcal C_i^2.
\end{align}
Thus, $(\mathcal C_i^0,\mathcal C_i^1,\mathcal C_i^2)\in J_{\mathcal C_i^0,\mathcal C_i^2}$ since \eqref{cyclepathvtj} is satisfied.
Hence, starting from $\mathcal C_{{\lfloor\frac{L}{2}\rfloor}}^0=\{\xi^*\}$ for some $\xi^*\in\bar R_{\lfloor\frac{L}{2}\rfloor,K}(r,s)$, we depict a cycle-path vtj-connected to $\mathcal C_2^2=\{\hat\eta\}$ for some appropriate $\hat\eta\in\bar R_{2,K}$ as
\begin{align}\label{cyclepathpart1}
(\mathcal C_{\lfloor\frac{L}{2}\rfloor}^0,\mathcal C_{\lfloor\frac{L}{2}\rfloor-1}^1,\mathcal C_{\lfloor\frac{L}{2}\rfloor-1}^2\equiv\mathcal C_{\lfloor\frac{L}{2}\rfloor-1}^0,\dots,\mathcal C_3^2\equiv\mathcal C_2^0,\mathcal C_2^1,\mathcal C_2^2).
\end{align}
Let us now to describe a cycle-path 
\begin{align}\label{cyclepathfinal1}
(\bar{\mathcal C}_1,\dots,\bar{\mathcal C}_m)
\end{align} vtj-connected to $\mathbf r$ such that $\bar{\mathcal C}_1=\{\bar\eta^*\}$, with $\bar\eta^*\in\overline{\mathscr P}(\mathbf r,\mathbf s)$ which exists in view of Lemma \ref{intPQH} and which is chosen in such a way that it is defined by a spin update in a vertex of $\hat\eta$. Let us begin to note that any set from \eqref{defKbar} to \eqref{defEi} is an extended cycle, i.e., a maximal connected set of equielevated trivial cycles. 
Thanks to Lemmas \ref{lemmaprinboundHK}, \ref{lemmaprinboundDE} and \ref{lemmaprinboundfinal} we obtain that for every $i=m,\dots,n-1$, the cycle-path \eqref{cyclepathfinal1} is characterized by a sequence of cycles and extended cycles $\bar{\mathcal C}_1,\dots,\bar{\mathcal C}_m$ such that \[\bar{\mathcal C}_1,\dots,\bar{\mathcal C}_m\subset\overline{\mathcal K}(r,s)\cup\bigcup_{i=1}^{K-2} (\overline{\mathcal D}_i(r,s)\cup\overline{\mathcal E}_i(r,s))\cup\bar R_{1,1}(r,s)\]
and $\bar{\mathcal C}_m=\mathbf r$.
More precisely, we can say that $(\bar{\mathcal C}_1,\dots,\bar{\mathcal C}_m)\in J_{\bar{\mathcal C}_1,\bar{\mathcal C}_m}$.  Indeed, using Lemmas \ref{lemmaprinboundHK}, \ref{lemmaprinboundDE} and \ref{lemmaprinboundfinal} it follows that for any pair of consecutive cycles the condition \eqref{cyclepathvtj} is satisfied.\\
Arguing similarly, we construct a cycle-path vtj-connected from $\{\xi^*\}$ to $\mathbf s$. Indeed, first we construct the following 
\begin{align}\label{cyclepathpart2}
(\mathcal C_{\lfloor\frac{L}{2}\rfloor}^0,\mathcal C_{\lfloor\frac{L}{2}\rfloor+1}^1,\mathcal C_{\lfloor\frac{L}{2}\rfloor+1}^2\equiv\mathcal C_{\lfloor\frac{L}{2}\rfloor+1}^0,\dots,\mathcal C_{L-3}^2\equiv\mathcal C_{L-2}^0,\mathcal C_{L-2}^1,\mathcal C_{L-2}^2)\in J_{\{\xi^*\},\{\check\eta\}}
\end{align}
from $\{\xi^*\}$ to $\{\check\eta\}$ for some $\check\eta\in\bar R_{L-2,K}(r,s)$. Then, we define
\begin{align}\label{cyclepathfinal2}
(\tilde{\mathcal C}_1,\dots,\tilde{\mathcal C}_m)\in J_{\tilde{\mathcal C}_1,\tilde{\mathcal C}_m}
\end{align}
 for the first descent from $\{\tilde\eta^*\}$ to $\mathbf s$, where $\tilde\eta^*\in\widetilde{\mathscr P}(\mathbf r,\mathbf s)$ exists in view of Lemma \ref{intPQH} and it is chosen in such a way that it is defined by a spin update in a vertex of $\check\eta$. More precisely, using Lemma \ref{strisciastriscia}(b) and (d), for any $i=\lfloor\frac{L}{2}\rfloor+1,\dots,L-2\ $ we define a cycle-path $(\mathcal C_i^0,\mathcal C_i^1,\mathcal C_i^2)$ such that  
\begin{itemize}
\item[-] $\mathcal C_{i}^{0}=\{\eta_0\}$ for $\eta_0\in\bar R_{i-1,K}$, 
\item[-] $\mathcal C_{i}^{1}=\bigcup_{j=1}^{K-1} \{\eta_j\}$ for $\eta_j\in\bar B_{i,K}^j$,
\item[-]$\mathcal C_{i}^{2}=\{\eta_K\}$ for $\eta_K\in\bar R_{i,K}$,
\end{itemize}
where $\eta_1,\dots,\eta_K$ are chosen in such a way that there exists $v\in V$ such that $\eta_i:=\eta_{i-1}^{v,r}$. On the other hand, thanks to the symmetry of the model on $\Lambda$, we prove the analogue of Lemmas \ref{lemmaprinboundHK}, \ref{lemmaprinboundDE} and \ref{lemmaprinboundfinal} for $\widetilde{\mathcal K}(r,s), \widetilde{\mathcal D}_i(r,s), \widetilde{\mathcal E}_i(r,s)$ for any $i=1,\dots,K-2$. Thus, we construct the cycle-path vtj-connected $(\tilde{\mathcal C}_1,\dots,\tilde{\mathcal C}_m)$ from $\{\tilde\eta^*\}$ to $\mathbf s$ with  
 \[\tilde{\mathcal C}_1,\dots,\tilde{\mathcal C}_m\subset\widetilde{\mathcal K}(r,s)\cup\bigcup_{i=1}^{K-2} (\widetilde{\mathcal D}_i(r,s)\cup\widetilde{\mathcal E}_i(r,s))\cup\tilde R_{1,1}(r,s).\]
Finally, using reversibility and in view of the above construction, we construct a cycle-path vtj-connected from $\mathbf r$ to $\mathbf s$  as
\begin{align}\label{cyclepathvtjcomplete}
 (\bar{\mathcal C}_m&\equiv\mathbf r,\dots,\bar{\mathcal C}_1, \mathcal C_2^2, \mathcal C_2^1, \mathcal C_3^2\equiv\mathcal C_2^0, \dots,\mathcal C_{\lfloor\frac{L}{2}\rfloor-1}^2\equiv\mathcal C_{\lfloor\frac{L}{2}\rfloor-1}^0, \mathcal C_{\lfloor\frac{L}{2}\rfloor-1}^1, \mathcal C_{\lfloor\frac{L}{2}\rfloor}^0\equiv\{\xi^*\}, \notag  \\ &\mathcal C_{\lfloor\frac{L}{2}\rfloor+1}^1,\mathcal C_{\lfloor\frac{L}{2}\rfloor+1}^2\equiv\mathcal C_{\lfloor\frac{L}{2}\rfloor+1}^0,\dots,\mathcal C_{L-3}^2\equiv\mathcal C_{L-2}^0,\mathcal C_{L-2}^1,\mathcal C_{L-2}^2,\tilde{\mathcal C}_1,\dots,\tilde{\mathcal C}_m\equiv\mathbf s).
\end{align}

Hence, we conclude. Indeed, according to the discussion above and using \eqref{cyclepathvtjcomplete}, we depict the restricted-tube of typical paths between $\mathbf r$ and $\mathbf s$ as in \eqref{aligntuberes}. Instead, \eqref{timetuberes} follows by \cite[Lemma 3.13]{nardi2016hitting}. $\qed$\\

We are now ready to prove Theorems \ref{exittubestable} and \ref{exittubetarget}.\\

\textit{Proof of Theorem \ref{exittubestable}.} Using the description of the restricted-tube of typical paths between two stable states, we depict the tube $\mathfrak T_{\mathcal X^s\backslash\{\mathbf r\}}(\mathbf r)$ as in \eqref{aligntubestable}. Indeed, any cycle-path $(\mathcal C_1,\dots,\mathcal C_n)$, $n\in\mathbb N$, vtj-connected to $\mathcal X^s\backslash\{\mathbf r\}$ starting from $\mathbf r$ is described by a sequence of cycles belonging to some $\mathscr U_{\mathbf s}(\mathbf r)$ where $\mathbf s$ is the first stable configuration visited by that $(\mathcal C_1,\dots,\mathcal C_n)$  in $\mathcal X^s\backslash\{\mathbf r\}$. Instead, \eqref{timetubestable} follows by \cite[Lemma 3.13]{nardi2016hitting}. $\qed$\\

\textit{Proof of Theorem \ref{exittubetarget}.} As we have already noted in Section \ref{subgatetarget}, starting from $\mathbf r$ the process could intersect some stable configurations different from the target $\mathbf s$. It follows that, when we study the typical trajectories for the transition $\mathbf r\to\mathbf s$ we have to take into account the possible intermediate transitions between two stable states and such that do not intersect other stable configurations beyond them. For any of these intermediate transitions, we apply the result of the restricted-case. Hence, we depict $\mathfrak T_{\mathbf s}(\mathbf r)$ as in \eqref{aligntubetarget}. Instead, \eqref{timetubetarget} follows by \cite[Lemma 3.13]{nardi2016hitting}. $\qed$\\

\textbf{Proof of Lemmas \ref{lemmaprinboundHK}, \ref{lemmaprinboundDE}, \ref{lemmaprinboundfinal}}\label{prooflemmas}\\

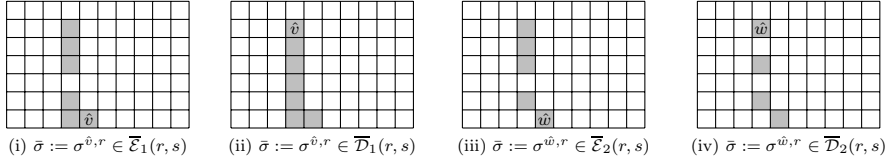
\begin{figure}[h!]
\begin{minipage}[c]{1\textwidth}
    \centering
    \makebox[0pt]{%
\begin{tikzpicture} [scale=0.8, transform shape]
\fill[lightgray] (0.9,0)rectangle(1.5,0.3) (0.9,0.3)rectangle(1.2,0.6) (0.9,0.9)rectangle(1.2,1.8);
\draw (1.35,0.15) node{\footnotesize{$\hat v$}};
\draw (1.5,0) node[below] {\footnotesize{(i) $\bar\sigma:=\sigma^{\hat v,r}\in\overline{\mathcal E}_1(r,s)$}};
\draw[step=0.3cm,color=black] (0,0) grid (3,2.1);
\end{tikzpicture}\ \ \
\begin{tikzpicture}[scale=0.8, transform shape]
\fill[lightgray] (0.9,0)rectangle(1.5,0.3) (0.9,0.3)rectangle(1.2,1.8);
\draw (1.05,1.65) node{\footnotesize{$\hat v$}};
\draw (1.5,0)  node[below] {\footnotesize{(ii) $\bar\sigma:=\sigma^{\hat v,r}\in\overline{\mathcal D}_1(r,s)$}};
\draw[step=0.3cm,color=black] (0,0) grid (3,2.1);
\end{tikzpicture}\ \ \
\begin{tikzpicture}[scale=0.8, transform shape]
\fill[lightgray] (1.2,0)rectangle(1.5,0.3) (0.9,0.3)rectangle(1.2,0.6) (0.9,0.9)rectangle(1.2,1.8);
\draw (1.35,0.15) node{\footnotesize{$\hat w$}};
\draw (1.5,0) node[below] {\footnotesize{(iii) $\bar\sigma:=\sigma^{\hat w,r}\in\overline{\mathcal E}_2(r,s)$}};
\draw[step=0.3cm,color=black] (0,0) grid (3,2.1);
\end{tikzpicture}\ \ \
\begin{tikzpicture}[scale=0.8, transform shape]
\fill[lightgray] (1.2,0)rectangle(1.5,0.3) (0.9,0.3)rectangle(1.2,0.6) (0.9,0.9)rectangle(1.2,1.2)(0.9,1.5)rectangle(1.2,1.8);
\draw (1.05,1.65) node{\footnotesize{$\hat w$}};
\draw (1.5,0) node[below] {\footnotesize{(iv) $\bar\sigma:=\sigma^{\hat w,r}\in\overline{\mathcal D}_2(r,s)$}};
\draw[step=0.3cm,color=black] (0,0) grid (3,2.1);
\end{tikzpicture}%
    }\par
     \end{minipage}
\caption{\label{examplesprincbound} Examples of $\sigma\in\overline{\mathcal K}(r,s)$ and $\bar\sigma\in\mathcal B(\overline{\mathcal K}(r,s)$. We color white the vertices with spin $r$ and gray those vertices whose spin is $s$.}
\end{figure}

\textit{Proof of Lemma  \ref{lemmaprinboundHK}} According to \eqref{principalboundary}, we describe the principal boundary of the extended cycle $\overline{\mathcal K}(r,s)$ by looking for those configurations $\bar\sigma\notin\overline{\mathcal K}(r,s)$ which communicate with some $\sigma\in\overline{\mathcal K}(r,s)$ such that $\sigma$ and $\bar\sigma$ satisfy either case (i) or case (ii) of Remark \ref{remprinbound}.\\
Let us start to consider case Remark \ref{remprinbound}(i). In view of \eqref{defKbar}, this case occurs only when $\sigma$ has a spin $s$ with three nearest-neighbor spins $r$ and $\bar\sigma$ is obtained from $\sigma$ by flipping from $s$ to $r$ this spin $s$. In particular, we note that $\sigma\in\overline{\mathcal K}(r,s)$  has a spin $s$ with three nearest-neighbor spins $r$ only when $\sigma\in\overline{\mathcal K}(r,s)\backslash[\overline{\mathcal Q}(\mathbf r,\mathbf s)\cup\overline{\mathscr P}(\mathbf r,\mathbf s)]$, i.e., when
\begin{align}
\sigma\in\overline{\mathscr H}(\mathbf r,\mathbf s)\cup&\{\sigma\in\mathcal X: \sigma(v)\in\{r,s\}\ \forall v\in V, H(\sigma)=2K+2+H(\mathbf r), \sigma\ \text{has at}\notag\\ &\text{ least two $s$-interacting clusters},\ \text{and}\ R(\text C^s(\sigma))=R_{2\times(K-1)}\}.\notag
\end{align}
Hence, consider $\sigma\in\overline{\mathcal K}(r,s)$ with a spin $s$, say on vertex $\hat v$, with three nearest-neighbor spins $r$ and one nearest-neighbor spin $s$ and define $\bar\sigma:=\sigma^{\hat v,r}$. We note that
\begin{itemize}
\item[-] for any $h=2,\dots,K-2$, if $\sigma\in\bar B_{1,K-1}^h(r,s)$, then $\bar\sigma\in\bar B_{1,K-1}^{h-1}(r,s)\subset\overline{\mathcal D}_1(r,s)$;
\item[-] if $\sigma\in\bar B_{1,K}^1(r,s)$, then $\bar\sigma\in\bar R_{1,K}(r,s)\subset\overline{\mathcal E}_1(r,s)$;
\item[-] if $\hat v$ has its unique nearest-neighbor spin $s$ on an adjacent column, then $\bar\sigma\in\overline{\mathcal E}_1(r,s)$, see Figure \ref{examplesprincbound}(i);
\item[-] if $\hat v$ and its nearest-neighbor spin $s$ lie on the same column, then $\bar\sigma\in\overline{\mathcal D}_1(r,s)$,  see Figure \ref{examplesprincbound}(ii).
\end{itemize}

\FloatBarrier
Next we move to consider case Remark \ref{remprinbound}(ii), that occurs only when $\sigma$ has a spin $s$, say on vertex $\hat w$, sourrounded by four spins $r$ and $\bar\sigma:=\sigma^{\hat w,r}$. This happens when \[\sigma\in\overline{\mathcal K}(r,s)\backslash[\overline{\mathscr H}(\mathbf r,\mathbf s)\cup\overline{\mathcal Q}(\mathbf r,\mathbf s)\cup{\mathscr P}(\mathbf r,\mathbf s)].\]
Then, we note that
\begin{itemize}
\item[-] if $R(\text C^s(\bar\sigma))=R_{1\times(K-2)}$, i.e., if $\hat w$ lies on a column where there are no other spins $s$, then $\bar\sigma\in\overline{\mathcal E}_2(r,s)$, see Figure \ref{examplesprincbound}(iii);
\item[-] if $R(\text C^s(\bar\sigma))=R_{2\times(K-3)}$, i.e., if $\hat w$ lies on a column where there are other spins $s$, then $\bar\sigma\in\overline{\mathcal D}_2(r,s)$, see Figure \ref{examplesprincbound}(iv).
\end{itemize}
$\qed$\\

\textit{Proof of Lemma \ref{lemmaprinboundDE}} Similarly to the proof of Lemma \ref{lemmaprinboundHK}, for any $i=1,\dots,K-4$ we describe the principal boundary of the extended cycles $\overline{\mathcal D}_i(r,s)$ and $\overline{\mathcal E}_i(r,s)$ by using Remark \ref{remprinbound}.\\
For any $i=1,\dots,K-4$, let us start to study the principal boundary of $\overline{\mathcal D}_i(r,s)$. Case Remark \ref{remprinbound}(i) occurs only when $\sigma\in\overline{\mathcal D}_i(r,s)$ has a spin $s$, say on vertex $\hat v$ with three nearest-neighbor spins $r$ and one nearest-neighbor spin $s$ and $\bar\sigma:=\sigma^{\hat v,r}$. Then, we note that
\begin{itemize}
\item[-] if $\hat v$ its unique nearest-neighbor spin $s$ on an adjacent column, then $R(\text C^s(\bar\sigma))=R_{1\times(K-(i+1))}$ and $\bar\sigma\in\overline{\mathcal E}_{i+1}(r,s)$, see Figure \ref{secondexamplesprincbound}(ii);
\item[-] if $\hat v$ and its nearest-neighbor spin $s$ lie on the same column, then $R(\text C^s(\bar\sigma))=R_{2\times(K-(i+2))}$ and $\bar\sigma\in\overline{\mathcal D}_{i+1}(r,s)$,  see Figure \ref{secondexamplesprincbound}(i).
\end{itemize}

Regarding case Remark \ref{remprinbound}(ii), it occurs only when $\sigma$ has a spin $s$, say on vertex $\hat w$, sourrounded by four spins $r$ and $\bar\sigma:=\sigma^{\hat w,r}$. It follows that, 
\begin{itemize}
\item[-] if $R(\text C^s(\bar\sigma))=R_{1\times(K-(i+2))}$, i.e., if $\hat w$ lies on a column where there are no other spins $s$, then $\bar\sigma\in\overline{\mathcal E}_{i+2}(r,s)$,  see Figure \ref{secondexamplesprincbound}(i);
\item[-] if $R(\text C^s(\bar\sigma))=R_{2\times(K-(i+3))}$, i.e., if $\hat w$ lies on a column where there are other spins $s$, then $\bar\sigma\in\overline{\mathcal D}_{i+2}(r,s)$, see Figure \ref{secondexamplesprincbound}(ii).
\end{itemize} 

For any $i=1,\dots,K-4$, next we move to describe the principal boundary of $\overline{\mathcal E}_{i}(r,s)$. Case Remark \ref{remprinbound}(i) occurs when $\sigma\in\overline{\mathcal E}_{i}(r,s)$ has a spin $s$, say on vertex $\hat v$, with three nearest-neighbor spins $r$ and one nearest-neighbor spin $s$ and $\bar\sigma:=\sigma^{\hat v,r}$. Hence, $\bar\sigma\in\overline{\mathcal E}_{i+1}(r,s)$, see Figure \ref{secondexamplesprincbound}(ii). Instead case Remark \ref{remprinbound}(ii) occurs only when $\sigma$ has a spin $s$, say on vertex $\hat w$, sourrounded by four spins $r$ and $\bar\sigma:=\sigma^{\hat w,r}$. Thus $\bar\sigma\in\overline{\mathcal E}_{i+2}(r,s)$, see Figure \ref{secondexamplesprincbound}(iv). $\qed$\\

\begin{figure}[h!]
\begin{minipage}[c]{1\textwidth}
    \centering
    \makebox[0pt]{%
\begin{tikzpicture}[scale=0.8, transform shape]
\fill[lightgray] (1.2,0)rectangle(1.5,0.3) (0.9,0.3)rectangle(1.2,1.2);
\draw (1.05,1.05) node{\footnotesize{$\hat v$}};
\draw (1.35,0.15) node{\footnotesize{$\hat w$}};
\draw (1.5,0) node[below] {\footnotesize{(i) $\bar\sigma:=\sigma^{\hat v,r}\in\overline{\mathcal D}_3(r, s)$}};
\draw (1.78,-0.4) node[below] {\footnotesize{$\bar\sigma:=\sigma^{\hat w,r}\in\overline{\mathcal E}_4(r,s)$}};
\draw[step=0.3cm,color=black] (0,0) grid (3,2.1);
\end{tikzpicture}\ \ \
\begin{tikzpicture}[scale=0.8, transform shape]
\fill[lightgray] (0.9,0)rectangle(1.5,0.3) (0.9,0.3)rectangle(1.2,0.6) (0.9,0.9)rectangle(1.2,1.2);
\draw (1.05,1.05) node{\footnotesize{$\hat w$}};
\draw (1.35,0.15) node{\footnotesize{$\hat v$}};
\draw (1.5,0) node[below]  {\footnotesize{(ii) $\bar\sigma:=\sigma^{\hat v,r}\in\overline{\mathcal E}_3(r,s)$}};
\draw (1.78,-0.4) node[below]{\footnotesize{ $\bar\sigma:=\sigma^{\hat w,r}\in\overline{\mathcal D}_4(r,s)$}};
\draw[step=0.3cm,color=black] (0,0) grid (3,2.1);
\end{tikzpicture}\ \ \
\begin{tikzpicture}[scale=0.8, transform shape]
\fill[lightgray] (0.9,0)rectangle(1.2,0.6) (0.9,0.9)rectangle(1.2,1.5);
\draw (1.05,0.15) node{\footnotesize{$\hat v$}};
\draw (1.5,0) node[below] {\footnotesize{(iv) $\bar\sigma:=\sigma^{\hat v,r}\in\overline{\mathcal E}_3(r,s)$}};
\draw (1.78,-0.4) node[below, white] {\footnotesize{$\bar\sigma:=\sigma^{w,r}\in\overline{\mathcal E}_3(r,s)$}};
\draw[step=0.3cm,color=black] (0,0) grid (3,2.1);
\end{tikzpicture}\ \ \
\begin{tikzpicture}[scale=0.8, transform shape]
\fill[lightgray] (0.9,0)rectangle(1.2,0.9) (0.9,1.2)rectangle(1.2,1.5);
\draw (1.05,1.35) node{\footnotesize{$\hat w$}};
\draw (1.5,0) node[below] {\footnotesize{(iii) $\bar\sigma:=\sigma^{\hat w,r}\in\overline{\mathcal E}_4(r,s)$}};
\draw (1.78,-0.4) node[below,white] {\footnotesize{$\bar\sigma:=\sigma^{w,r}\in\overline{\mathcal E}_4(r,s)$}};
\draw[step=0.3cm,color=black] (0,0) grid (3,2.1);
\end{tikzpicture}%
    }\par
     \end{minipage}
\caption{\label{secondexamplesprincbound} Examples of $\sigma\in\overline{\mathcal D}_2(r,s)$ and $\bar\sigma\in\mathcal B(\overline{\mathcal D}_2(r,s))$ in (i) and (ii); examples of $\sigma\in\overline{\mathcal E}_2(r,s)$ and $\bar\sigma\in\mathcal B(\overline{\mathcal E}_2(r,s))$ in (iii) and (iv). We color white the vertices with spin $r$ and gray those vertices whose spin is $s$.}
\end{figure}
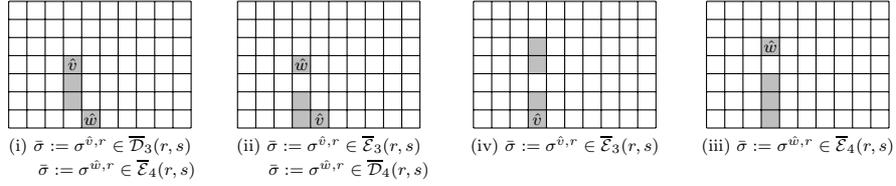\FloatBarrier

\textit{Proof of Lemma \ref{lemmaprinboundfinal}} First of all, we note that from \eqref{defDi} and \eqref{defEi} we have
\begin{align}\label{setKmeno2}
\overline{\mathcal D}_{K-2}(r,s)=\bar R_{2,1}(r, s)\ \ \ \text{and}\ \ \ \overline{\mathcal E}_{K-2}(r,s)=\bar R_{1,2}(r,s).
\end{align}
For $i=K-3,K-2$ once again we describe the principal boundary of the extended cycles $\overline{\mathcal D}_i(r,s)$ and $\overline{\mathcal E}_i(r,s)$ by using Remark \ref{remprinbound}. Let us begin to study the principal boundary of $\overline{\mathcal D}_{K-3}(r,s)$. Case Remark \ref{remprinbound}(i) takes place if $\sigma$ has a spin $s$, say on vertex $\hat v$, with three nearest-neighbor spins $r$ and one nearest-neighbor spin $s$ and $\bar\sigma:=\sigma^{\hat v,r}$. Hence, it occurs only when $\sigma\in\bar B_{1,2}^1(r, s)$ and either $\bar\sigma\in\overline{\mathcal D}_{K-2}(r,s)$ or $\bar\sigma\in\overline{\mathcal E}_{K-2}(r,s)$, see Figure \ref{thirdexamplesprincbound}(i) and (ii). 
Instead case Remark \ref{remprinbound}(ii) occurs when $\sigma$ has only two spins $s$ and they lie on the diagonal of a rectangle $R_{2\times 2}$, i.e., when $\bar\sigma$ is obtained by flipping from $s$ to $r$ one of these two spins $s$ and $\bar\sigma\in\bar R_{1,1}(r,s)$, see Figure \ref{thirdexamplesprincbound}(iii).\\

Next we move to describe the principal boundary of $\overline{\mathcal E}_{K-3}(r,s)$. Case Remark \ref{remprinbound}(ii) occurs when $\sigma\in\overline{\mathcal E}_{K-3}(r,s)$ has a spin $s$, say on vertex $\hat v$, with three nearest-neighbor spins $r$ and one nearest-neighbor spin $s$ and $\bar\sigma:=\sigma^{\hat v,r}$.  Hence, when $\sigma\in\bar R_{1,3}(r,s)$ and $\bar\sigma\in\overline{\mathcal E}_{K-2}(r,s)$. Finally, case Remark \ref{remprinbound}(ii) is verified when $\sigma\in\overline{\mathcal E}_{K-3}(r,s)$ has two spins $s$ with four nearest-neighbor spins $r$ and one of them is flipped to $r$, i.e., when $\bar\sigma\in\bar R_{1,1}(r,s)$, see Figure \ref{thirdexamplesprincbound}(iv).\\
\begin{figure}[h!]
\begin{minipage}[c]{1\textwidth}
    \centering
    \makebox[0pt]{%
    \begin{tikzpicture}[scale=0.8, transform shape]
\fill[lightgray] (0.9,0)rectangle(1.5,0.3) (0.9,0.3)rectangle(1.2,0.6);
\draw (1.05,0.45) node{\footnotesize{$\hat v$}};
\draw (1.5,0) node[below] {\footnotesize{(i) $\bar\sigma:=\sigma^{\hat v,r}\in\overline{\mathcal D}_{K-2}(r,s)$}};
\draw (1.78,-0.4) node[below,white] {\footnotesize{$\bar\sigma:=\sigma^{w,r}\in\overline{\mathcal E}_{K-2}(r,s)$}};
\draw[step=0.3cm,color=black] (0,0) grid (3,2.1);
\end{tikzpicture}
\begin{tikzpicture}[scale=0.8, transform shape]
\fill[lightgray] (0.9,0)rectangle(1.5,0.3) (0.9,0.3)rectangle(1.2,0.6);
\draw (1.35,0.15) node{\footnotesize{$\hat v$}};
\draw (1.78,-0.4) node[below,white] {\footnotesize{$\bar\sigma:=\sigma^{v,r}\in\overline{\mathcal D}_{K-2}(r,s)$}};
\draw (1.5,0) node[below] {\footnotesize{(ii) $\bar\sigma:=\sigma^{\hat v,r}\in\overline{\mathcal E}_{K-2}(r,s)$}};
\draw[step=0.3cm,color=black] (0,0) grid (3,2.1);
\end{tikzpicture}
\begin{tikzpicture}[scale=0.8, transform shape]
\fill[lightgray] (1.2,0)rectangle(1.5,0.3) (0.9,0.3)rectangle(1.2,0.6);
\draw (1.05,0.45) node{\footnotesize{$\hat w$}};
\draw (1.5,0) node[below] {\footnotesize{(iii) $\bar\sigma:=\sigma^{\hat w,r}\in\bar R_{1,1}(r,s)$}};
\draw (1.78,-0.4) node[below,white] {\footnotesize{$\bar\sigma:=\sigma^{w,r}\in\bar R_{1,1}(r,s)$}};
\draw[step=0.3cm,color=black] (0,0) grid (3,2.1);
\end{tikzpicture}
\begin{tikzpicture}[scale=0.8, transform shape]
\fill[lightgray] (1.2,0)rectangle(1.5,0.3) (1.2,0.6)rectangle(1.5,0.9);
\draw  (1.35,0.15) node{\footnotesize{$\hat w$}};
\draw (1.5,0) node[below] {\footnotesize{(iv) $\bar\sigma:=\sigma^{\hat w,r}\in\bar R_{1,1}(r,s)$}};
\draw (1.78,-0.4) node[below,white] {\footnotesize{$\bar\sigma:=\sigma^{w,r}\in\bar R_{1,1}(r,s)$}};
\draw[step=0.3cm,color=black] (0,0) grid (3,2.1);
\end{tikzpicture}%
    }\par
     \end{minipage}
\caption{\label{thirdexamplesprincbound} Examples of $\sigma\in\overline{\mathcal D}_{K-3}(r,s)$ and $\bar\sigma\in\mathcal B(\overline{\mathcal D}_{K-3}(r,s))$. We color white the vertices with spin $r$ and gray those vertices whose spin is $s$.}
\end{figure}
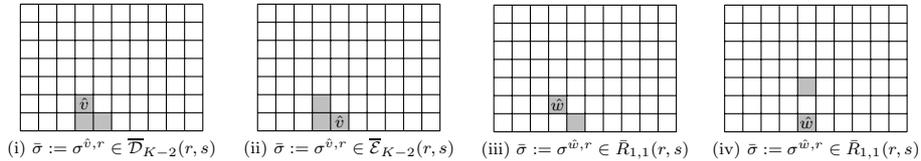\FloatBarrier

To conclude it is enough to see that \eqref{claimlastlemma} follows by Remark \ref{remprinbound} and \eqref{setKmeno2}. Indeed, both $\overline{\mathcal D}_{K-2}(r,s)$ and $\overline{\mathcal E}_{K-2}(r,s)$ are characterized by configurations in which there are two spins $s$ with three nearest-neighbor spins $r$ and by flipping from $s$ to $r$ one of these spins we obtain a configuration belonging to their principal boundary.
$\qed$\\

\appendix
\section{Appendix}\label{sectionappendix}
\subsection{Addition material for the Subsection \ref{proofappendixsub}}\label{appendixproof1}
In this subsection we give more details to the proof of items (b), (c) and (d) of Theorem \ref{mingatesNOcond}.\\

\textit{Proof of Theorem \ref{mingatesNOcond}.} Next we move to the proof that the sets depicted in (b) of Theorem \ref{mingatesNOcond} are minimal gates for the transition $\mathbf r\to\mathcal X^s\backslash\{\mathbf r\}$. Consider any $\omega\in\Omega_{\mathbf r,\mathcal X^s\backslash\{\mathbf r\}}^{opt}$ and let $\mathbf s\in\mathcal X^s\backslash\{\mathbf r\}$ be the first configuration visited by $\omega$ in $\mathcal X^s\backslash\{\mathbf r\}$. From Theorem \ref{mingatescond}(b) we have $\omega\cap\overline{\mathcal Q}(\mathbf r,\mathbf s)\neq\varnothing$ and $\omega\cap\widetilde{\mathcal Q}(\mathbf r,\mathbf s)\neq\varnothing$.  Thus, \[\omega\cap\biggl(\bigcup_{\mathbf t\in\mathcal X^s\backslash \{\mathbf r\}}\overline{\mathcal Q}(\mathbf r,\mathbf t)\biggr)\neq\varnothing\ \ \ \text{and}\ \ \ \omega\cap\biggl(\bigcup_{\mathbf t\in\mathcal X^s\backslash \{\mathbf r\}}\widetilde{\mathcal Q}(\mathbf r,\mathbf t)\biggr)\neq\varnothing\]
and (i) is verified. \\
Now consider $\eta\in\bigcup_{\mathbf t\in\mathcal X^s\backslash \{\mathbf r\}}\overline{\mathcal Q}(\mathbf r,\mathbf t)$. There exists $\mathbf s\in\mathcal X^s\backslash\{\mathbf r\}$ such that $\eta\in\overline{\mathcal Q}(\mathbf r,\mathbf s)$.  Let $\omega' \in\Omega_{\mathbf r,\mathcal X^s\backslash\{\mathbf r\}}^{opt}$ be the optimal path from $\mathbf r$ to $\mathbf s\in\mathcal X^s\backslash\{\mathbf r\}$ constructed in the proof of Theorem \ref{mingatescond}(b), such that $\omega'\cap(\bigcup_{\mathbf t\in\mathcal X^s\backslash \{\mathbf r\}}\overline{\mathcal Q}(\mathbf r,\mathbf t)\backslash\{\eta\})$\\$=\varnothing$ and $\omega'\cap\mathcal X^s\backslash\{\mathbf r,\mathbf s\}=\varnothing$. Hence, (ii) is verified for $\bigcup_{\mathbf t\in\mathcal X^s\backslash \{\mathbf r\}}\overline{\mathcal Q}(\mathbf r,\mathbf t)$. By the symmetry of the model, we can argue in the same way to prove that (ii) holds also for $\bigcup_{\mathbf t\in\mathcal X^s\backslash \{\mathbf r\}}\widetilde{\mathcal Q}(\mathbf r,\mathbf t)$.\\

Now we move to showing that the sets depicted in (c) of Theorem \ref{mingatesNOcond} are minimal gates for the transition $\mathbf r\to\mathcal X^s\backslash\{\mathbf r\}$. Consider any $\omega\in\Omega_{\mathbf r,\mathcal X^s\backslash\{\mathbf r\}}^{opt}$ and let $\mathbf s\in\mathcal X^s\backslash\{\mathbf r\}$ be the first configuration visited by $\omega$ in $\mathcal X^s\backslash\{\mathbf r\}$. From Theorem \ref{mingatescond}(c) we have $\omega\cap\overline{\mathscr H}(\mathbf r,\mathbf s)\neq\varnothing$ and $\omega\cap\widetilde{\mathscr H}(\mathbf r,\mathbf s)\neq\varnothing$. Thus, \[\omega\cap\biggl(\bigcup_{\mathbf t\in\mathcal X^s\backslash \{\mathbf r\}}\overline{\mathscr H}(\mathbf r,\mathbf t)\biggr)\neq\varnothing\ \ \ \text{and}\ \ \ \omega\cap\biggl(\bigcup_{\mathbf s\in\mathcal X^s\backslash \{\mathbf r\}}\widetilde{\mathscr H}(\mathbf r,\mathbf t)\biggr)\neq\varnothing\]
and (i) is verified. \\
Now consider $\eta\in\bigcup_{\mathbf t\in\mathcal X^s\backslash \{\mathbf r\}}\overline{\mathscr H}(\mathbf r,\mathbf t)$. There exists $\mathbf s\in\mathcal X^s\backslash\{\mathbf r\}$ such that $\eta\in\overline{\mathscr H}(\mathbf r,\mathbf s)$.  Let $\omega' \in\Omega_{\mathbf r,\mathcal X^s\backslash\{\mathbf r\}}^{opt}$ be the optimal path from $\mathbf r$ to $\mathbf s\in\mathcal X^s\backslash\{\mathbf r\}$ constructed in the proof of Theorem \ref{mingatescond}(c), such that $\omega'\cap(\bigcup_{\mathbf t\in\mathcal X^s\backslash \{\mathbf r\}}\overline{\mathscr H}(\mathbf r,\mathbf t)\backslash\{\eta\})$\\$=\varnothing$ and $\omega'\cap\mathcal X^s\backslash\{\mathbf r,\mathbf s\}=\varnothing$. Hence, (ii) is verified for $\bigcup_{\mathbf t\in\mathcal X^s\backslash \{\mathbf r\}}\overline{\mathscr H}(\mathbf r,\mathbf t)$. By the symmetry of the model, we can argue in the same way to prove that (ii) holds also for $\bigcup_{\mathbf t\in\mathcal X^s\backslash \{\mathbf r\}}\widetilde{\mathscr H}(\mathbf r,\mathbf t)$.\\
Finally, we prove that also the sets depicted in (d) of Theorem \ref{mingatesNOcond} are minimal gates for the transition $\mathbf r\to\mathcal X^s\backslash\{\mathbf r\}$. Consider any $\omega\in\Omega_{\mathbf r,\mathcal X^s\backslash\{\mathbf r\}}^{opt}$ and let $\mathbf s\in\mathcal X^s\backslash\{\mathbf r\}$ be the first configuration visited by $\omega$ in $\mathcal X^s\backslash\{\mathbf r\}$.
From Theorem \ref{mingatescond}(d) for any $j=2,\dots,L-3$ and any $h=1,\dots,K-1$, we have
\[\omega\cap\mathcal W_j^{(h)}(\mathbf r,\mathbf s)\neq\varnothing\]
and it follows that
\[\omega\cap\biggl(\bigcup_{\mathbf t\in\mathcal X^s\backslash \{\mathbf r\}}\mathcal W_j^{(h)}(\mathbf r,\mathbf t)\biggr)\neq\varnothing\] and (i) is verified. \\
For any $j=2,\dots,L-3$ and $h=1,\dots,K-1$, consider $\eta\in\bigcup_{\mathbf t\in\mathcal X^s\backslash \{\mathbf r\}} W_j^{(h)}(\mathbf r,\mathbf t)$, thus there exists $\mathbf s\in\mathcal X^s\backslash\{\mathbf r\}$ such that $\eta\in W_j^{(h)}(\mathbf r,\mathbf s)$.  Let $\omega' \in\Omega_{\mathbf r,\mathcal X^s\backslash\{\mathbf r\}}^{opt}$ be the optimal path from $\mathbf r$ to $\mathbf s\in\mathcal X^s\backslash\{\mathbf r\}$ constructed in the proof of Theorem \ref{mingatescond}(d) such that $\omega'\cap(\bigcup_{\mathbf t\in\mathcal X^s\backslash \{\mathbf r\}} W_j^{(h)}(\mathbf r,\mathbf t)\backslash\{\eta\})=\varnothing$ and $\omega'\cap\mathcal X^s\backslash\{\mathbf r,\mathbf s\}=\varnothing$. Thus, (ii) is verified and $\bigcup_{\mathbf t\in\mathcal X^s\backslash \{\mathbf r\}}\mathcal W_j^{(h)}(\mathbf r,\mathbf t)$ is a minimal gate for the transition from $\mathbf r$ to $\mathcal X^s\backslash\{\mathbf r\}$. $\qed$ 
\subsection{Addition material for the Subsection \ref{proofmingateappendix}}\label{appendixproof2}
In this subsection we give more details to the proof of items (b), (c) and (d) of Theorem \ref{mingatessingh0}.\\

\textit{Proof of Theorem \ref{mingatessingh0}.} Next we move to prove that the sets depicted in (b) are minimal gates for the transition from $\mathbf r$ to $\mathbf s$. First, we focus on the sets
\begin{align}\label{Qmingatefirstrs}
\bigcup_{\mathbf t\in\mathcal X^s\backslash\{\mathbf r\}} \overline{\mathcal Q}(\mathbf r,\mathbf t)\ \ \ \text{and}\ \bigcup_{\mathbf t\in\mathcal X^s\backslash\{\mathbf r\}} \widetilde{\mathcal Q}(\mathbf r,\mathbf t).
\end{align} Let $\omega\in\Omega_{\mathbf r,\mathbf s}^{opt}$ and let $\mathbf z\in\mathcal X^s\backslash\{\mathbf r\}$ be the first stable configuration visited by $\omega$ after $\mathbf r$. By Theorem \ref{mingatescond}(b), we get $\omega\cap \overline{\mathcal Q}(\mathbf r,\mathbf z)\neq\varnothing$ and $\omega\cap\widetilde{\mathcal Q}(\mathbf r,\mathbf z)\neq\varnothing$. Thus,
$$\omega\cap(\bigcup_{\mathbf t\in\mathcal X^s\backslash\{\mathbf r\}} \overline{\mathcal Q}(\mathbf r,\mathbf t))\neq\varnothing\ \ \ \text{and}\ \ \ \omega\cap(\bigcup_{\mathbf t\in\mathcal X^s\backslash\{\mathbf r\}} \widetilde{\mathcal Q}(\mathbf r,\mathbf t))\neq\varnothing$$
and (i) is satisfied. Let us now prove that for any $\eta\in\bigcup_{\mathbf t\in\mathcal X^s\backslash\{\mathbf r\}} \overline{\mathcal Q}(\mathbf r,\mathbf t)$ there exists $\omega'\in\Omega_{\mathbf r,\mathbf s}^{opt}$ that satisfies (ii). There exists $\mathbf z\in\mathcal X^s\backslash\{\mathbf r\}$ such that $\eta\in \overline{\mathcal Q}(\mathbf r,\mathbf z)$. Let $\omega^{(1)}:\mathbf r\to\mathbf z$ be the path constructed in the proof of Theorem \ref{mingatescond}(b) such that $\omega^{(1)}\cap\mathcal X^s\backslash\{\mathbf r,\mathbf z\}=\varnothing$. Let $\omega^{(2)}:\mathbf z\to\mathbf s$ be the reference path given in Definition \ref{remarkrepathNZ}. If $\mathbf z=\mathbf s$, then we define the optimal path $\omega'$ as the path $\omega^{(1)}$. Otherwise, we define $\omega'$ as the concatenation of the paths $\omega^{(1)}$ and $\omega^{(2)}$. Concerning $\bigcup_{\mathbf t\in\mathcal X^s\backslash\{\mathbf r\}} \widetilde{\mathcal Q}(\mathbf r,\mathbf t)$, we argue similarly to prove (ii) using the symmetry of the model on $\Lambda$. In both cases, (ii) is verified and the sets in \eqref{Qmingatefirstrs} are minimal gates for the transition from $\mathbf r$ to $\mathbf s$.\\
Second, let us focus on the sets 
\begin{align}\label{Qmingatesecondrs}
\bigcup_{\mathbf t\in\mathcal X^s\backslash\{\mathbf s\}} \overline{\mathcal Q}(\mathbf t,\mathbf s)\ \ \ \text{and}\ \bigcup_{\mathbf t\in\mathcal X^s\backslash\{\mathbf s\}} \widetilde{\mathcal Q}(\mathbf t,\mathbf s).
\end{align}
 Let $\omega\in\Omega_{\mathbf r,\mathbf s}^{opt}$ and let $\mathbf z\in\mathcal X^s\backslash\{\mathbf s\}$ be the last stable configuration visited by $\omega$ before hitting $\mathbf s$. By Theorem \ref{mingatescond}(b), we get $\omega\cap\overline{\mathcal Q}(\mathbf z,\mathbf s)\neq\varnothing$ and $\omega\cap\widetilde{\mathcal Q}(\mathbf z,\mathbf s)\neq\varnothing$. Thus,
$$\omega\cap(\bigcup_{\mathbf t\in\mathcal X^s\backslash\{\mathbf s\}} \overline{\mathcal Q}(\mathbf t,\mathbf s))\neq\varnothing\ \ \ \text{and}\ \ \ \omega\cap(\bigcup_{\mathbf t\in\mathcal X^s\backslash\{\mathbf s\}} \widetilde{\mathcal Q}(\mathbf t,\mathbf s))\neq\varnothing$$
and (i) is satisfied. Let us now prove that for any $\eta\in\bigcup_{\mathbf t\in\mathcal X^s\backslash\{\mathbf s\}} \overline{\mathcal Q}(\mathbf t,\mathbf s)$ there exists $\omega'\in\Omega_{\mathbf r,\mathbf s}^{opt}$ that satisfies (ii). There exists $\mathbf z\in\mathcal X^s\backslash\{\mathbf s\}$ such that $\eta\in\overline{\mathcal Q}(\mathbf z,\mathbf s)$. Let $\omega^{(1)}:\mathbf r\to\mathbf z$ be the reference path given in Definition \ref{remarkrepathNZ}. Let $\omega^{(2)}:\mathbf z\to\mathbf s$ be
the path constructed in the proof of Theorem \ref{mingatescond}(b) such that $\omega^{(2)}\cap\mathcal X^s\backslash\{\mathbf z,\mathbf s\}=\varnothing$. Thus, if $\mathbf z=\mathbf r$ we define the optimal path $\omega'$ as the path $\omega^{(2)}$. Otherwise, we define $\omega'$ as the concatenation of the paths $\omega^{(1)}$ and $\omega^{(2)}$. Concerning $\bigcup_{\mathbf t\in\mathcal X^s\backslash\{\mathbf s\}} \widetilde{\mathcal Q}(\mathbf t,\mathbf s)$,  we argue similarly to prove (ii) using the symmetry of the model on $\Lambda$. In both cases, (ii) is verified and the sets in \eqref{Qmingatesecondrs} are minimal gates for the transition from $\mathbf r$ to $\mathbf s$.\\
Similarly, we prove that the sets depicted in (c) are minimal gates for the transition from $\mathbf r$ to $\mathbf s$. First, we focus on the sets
\begin{align}\label{Hmingatefirstrs}
\bigcup_{\mathbf t\in\mathcal X^s\backslash\{\mathbf r\}} \overline{\mathscr H}(\mathbf r,\mathbf t)\ \ \ \text{and}\ \bigcup_{\mathbf t\in\mathcal X^s\backslash\{\mathbf r\}} \widetilde{\mathscr H}(\mathbf r,\mathbf t).
\end{align} Let $\omega\in\Omega_{\mathbf r,\mathbf s}^{opt}$ and let $\mathbf z\in\mathcal X^s\backslash\{\mathbf r\}$ be the first stable configuration visited by $\omega$ after $\mathbf r$. By Theorem \ref{mingatescond}(c), we get $\omega\cap \overline{\mathscr H}(\mathbf r,\mathbf z)\neq\varnothing$ and $\omega\cap\widetilde{\mathscr H}(\mathbf r,\mathbf z)\neq\varnothing$. Thus,
$$\omega\cap(\bigcup_{\mathbf t\in\mathcal X^s\backslash\{\mathbf r\}} \overline{\mathscr H}(\mathbf r,\mathbf t))\neq\varnothing\ \ \ \text{and}\ \ \ \omega\cap(\bigcup_{\mathbf t\in\mathcal X^s\backslash\{\mathbf r\}} \widetilde{\mathscr H}(\mathbf r,\mathbf t))\neq\varnothing$$
and (i) is satisfied. Let us now prove that for any $\eta\in\bigcup_{\mathbf t\in\mathcal X^s\backslash\{\mathbf r\}} \overline{\mathscr H}(\mathbf r,\mathbf t)$ there exists $\omega'\in\Omega_{\mathbf r,\mathbf s}^{opt}$ that satisfies (ii). There exists $\mathbf z\in\mathcal X^s\backslash\{\mathbf r\}$ such that $\eta\in \overline{\mathscr H}(\mathbf r,\mathbf z)$. Let $\omega^{(1)}:\mathbf r\to\mathbf z$ be the path constructed in the proof of Theorem \ref{mingatescond}(c) such that $\omega^{(1)}\cap\mathcal X^s\backslash\{\mathbf r,\mathbf z\}=\varnothing$. Let $\omega^{(2)}:\mathbf z\to\mathbf s$ be the reference path given in Definition \ref{remarkrepathNZ}. If $\mathbf z=\mathbf s$, then we define the optimal path $\omega'$ as the path $\omega^{(1)}$. Otherwise, we define $\omega'$ as the concatenation of the paths $\omega^{(1)}$ and $\omega^{(2)}$. Concerning $\bigcup_{\mathbf t\in\mathcal X^s\backslash\{\mathbf r\}} \widetilde{\mathscr H}(\mathbf r,\mathbf t)$, we argue similarly to prove (ii) using the symmetry of the model on $\Lambda$. In both cases, (ii) is verified and the sets in \eqref{Hmingatefirstrs} are minimal gates for the transition $\mathbf r\to\mathbf s$.\\
Second, let us focus on the sets 
\begin{align}\label{Hmingatesecondrs}
\bigcup_{\mathbf t\in\mathcal X^s\backslash\{\mathbf s\}} \overline{\mathscr H}(\mathbf t,\mathbf s)\ \ \ \text{and}\ \bigcup_{\mathbf t\in\mathcal X^s\backslash\{\mathbf s\}} \widetilde{\mathscr H}(\mathbf t,\mathbf s).
\end{align}
Let $\omega\in\Omega_{\mathbf r,\mathbf s}^{opt}$ and let $\mathbf z\in\mathcal X^s\backslash\{\mathbf s\}$ be the last stable configuration visited by $\omega$ before hitting $\mathbf s$. By Theorem \ref{mingatescond}(c), we get $\omega\cap\overline{\mathscr H}(\mathbf z,\mathbf s)\neq\varnothing$ and $\omega\cap\widetilde{\mathscr H}(\mathbf z,\mathbf s)\neq\varnothing$. Thus,
$$\omega\cap(\bigcup_{\mathbf t\in\mathcal X^s\backslash\{\mathbf s\}} \overline{\mathscr H}(\mathbf t,\mathbf s))\neq\varnothing\ \ \ \text{and}\ \ \ \omega\cap(\bigcup_{\mathbf t\in\mathcal X^s\backslash\{\mathbf s\}} \widetilde{\mathscr H}(\mathbf t,\mathbf s))\neq\varnothing$$
and (i) is satisfied. Let us now prove that for any $\eta\in\bigcup_{\mathbf t\in\mathcal X^s\backslash\{\mathbf s\}} \overline{\mathscr H}(\mathbf t,\mathbf s)$ there exists $\omega'\in\Omega_{\mathbf r,\mathbf s}^{opt}$ that satisfies (ii). There exists $\mathbf z\in\mathcal X^s\backslash\{\mathbf s\}$ such that $\eta\in\overline{\mathscr H}(\mathbf z,\mathbf s)$. Let $\omega^{(1)}:\mathbf r\to\mathbf z$ be the reference path given in Definition \ref{remarkrepathNZ}. Let $\omega^{(2)}:\mathbf z\to\mathbf s$ be
the path constructed in the proof of Theorem \ref{mingatescond}(c) such that $\omega^{(2)}\cap\mathcal X^s\backslash\{\mathbf z,\mathbf s\}=\varnothing$. If $\mathbf z=\mathbf r$, then we define the optimal path $\omega'$ as the path $\omega^{(2)}$. Otherwise, we define $\omega'$ as the concatenation of the paths $\omega^{(1)}$ and $\omega^{(2)}$. Concerning $\bigcup_{\mathbf t\in\mathcal X^s\backslash\{\mathbf s\}} \widetilde{\mathscr H}(\mathbf t,\mathbf s)$, we argue similarly to prove (ii) using the symmetry of the model on $\Lambda$. In both cases, (ii) is verified and the sets in \eqref{Hmingatesecondrs} are minimal gates for the transition $\mathbf r\to\mathbf s$.\\ 
Finally, we prove that the sets depicted in (d) are minimal gates for the transition from $\mathbf r$ to $\mathbf s$. First we consider the set 
\begin{align}\label{aligngatersfix}
\bigcup_{\mathbf t\in\mathcal X^s\backslash\{\mathbf r\}} \mathcal W_j^{(h)}(\mathbf r,\mathbf t),\ \text{for any $j=2,\dots,L-3$ and any $h=1,\dots,K-1$}.
\end{align}
Let $\omega\in\Omega_{\mathbf r,\mathbf s}^{opt}$ and let $\mathbf z\in\mathcal X^s\backslash\{\mathbf r\}$ be the first stable configuration visited by $\omega$ after $\mathbf r$. By Theorem \ref{mingatescond}(d), for any $j=2,\dots,L-3$ and any $h=1,\dots,K-1$, we get $\omega\cap\mathcal W_j^{(h)}(\mathbf r,\mathbf z)\neq\varnothing$. Thus, for any $j=2,\dots,L-3$ and any $h=1,\dots,K-1$,
$$\omega\cap(\bigcup_{\mathbf t\in\mathcal X^s\backslash\{\mathbf r\}} \mathcal W_j^{(h)}(\mathbf r,\mathbf t))\neq\varnothing$$
and (i) is verified. Let us now prove that for any $\eta\in\bigcup_{\mathbf t\in\mathcal X^s\backslash\{\mathbf r\}} \mathcal W_j^{(h)}(\mathbf r,\mathbf t)$ there exists $\omega'\in\Omega_{\mathbf r,\mathbf s}^{opt}$ that satisfies (ii). For any $j=2,\dots,L-3$ and any $h=1,\dots,K-1$, there exists $\mathbf z\in\mathcal X^s\backslash\{\mathbf r\}$ such that $\eta\in\mathcal W_j^{(h)}(\mathbf r,\mathbf z)$. Let $\omega^{(1)}:\mathbf r\to\mathbf z$ be the path constructed in the proof of Theorem \ref{mingatescond}(d) such that $\omega^{(1)}\cap\mathcal X^s\backslash\{\mathbf r,\mathbf z\}=\varnothing$. Let $\omega^{(2)}:\mathbf z\to\mathbf s$ be the reference path given in Definition \ref{remarkrepathNZ}. If $\mathbf z=\mathbf s$, then we define the optimal path $\omega'$ as the path $\omega^{(1)}$. Otherwise, we define $\omega'$ as the concatenation of the paths $\omega^{(1)}$ and $\omega^{(2)}$. Thus (ii) is satisfied the set \eqref{aligngatersfix} is a minimal gate for the transition $\mathbf r\to\mathbf s$.\\
Second, let us focus on the set
\begin{align}
\bigcup_{\mathbf t\in\mathcal X^s\backslash\{\mathbf s\}}  \mathcal W_j^{(h)}(\mathbf t,\mathbf s),\ \text{for any $j=2,\dots,L-3$ and any $h=1,\dots,K-1$}.
\end{align}
Let $\omega\in\Omega_{\mathbf r,\mathbf s}^{opt}$ and let $\mathbf z\in\mathcal X^s\backslash\{\mathbf s\}$ be the last stable configuration visited by $\omega$ before hitting $\mathbf s$. By Theorem \ref{mingatescond}(d), we get $\omega\cap \mathcal W_j^{(h)}(\mathbf z,\mathbf s)\neq\varnothing$ and $\omega\cap \mathcal W_j^{(h)}(\mathbf z,\mathbf s)\neq\varnothing$. Thus,
$$\omega\cap(\bigcup_{\mathbf t\in\mathcal X^s\backslash\{\mathbf s\}}  \mathcal W_j^{(h)}(\mathbf t,\mathbf s))\neq\varnothing\ \ \ \text{and}\ \ \ \omega\cap(\bigcup_{\mathbf t\in\mathcal X^s\backslash\{\mathbf s\}}  \mathcal W_j^{(h)}(\mathbf t,\mathbf s))\neq\varnothing$$
and (i) is satisfied. Let us now prove that for any $\eta\in\bigcup_{\mathbf t\in\mathcal X^s\backslash\{\mathbf s\}}  \mathcal W_j^{(h)}(\mathbf t,\mathbf s)$ there exists $\omega'\in\Omega_{\mathbf r,\mathbf s}^{opt}$ that satisfies (ii). There exists $\mathbf z\in\mathcal X^s\backslash\{\mathbf s\}$ such that $\eta\in \mathcal W_j^{(h)}(\mathbf z,\mathbf s)$. Let $\omega^{(1)}:\mathbf r\to\mathbf z$ be the reference path given in Definition \ref{remarkrepathNZ}. Let $\omega^{(2)}:\mathbf z\to\mathbf s$ be
the path constructed in the proof of Theorem \ref{mingatescond}(c) such that $\omega^{(2)}\cap\mathcal X^s\backslash\{\mathbf z,\mathbf s\}=\varnothing$. Thus, if $\mathbf z=\mathbf r$ we define the optimal path $\omega'$ as the path $\omega^{(2)}$. Otherwise, we define $\omega'$ as the concatenation of the paths $\omega^{(1)}$ and $\omega^{(2)}$. Thus (ii) is satisfied and the sets depicted in (d) are minimal gates for the transition $\mathbf r\to\mathbf s$.
$\qed$
\section{Table of notations}
\vspace{-3mm}
\renewcommand{\arraystretch}{1.2}
\begin{longtable}{p{1.5cm}p{12cm}}

\textbf{Symbol} & \textbf{Meaning} \\

$H$  & Hamiltonian function \\

$\mathcal X$ & State space \\

$\mathcal X^s$ & Set of global minima of $H$ \\

$\mathcal X^m$ & Set of metastable configurations \\

$S$ & Set of all possible spin values \\

$N_s(\sigma)$ & Number of spins $s$ in configuration $\sigma$ \\

$\mathcal X(r,s)$ & Set of configurations with all spins either $r$ or $s$ \\

$\Omega_{\sigma,\sigma'}$ & Set of all the paths from $\sigma$ to $\sigma'$ \\

$\Omega_{\sigma,\sigma'}^{opt}$ & Set of all the optimal paths from $\sigma$ to $\sigma'$ \\

$\Phi_\omega$ & $\max_{i=0,\dots,n} H(\omega_i)$, height of the path $\omega$ \\ 

$\Phi(\sigma,\sigma')$ & $\min_{\omega:\sigma \to \sigma'} \Phi_\omega$, communication height between $\sigma$ and $\sigma'$\\

$\tau^\sigma_{\mathcal A}$ & First hitting time of the subset $\mathcal A\subset\mathcal X$ starting from $\sigma$ \\

$\mathcal S(\sigma,\sigma')$ & Set of essential saddles for the transition from $\sigma$ to $\sigma'$ \\

$\mathcal W(\sigma,\sigma')$ & Gate for the transition from $\sigma$ to $\sigma'$ \\

$\mathcal W_\text{RES}(\mathbf s,\mathbf s')$ & Restricted-gate for the transition from $\mathbf s\in\mathcal X^s$ and $\mathbf s'\in\mathcal X^s\backslash\{\mathbf s\}$ following a path that does not intersect $\mathcal X^s\backslash\{\mathbf s,\mathbf s'\}$ \\

$\mathcal G(\sigma,\sigma')$ & Union of all minimal gates for the transition from $\sigma$ to $\sigma'$ \\

$\mathcal F(\mathbf s,\mathbf s')$ & Union of all minimal restricted-gate for the transition $\sigma\to\sigma'$ \\

$\mathcal C_{\mathcal A}(\sigma)$ & Initial cycle for the transition $\sigma\to\mathcal A$\\

$\mathscr C(\mathcal X)$ & Set of all the cycles and extended cycles of $\mathcal X$\\

$\mathcal M(\mathcal A)$ & Collection of maximal cycles and extended cycle that partitions $\mathcal A\subset\mathcal X$\\

$\mathcal B(\mathcal C)$ & Principal boundary of the cycle $\mathcal C$\\

$\partial^{np}\mathcal C$ & Non-principal boundary of the cycle $\mathcal C$\\

$T_{\mathcal A}(\sigma)$  & $\{\eta\in\mathcal X|\ \exists\omega\in\Omega_{\sigma,A}^{\text{vtj}}:\ \eta\in\omega\}$, tube of typical paths from $\sigma$ to $\mathcal A$\\

$\mathfrak T_{\mathcal A}(\sigma)$ & $\{\mathcal C\in\mathcal M(\mathcal C_{\mathcal A}^+(\sigma)\backslash\mathcal A)| \exists(\mathcal C_1,\dots,\mathcal C_n)\in J_{\mathcal C_{\mathcal A}(\sigma),\mathcal A}, \text{and}\ \exists j\in\{1,\dots,n\}: \mathcal C_j=\mathcal C\}$, tube of typical paths from $\sigma$ to $\mathcal A$\\

$\mathcal U_{\sigma'}(\sigma)$ & $\{\eta\in\mathcal X|\ \exists\omega\in\Omega_{\sigma,\sigma'}^{\text{vtj}}\ \text{s.t.}\ \omega\cap\mathcal X^s\backslash\{\sigma,\sigma'\}=\varnothing\ \text{and}\ \eta\in\omega\}$, restricted-tube of typical paths from $\sigma$ to $\sigma'$\\

$\mathscr U_{\sigma'}(\sigma)$ & $\{\mathcal C\in\mathcal M(\mathcal C_{\{\sigma'\}}^+(\sigma)\backslash\{\sigma'\})|\exists(\mathcal C_1,\dots,\mathcal C_m)\in J_{\mathcal C_{\sigma'}(\sigma),\{\sigma'\}} \text{s. t.} \bigcup_{i=1}^m \mathcal C_i\cap\mathcal X^s\backslash\{\sigma,\sigma'\}=\varnothing\ \text{and}\ \exists j\in\{1,\dots,n\}: \mathcal C_j=\mathcal C\}$, restricted-tube of typical paths from $\sigma$ to $\sigma'$\\

$C^s(\sigma)$ & Union of all the $s$-cluster in the configuration $\sigma$ for some $s\in S$ \\

$R(C^s(\sigma))$ & Smallest rectangle surrounding $C^s(\sigma)$ \\

$R_{\ell_1\times\ell_2}$ & Rectangle in $\mathbb R^2$ with sides of length $\ell_1$ and $\ell_2$ \\

$\bar R_{a,b}(r, s)$ & Set of configurations with all spins $r$, except those, which are $s$, in a rectangle $a\times b$ \\

$\tilde R_{a,b}(r, s)$ & Set of configurations with all spins $s$, except those, which are $r$, in a rectangle $a\times b$ \\

$\bar B_{a,b}^l(r, s)$ & Set of configurations with all spins $r$, except those, which are $s$, in a rectangle $a\times b$ with a bar $1\times l$ adjacent to one of the sides of length $b$, with $1\le l\le b-1$ \\

$\tilde B_{a,b}^l(r, s)$ & Set of configurations with all spins $s$, except those, which are $r$, in a rectangle $a\times b$ with a bar $1\times l$ adjacent to one of the sides of length $b$, with $1\le l\le b-1$ \\

$\overline{\mathscr H}(\mathbf r,\mathbf s)$ & $\bigcup_{i=1}^{K-3}\left[\bar B_{1,K}^i(r, s)\cup\bigcup_{h=i+1}^{K-2} \bar B_{1,K-1}^i(r, s)\right]$ \\

$\widetilde{\mathscr H}(\mathbf r,\mathbf s)$ & $\bigcup_{i=1}^{K-3}\left[\tilde B_{1,K}^i(r, s)\cup\bigcup_{h=i+1}^{K-2} \tilde B_{1,K-1}^i(r, s)\right]$ \\

$\overline{\mathcal Q}(\mathbf r,\mathbf s)$ & $\bar R_{2,K-1}(r, s)\cup\bar B_{1,K}^{K-2}(r, s)$ \\

$\widetilde{\mathcal Q}(\mathbf r,\mathbf s)$ & $\tilde R_{2,K-1}(r, s)\cup\tilde B_{1,K}^{K-2}(r, s)$ \\

$\overline{\mathscr P}(\mathbf r,\mathbf s)$ & $\bar B_{1,K}^{K-1}(r, s)$\\

$\widetilde{\mathscr P}(\mathbf r,\mathbf s)$ & $\tilde B_{1,K}^{K-1}(r, s)$ \\

$\mathcal{W}_j^{(h)}(\mathbf r,\mathbf s)$ & $\bar B_{j,K}^h(r, s)=\tilde B_{L-j-1,K}^{K-h}(r, s)\ \ \ \text{for}\ j=2,\dots,L-3$  \\

$\mathcal{W}_j(\mathbf r,\mathbf s)$ & $\bigcup_{h=1}^{K-1} \mathcal{W}_j^{(h)}(\mathbf r,\mathbf s)$ \\

$\overline{\mathcal K}(r,s)$ & $\{\sigma\in\mathcal X(r,s):\ H(\sigma)=2K+2+H(\mathbf r)$,  $\sigma$ has a $s$-cluster or more $s$-interacting clusters, $R(\text C^s(\sigma))=R_{2\times(K-1)}\}\cup\overline{\mathcal Q}(\mathbf r,\mathbf s)\cup\overline{\mathscr P}(\mathbf r,\mathbf s)$ \\

$\overline{\mathcal D}_1(r,s)$ & $\{\sigma\in\mathcal X(r,s):\ H(\sigma)=2K+H(\mathbf r)$,  $\sigma$ has either a $s$-cluster or more  $s$-interacting clusters such that $R(\text C^s(\sigma))=R_{2\times(K-2)}\}$  \\

$\overline{\mathcal E}_1(r,s)$ & $\{\sigma\in\mathcal X(r,s):\ H(\sigma)=2K+H(\mathbf r)$,  $\sigma$ has either a $s$-cluster or more $s$-interacting clusters such that $R(\text C^s(\sigma))=R_{1\times(K-1)}\}\cup\bar R_{1,K}(r, s)$. \\

$\overline{\mathcal D}_i(r,s)$ & $\{\sigma\in\mathcal X(r,s):\ H(\sigma)=2K-2i+2+H(\mathbf r)$, $\sigma$ has either a $s$-cluster or more $s$-interacting clusters such that $R(\text C^s(\sigma))=R_{2\times(K-(i+1))}\}$ \\

$\overline{\mathcal E}_i(r,s)$ & $\{\sigma\in\mathcal X(r,s):\ H(\sigma)=2K-2i+2+H(\mathbf r)$, $\sigma$ has either a $s$-cluster or more $s$-interacting clusters such that $R(\text C^s(\sigma))=R_{1\times(K-i)}\}$ \\

\end{longtable}
\bibliographystyle{abbrv}
\bibliography{mybib}
\end{document}